\pageno=1

\input amstex

\def\1{\hbox{\rm\rlap {1}\hskip .03in{\rm I}}} \def\Int{{\text {Int}}}
 
 \def\tr{{\text {tr}}}
\def\card{{\text {card}}}  \def\Z{\bold Z}

\hsize=12.5cm \vsize=19.0cm \parindent=0.5cm \parskip=0pt \baselineskip=12pt
\topskip=12pt \def\skipaline{\vskip 12pt plus 1pt} 

 \def\Hom{{\text {Hom}}} 
\def\vect{{\text {vect}}}
 \def\Eul{{\text {Eul}}}
 
\def\rank{{\text {rank}}}

\def\sign{{\text {sign}}}

\def\inc {{\text {in}}}

\def\Z {{\text {Z}}}
\def\char {{\text {char}}}
\def\red {{\text {red}}}

 \def\sin {{\text {sin}}}

\def\Ann {{\text {Ann}}}

\def\mod {{\text {mod}}}

 \def\dim  {{\text {dim}}}

\def\mod {{\text {mod}}\,}

\def\can {{\text {can}}}

\def\Coker {{\text {Coker}}}
 \def\Tors {{\text {Tors}}}

\def\ad {{\text {ad}}}

 \def\Ker {{\text {Ker}}} \def\Im {{\text {Im}}}
 \def\det {{\text {det}}}

 \centerline {  \bf  Surgery formula for torsions and Seiberg-Witten invariants of 3-manifolds}

\skipaline
\centerline {Vladimir  Turaev}

\skipaline
 \skipaline
 \skipaline
 \centerline {  \bf   Abstract}
\skipaline
 We give a  surgery formula for the torsions and Seiberg-Witten invariants associated with $Spin^c$-structures on 3-manifolds. 
We use the technique   of Reidemeister-type torsions and their refinements. 

   \skipaline
\skipaline
\centerline {\bf Contents}
\skipaline

\noindent  {\bf Introduction}

\noindent {\bf 1. Smooth Euler structures on 3-manifolds} 

\noindent {\bf 2.  Torsions of chain complexes}

\noindent  {\bf 3.  Homology orientations of 3-manifolds}

\noindent {\bf 4.  Combinatorial Euler structures on 3-manifolds}

\noindent {\bf 5.  The Torres formula for torsions}

\noindent {\bf 6.  Additivity of   torsions}

\noindent  {\bf 7.  The torsion $\tau$}

\noindent {\bf 8.   The Alexander-Conway function and derived invariants}

\noindent {\bf 9.   A surgery formula for  $\varphi$-torsions}

\noindent {\bf 10.  A surgery formula for the Alexander polynomial}

\noindent {\bf 11.   A surgery formula for  $\tau(M)$ in the case $b_1(M)\geq 1$} 

\noindent {\bf 12.   A surgery formula for the Seiberg-Witten invariants  of 3-manifolds}

\noindent {\bf   Appendix 1.  A surgery formula for  rational homology spheres}  

\noindent {\bf   Appendix 2. Computation of $\omega_L^M$}

\noindent {\bf   Appendix 3. Corrections and additions to [Tu5], [Tu6]}

\skipaline

\skipaline

\skipaline
\centerline {\bf    Introduction}

\skipaline

In 1976 the author introduced a combinatorial torsion-type invariant, $\tau$,  of compact  PL-manifolds of any dimension 
(see [Tu2]). For a compact 3-manifold $M$, the invariant $\tau(M)$ lies in the   group ring $\bold Z [H_1(M)]$ if  
$b_1(M)\geq 2$ and  lies in a certain extension of this  ring    if 
$b_1(M)=0,1$.  The study of    $\tau(M)$    was motivated by its connections with the Alexander-Fox
invariants of
$M$ including the   Alexander polynomial   $\Delta(M)$.

The definition of $\tau $ contains  an indeterminacy so that $\tau(M)$  is defined only up to
multiplication  by $\pm 1$ and by elements of $H_1(M)$.  In  [Tu3], [Tu4] the author introduced a refined version 
$\tau(M,e, \omega)$ of $\tau(M)$ depending on the choice of  a so-called Euler structure $e$ on $M$ and a homology orientation
 $\omega$ of 
$M$.  The invariant $\tau(M,e, \omega)$ has no indeterminacy  and   $\tau(M)=\pm H_1(M)\, \tau(M,e, \omega)$ for all $e,
\omega$.

The Seiberg-Witten invariant of a  closed oriented
3-manifold
$M$ with
$b_1(M)
\geq 1$ is a numerical function, $ SW=SW(M)$,  on the set of  $Spin^c$-structures on $M$, see for instance [FS],  [HL],
[Li], [MT], [Mo], [MOY], [OT].  This function also depends on the  choice  of  a homology orientation, $\omega$, 
  of 
$M$. For a $Spin^c$-structure $e$ on $M$, the integer  $SW(e,\omega)$ is the   algebraic    number of solutions, called
monopoles,   to a certain system of differential equations associated with $e$. This number   coincides with the  4-dimensional 
SW-invariant  of  the    $Spin^c$-structure $e \times 1$ on   $M\times S^1$. 

The invariants $ SW (M) $ and
$\tau  (M)$  turn out to be 
   equivalent (at least up to sign).  The first step in this direction was made by Meng and Taubes [MT] (see
also [FS]) who observed that $SW(M)$ determines the Alexander polynomial   $\Delta(M)$. The equivalence between $SW(M)$ 
and
$\tau(M)$ was established in [Tu5], [Tu6] 
where the SW-invariants
of $Spin^c$-structures on $M$ are identified with the   coefficients
in the expansion of $\tau(M)$ as an element of the group ring.
This involves
 an identification of   the Euler structures on $M$ with the $Spin^c$-structures  on $M$.

The definition of $\tau $ is based on the   methods of the theory of torsions, specifically, triangulations, 
chain complexes, etc.  The definition of  $SW $ is analytical. These definitions     are  not always suitable for explicit
computations.
The aim of this paper is to give     surgery formulas for $\tau $  and $SW$ suitable for   
computations.

We first give a surgery description of Euler structures  (=  $Spin^c$-structures) on  closed oriented 3-manifolds. To this end      
we introduce a notion of  a charge.
 A {\sl charge} 
on an oriented link
$L=L_1\cup ... \cup L_m$ in
$S^3$
is  an
$m$-tuple  $(k_1,...,k_m) \in \bold Z^m$ such that for all
$i=1,...,m$,  
$$k_i\equiv 1+ \sum_{j\neq i} lk (L_i, L_j) (\mod 2)  \eqno (0.a)$$
where $lk$ is the linking number in $S^3$. 
We   show that  a charge $k$ on $L$    determines  an Euler structure,  $e_k^M$, on any 3-manifold $M$  obtained by surgery 
on
$L$. 
Similarly, the orientation of $L$ induces  a homology orientation, $\omega_L^M$,  of $M$.

 Our main  formula computes
$\tau(M, e_k^M, \omega_L^M)$   in terms of the framing and  linking numbers of the components of $L$, and the
Alexander-Conway polynomials of $L$ and its sublinks.  
This  implies  a  surgery formula  for    $\pm  SW( e_k^M, \omega_L^M) $.
Thus, the algebraic number of monopoles can be computed (at least up to sign) in terms of classical  link invariants.

For the sake of introduction, we state  our surgery formula for  SW  in the case of    3-manifolds  with $b_1 \geq 2$ obtained  by
surgery on   algebraically split    links.
  Recall  that the
Alexander-Conway polynomial  $\nabla_L $ 
 of an oriented link
$L=L_1\cup ... \cup L_m\subset S^3$ with  $m\geq 2$  is a Laurent polynomial, i.e., an element of  
$\bold Z [t_1^{\pm 1},..., t_m^{\pm 1}]$.
If   $L$ is  algebraically split, i.e.,  if  $lk(L_i,L_j)=0$ for all $i\neq j$,  then   
$\nabla_L $ is   divisible by $ \prod_{i=1}^m (t_i^2-1)$ in $\bold Z [t_1^{\pm 1},..., t_m^{\pm 1}]$. 
 We have a  finite  expansion  
 $$\nabla_L / \prod_{i=1}^m (t_i^2-1)=\sum_{l=(l_1,...,l_m)\in \bold Z^m} z_l (L) \,\,  t_1^{l_1}... t_m^{l_m}$$
with $z_l (L)\in \bold Z$. 
 
Let    $M$ be a closed  oriented 3-manifold with $b_1(M)\geq 2$  obtained  by surgery on a      framed
oriented  algebraically split  link
$L =L_1\cup ... \cup L_m\subset S^3$.   
Let $f=(f_1,...,f_m)$ be the tuple of  the   framing numbers of $L_1,...,L_m$. 
Denote by
$J_0$  the  set of  all $j\in \{1,...,m\}$   such that  
$f_j=0$. For a   set $J\subset \{1,...,m\}$ we  denote the link 
 $\cup_{j\in  J} L_j$ by $L^J$.  Put   $\overline J=\{1,...,m\} \backslash J$ and $\vert J  \vert=\card (J)$.
Then for any charge  $k=(k_1,...,k_m)$ on $L$, 
$$  \pm  \, SW(e_k^M, \omega_L^M)   \eqno
(0.b)
$$
$$=  \sum_{J_0\subset J \subset \{1,...,m\}}\,\, (-1)^{\vert J  \vert}\, \prod_{j\in \overline J}  \sign (f_j)  \,\,  
  \sum_{l\in \bold Z^{J}, l\equiv -k (\mod 2f)} z_{l} (L^{J}). $$ 
 Here the sum goes over  
 all sets   $J\subset \{1,...,m\}$ containing $J_0$. The sign   $\sign (f_j)=\pm   1$   of $f_j$   is
well defined  since $f_j\neq 0$ for $j\in \overline J$. The formula
$l\in \bold Z^{ J}, l\equiv -k (\mod 2f)$   means that  
$l$ runs over all tuples of integers numerated by elements of $J$   such that $l_j\equiv -k_j (\mod 2f_j)$ for all
$j\in J$.   By
$\vert J  \vert\geq \vert J_0  \vert=b_1(M)\geq 2$, the algebraically split link $L^{J}$  has  $\geq
2$ components  so that 
$z_{l} (L^{J})$ is a well defined  integer.  Only a finite number of these integers  are non-zero and therefore
the sum on the right-hand side of (0.b) is finite. This sum  obviously depends only on $k (\mod 2f)$; the Euler structure 
 $e_k^M$ also depends  only on $k (\mod 2f)$.
The  precise sign   in  front of $SW$ in  (0.b) is 
unknown to the author. This sign is   the same for all charges
$k$ on
$L$.

For   a link  $L =L_1\cup ... \cup L_m$ which is not algebraically split,  the polynomial $\nabla_L $ can be   divided by $
\prod_{i=1}^m (t_i^2-1)$ in a certain quotient  of   $\bold Z [t_1^{\pm 1},..., t_m^{\pm 1}]$.  This leads to a  general  surgery
formula for  an arbitrary $L$.

We give also a surgery formula for the Alexander polynomial  $\Delta(M)$.  It is simpler than the surgery
formulas for $\tau $ and $SW$, and we establish it first.

The methods of this paper  yield similar surgery formulas for the invariants $\tau, \Delta, SW$ of  link exteriors in closed oriented
3-manifolds.  

 The invariants $\tau, \Delta, SW$ are closely related to the Casson-Walker-Lescop invariant of 3-manifolds. We shall
discuss these relations and other properties of   $\tau, \Delta, SW$ in another place.

 The   organisation of the paper should be clear from   Contents   above.

\skipaline
\skipaline \centerline {\bf 1.  Smooth Euler structures on 3-manifolds}

\skipaline \noindent {\bf   1.1.   Euler structures  on  manifolds.
}  We recall briefly the theory of (smooth) Euler structures on manifolds following [Tu4].  By a {\it relative vector
field} on a compact    (smooth) manifold 
$M$ we mean a {\it {nonsingular}}     tangent vector field on $M$  directed outside $M$ on
$\partial M$ (transversely to $\partial M$).  
Relative vector fields $u$ and $v$ on $M$ are called {\it homologous} if for  every connected component $M_0$ of $M$ and  a
point
$x\in M_0$ the restrictions of
$u$ and
$v$ to
$M_0\backslash
\{x\}$ are homotopic in the class of  nonsingular  vector fields on $M_0\backslash
\{x\}$   directed outside
$M_0$ on
$\partial M_0$.  
The homology class of a relative vector field  $u$ on $M$ is  called a (smooth)  {\it Euler
structure} on $M$ and denoted by $[u]$.  The set of Euler
structures on
$M$ is denoted by
$\vect (M)$. This set is non-void iff for each connected  component $M_0$ of $M$ we have   $\chi (M_0)=0$. This condition is 
satisfied for instance if
$M$ is a compact 3-manifold whose boundary consists of tori.

The group $H_1(M)$  acts  on  $\vect (M)$ as follows.  (This action as well as
the group operation in $H_1(M)$ will be written multiplicatively.)
  Let  $u$
be a relative vector field on $M$.  Let  $h\in H_1(M)$ be  represented  by an oriented simple closed curve $l\subset M\backslash
\partial M$.  Let $V$ be a regular neighborhood of $l$. We endow $V$ with
coordinates $(\theta, r)$, 
where $\theta \in \bold R/ 2 \pi \bold Z$ is the
coordinate along $l$ and $r$ is the radial coordinate on the discs 
transversal to $l$.
Choose a Riemannian metric on $V$ such that the tangent vector fields 
$\frac{\partial }{\partial  \theta}$ and $\frac{\partial }{\partial  r}$ are orthogonal everywhere and 
the Euclidean norm $\vert r\vert\geq 0$  measures the distance of a point of $V$  from $l$.    We
can  assume that $\vert r\vert_{\vert_{\partial  V}} = 1$, and after applying a homotopy to
$u$ we can assume that $u_{\vert_V} = - \frac{\partial }{\partial  \theta}$. 
Then $h [u]$ is represented by
the relative vector field on $M$ which is equal to $u$ on $M
\setminus V$ and is equal to $\cos (\vert r\vert \pi) \, \frac{\partial }{\partial  \theta} + \sin (\vert r\vert\pi) \,
\frac{\partial }{\partial  r}$ on $V$. 

If  $\vect (M) \neq \emptyset $ then the action of 
$H_1(M)$  on  $\vect (M)$ is free and transitive.
Thus,  for any 
  two relative vector fields   $u,v$ on $M$,   there is a unique
  $h\in H_1(M)$ such that $ [u]=h[v]$.  This  $h$ is denoted by $[u]/[v]$.  Clearly, $[v]/[u]=([u]/[v])^{-1}$.

For the torus $T=S^1\times S^1$, the set    $\vect (T)$ has a distinguished element defined as follows.
Consider a  non-singular vector field 
$w$ on $T$ whose flow lines are the circles
$x \times S^1$ with $  x\in S^1$. The class $[w]\in \vect (T)$ depends neither on the orientation of the circles  nor on the
choice of the splitting of $T$ as a product of two circles. This follows from the easy observation that $[w]$ is invariant under the
Dehn twists along  $S^1 \times x$ and $ x \times  S^1$.  In fact $[w]$  is the only element of
$\vect (T)$ invariant under all self-homeomorphisms of $T$.  
The formula
$h\mapsto h [w]$ defines a natural bijection
$H_1(T)
\to \vect (T)$.  We  call a non-singular vector field  on $T$  {\it canonical}
if it is homologous to $w$.

\skipaline \noindent {\bf   1.2. The Chern class.
}  Let $M$ be a compact 3-manifold whose boundary consists of tori.  For any $u\in \vect (M)$, the Chern class
$c(u) \in H_1(M)$ can be  defined as follows: consider  a  non-singular tangent vector field on $\partial M$
whose restrictions to all components are canonical;   the first obstruction to extending this vector field
to $M$ transversely to $u$ lies in $H^2(M,\partial M)$ and $c(u)\in H_1(M)$ is its  Poincar\'e  dual. We however shall use a
somewhat different definition of $c(u)$.

For a relative vector field $u$ on $M$, we define an opposite  vector field as follows. Deforming if necessary $u$  in a regular
neighborhood $\partial M\times [0,1] \subset M$
of $\partial M=\partial M\times 0$ we can assume that
$u$ is tangent to the lines $x\times [0,1]$ for $x\in \partial M$ and directed from $1$ to $0$.
Consider the vector field   on $M$ which is equal to $-u$ on the complement of  
$\partial M\times [0,1]$ in $M$ and is given by the formula $(x,t)\mapsto    \cos ( t\pi) u(x,t) +\sin (t\pi) w(x)$
on $\partial M\times [0,1]$ where $x\in \partial M, t\in [0,1]$ and $w$ is the canonical   vector field on the tori forming
$\partial M$. It is clear that this is a relative vector field on $M$. Its class in $\vect (M)$ depends only on
$e=[u]$.  We denote this class by $e^{-1}$.
It  follows from definitions that $(he)^{-1}= h^{-1} e^{-1}$ for any $h\in H_1(M), e\in \vect (M)$.

For $e\in \vect (M)$, set  $c(e)=e/e^{-1}$. We have the identities 
$$c(he)=(he)/(h^{-1}e^{-1}) =h^2(e/e^{-1})=h^2 c(e).$$ This formula has   useful applications. It implies that
the negation is an involution on $\vect (M)$. Indeed, we have
$$c(e^{-1})= c((e^{-1}/e) e)=(e^{-1}/e)^2 \,c(e)  = c(e)^{-2} c(e) = c(e)^{-1} \eqno (1.2.a)$$
and therefore
$$ (e^{-1})^{-1}= (c(e^{-1}))^{-1} e^{-1}= (c(e)^{-1})^{-1} e^{-1}=c(e) e^{-1}=e.$$
The formula $c(he)=h^2 c(e)$ shows   that the Chern class  map 
$e\mapsto c(e):\vect (M) \to H_1(M)$ is injective provided $H_1(M)$ has no 2-torsion.  
The formula $c(he)=h^2 c(e)$ implies also that the class $c(e) (\mod 2) \in H_1(M;\bold Z/2\bold Z)$
does not depend on
$e$. This class can be computed as follows.

\skipaline \noindent {\bf    1.3.  Lemma. } 
 {\sl  
Let $M$ be a compact   orientable 3-manifold whose boundary consists of tori.
Let $\Sigma\subset M$ be a compact  embedded surface  such that   $\partial \Sigma=\partial M\cap \Sigma$ 
  and  all components of $\partial \Sigma$ are non-contractible  on $\partial M$.
Then for any }  $e\in \vect (M)$  {\sl  we have }
$$c(e) \cdot  \Sigma
=b_0(\partial \Sigma) (\mod 2)$$
 {\sl  where $c(e) \cdot  \Sigma$ is the intersection number of    $c(e)$ with $\Sigma$  modulo 2 and $b_0(\partial \Sigma)$ is
the number of components of $\partial \Sigma$.}
\skipaline
 This lemma completely determines $c(e)  (\mod 2)$ since the  homology classes of
surfaces
$\Sigma$ as in the lemma generate the group $H_2(M, \partial M; \bold Z/2\bold Z)$ dual to $H_1(M; \bold Z/2\bold Z)$. We
give a proof of this lemma at the end of Sect.\  1.

\skipaline \noindent {\bf    1.4.   Example. }   
A   {\it  solid torus} is  the product $ S^1\times D^2$
where
$S^1$ is a  circle and $D^2$ is  a  closed 2-disc.
A    solid torus $\Z=S^1\times D^2$ endowed with a  generator,   $h_{\Z}$, of $ H_1(\Z)= \bold Z$ is said to be   
{\it directed}.   In  other words, $\Z $ is directed if its core circle $S^1\times  pt$ with $pt\in D^2$ is oriented.   Since
$H_1(\Z)$ has no torsion,  the Chern class map $ \vect (\Z) \to H_1(\Z)$ is injective.
Applying  Lemma 1.3  to  the meridional disc   $\Sigma=x\times D^2\subset \Z$ with $x\in S^1$ we obtain 
that the image of this map  
consists of the odd powers of the generator. Thus,  on a directed solid torus $\Z$ there is a unique Euler structure $e$ such
that
$c(e)=h_{\Z}^{-1}$. We call this $e \in \vect (\Z)$ the {\it distinguished Euler structure} on $\Z$ and denote it by $e_{\Z}$.

  \skipaline \noindent {\bf    1.5. Gluing of Euler structures. }
Let $M$ be a compact  3-manifold whose boundary consists of tori.
  Let  $T\subset M\backslash \partial M$ be a finite system of disjoint embedded 2-tori  splitting  $M$ into two 3-manifolds
$M_0$ and $M_1$. We define  now a gluing map $$\cup: \vect (M_0) \times \vect (M_1)\to \vect (M). \eqno (1.5.a)$$  Let
$u_i$ be a relative vector field on $M_i$ for $i=0,1$. We can  identify a regular neighborhood of $T$ in $M$ with
$T\times [0,1]$ so that $$T=T\times (1/2),   \,\,\,    \,\,\,  T\times [0,1/2] \subset M_0,  \,\,\,    \,\,\,  T\times [1/2,1]
\subset M_1.$$ We deform  $u_0$ so that it   is tangent to the lines $x\times [0,1/2]$ for $x\in T$  and is directed from $0$ to
$1/2$. Similarly, we 
 deform  $u_1$ so that it   is tangent to the lines $x\times [1/2,1]$ for $x\in T$  and is directed from $1$ to $1/2$.
 Let  $w$ be a non-singular  tangent vector field on $T$ whose restriction to each component of $T$ is canonical in the sense of
Sect.\   1.1. We define a  relative  vector field $u_0\cup u_1$ on $M$ as follows:  on $M_0\backslash (T\times
[0,1/2])$ it is equal to 
$u_0$, on $M_1\backslash (T\times [1/2,1])$ it is equal to  $u_1$,
at any point  $(x,t)\in T\times [0,1/2]$ it equals  $\cos (t\pi) u_0(x,t)+\sin (t\pi) w(x)$, 
at any point  $(x,t)\in T\times [1/2,1]$ it equals  $-\cos (t\pi) u_1(x,t)+\sin (t\pi) w(x)$.
The class $[u_0\cup u_1]\in \vect (M)$ depends only on $[u_0], [u_1]$. 
This yields a well-defined map (1.5.a). 

It follows  from definitions that
for any $  h_i\in H_1(M_i), e_i\in \vect (M_i)$ with $i=0,1$, 
$$ h_0 e_0 \cup  h_1e_1 = (\inc_0(h_0) \inc_1(h_1))\, (e_0\cup e_1)\eqno (1.5.b)$$
where $\inc_i$ is the inclusion homomorphism $H_1(M_i) \to H_1(M)$.  As an exercise, the reader may check that
$$(e_0\cup e_1)^{-1}= e_0^{-1} \cup e_1^{-1} $$
(hint: use that $w$ is homotopic to $-w$).
This implies that
$$c(e_0\cup e_1)=(e_0\cup e_1)/ (e_0\cup e_1)^{-1}
$$
$$=
(c(e_0) e_0^{-1} \cup
c(e_1) e_1^{-1})/(e_0^{-1} \cup e_1^{-1})= \inc_0(c(e_0)) \,\inc_1( c (e_1)).$$

 We shall use the gluing of Euler  structures in the following setting.   Consider a 3-manifold  $M$ obtained 
from a
compact 3-manifold $E$ by gluing $n$
directed solid tori
$\Z_1,...,\Z_n$ along $n$ tori in $\partial E$. Gluing to   $e\in \vect (M)$ the distinguished Euler structures on these solid tori we
obtain an Euler structure, $e^M$,  on $M$:
$$e^M=e\cup e_{ \Z_1} \cup...\cup e_{ \Z_n}\in \vect  (M).$$
 Clearly, $$c(e^M)= \inc (c(e)) \prod_{j=1}^n
h^{-1}_j \eqno (1.5.c)$$  where $\inc:H_1(E)\to H_1(M)$ is the inclusion homomorphism and $h_1,...,h_n\in H_1(M)$ are the
homology classes  of the oriented core circles of $\Z_1,...,\Z_n$. We have
$(e^M)^{-1}= \prod_{j=1}^n
h_j (e^{-1})^M$ and    $(he)^M= \inc (h)\, e^M$ for any $h\in H_1(E)$.

\skipaline \noindent {\bf 1.6.    Charges on  links and Euler structures on link exteriors.} The Euler structures on link exteriors
can be described in   terms  of   charges.  Let $L$
be  an oriented link
  in an oriented  3-dimensional integral homology sphere $N$.  A {\sl charge} on $L$ is an integer-valued function $k$ on the set
of components of $L$ such that for every component $\ell$ of $L$ we have $k(\ell)
\equiv 1+   lk (\ell, L\backslash \ell) (\mod 2)$ where $lk$ denotes the linking number   in $N$.

Consider the exterior  $E$ of $L$ in $N$.  There is a canonical bijection between the set $\vect (E)$ and the set of charges on
$L$. It is defined as follows.  For a component $\ell$ of $L$ denote by  $t_\ell\in H_1(E)$ the homology class of a meridian of
$\ell$.
 (By a meridian of an oriented knot we always mean  a  {\sl canonically oriented} meridian whose  linking number with the knot
equals
$+1$). 
For $e\in \vect (E) $, the Chern class $c(e)\in H_1(E)$ can be  uniquely expanded as $\prod_\ell  \,t_\ell^{k (\ell)}$  
  where $\ell$ runs over the components of $L$ and 
$k(\ell)\in \bold Z$.  We claim that (i) the function $\ell \mapsto k(\ell)$ is a charge on $L$ and (ii) the resulting map from
$\vect (E)$ to the set of charges on $L$ is bijective. 
To prove (i),  consider  a component 
$\ell$ of $ L$ and its   Seifert
surface
$\Sigma_\ell\subset N$  intersecting   $L\backslash \ell$ transversely.  Then the surface 
$\Sigma=\Sigma_\ell\cap E\subset E$ satisfies the conditions of Lemma 1.3. Note that 
$t_\ell \cdot  \Sigma =1 (\mod 2)$ and $t_{\ell'} \cdot  \Sigma = 0$ for any component $\ell'$ of $L\backslash \ell$.   This
and Lemma 1.3  imply that $$k (\ell) \equiv c(e) \cdot  \Sigma \equiv b_0(\partial \Sigma)\equiv 1+ lk (\ell, L\backslash \ell) 
(\mod 2).$$ Claim (ii) is essentially obvious: the injectivity follows from the injectivity  of the Chern class map $c: \vect (E) \to
H_1(E)$ and the surjectivity follows from the identity
$c(he) =h^2 c(e)$ where $e\in \vect (E)$, $h\in H_1(E)$.

For a charge $k$ on $L$, denote by $e_k$ the unique Euler structure on $E$ such that $c(e_k)=
\prod_\ell  \,t_\ell^{k (\ell)}$.  Formula (1.2.a) implies that
$(e_k)^{-1}=e_{-k}$ where $-k$ is the charge on $L$ defined by $(-k)(\ell)=-k(\ell)$. For any    charges   $k,k'$  on $L$, 
$$e_{k'}=\prod_\ell  \,t_\ell^{(k' (\ell) -k(\ell))/2}\, e_k.$$

 If  
$L=L_1\cup ...\cup L_m$ is ordered then a charge on $L$ is just an
$m$-tuple  $(k_1,...,k_m) \in \bold Z^m$ satisfying formula (0.a) for all
$i=1,...,m$.
For a set $I\subset \{1,...,m\}$, we   denote by $L^I$ the sublink  $\{L_i\}_{i\in I}$ of $L$.
A charge $k$ on $L$ induces a charge $k^I$ on $L^I$ by 
$$k^I(\ell)=k(\ell)-  lk (\ell, L\backslash L^{I})=k(\ell)-  lk (\ell, L^{\overline I})$$
where $\ell$ runs over  the components of $L^I$ and    $\overline I=\{1,...,m\} \backslash I$. 
The corresponding map $\vect (E) \to \vect (E^I)$ (where $E^I$ is the exterior of $L^I$) is described as follows.
We can  obtain $E^I$ from $E$ by gluing regular neighborhoods of the components of $L^{\overline I}$.    These regular
neighborhoods are directed solid tori: the orientation of their  cores   is determined by the   orientation of $L$.  Gluing to
$e_k\in \vect (E)$ the distinguished Euler structures on these  directed solid tori we obtain   $e_{k^I}\in \vect (E^I)$,  so that
in the notation of Sect.\ 1.5 we have 
$  e_{k^I}= (e_k)^{E^I} $.
This follows from   (1.5.c) applied to the inclusion $E\subset E^I$. (Here $h_1,..., h_n$ are the homology classes of
  the components of $L^{\overline I}$  in
$E^I$.)

\skipaline \noindent {\bf 1.7.  Surgery presentation of Euler structures.}  Let $M$ be a closed   3-manifold
obtained by Dehn surgery along an oriented  link $L=L_1\cup ...\cup L_m$ in an oriented 
 3-dimensional integral homology sphere  $N$.    This means that $M$ is obtained by
gluing    $m$ directed solid tori $ \Z_1,..., \Z_m$ to the exterior $E$ of $L$.  
Consider a charge $k$ on $L$. 
Gluing to $e_k\in \vect (E)$ the   distinguished Euler structures on $ \Z_1,...\Z_m$
we obtain  an Euler structure $e_k^M= (e_k)^M$ on $M$.
If  $k,k'\in \bold Z^{m}$ are two charges on $L$ then $k'-k=2h$ with $h\in \bold Z^m=H_1(E)$ and
$ e_{k'} =h   e_k,  e_{k'}^M =\inc (h) e_k^M$ where $\inc:H_1(E)\to H_1(M)$ is the inclusion homomorphism. 
This implies that  $ e_{k'}^M=e_k^M$ if and only if $h=(k'-k)/2\in \Ker (\inc)$. Since  $\inc:H_1(E)\to H_1(M)$  is surjective we
obtain  that  all  Euler structures on
$M$ have the form $e_k^M$ for a charge $k$ on $M$. 

For an  integral surgery on $L$ determined by a framing of $L$, we shall always use the following  structure of a  directed
solid tori  on   $ \Z_1,...,\Z_m$.  Let
$U_1,...,U_m$ be disjoint closed regular neighborhoods of $L_1,...,L_m$ so that
$E=N\backslash
\cup_{i=1}^m \Int(U_i)$. Let
$ \Z_i$ be the solid tori  glued to $\partial  U_i\subset  \partial E$ to form $M$ where $i=1,...,m$. 
A meridian of $L_i$ lying in $\partial  U_i$ is a core of $\Z_i$. The canonical orientation of this meridian induced by the
orientation of $L_i$ makes $\Z_i$ directed.

As an exercise, the reader may check that $(e_k^M)^{-1}=e_{2-k}^M$ where
$2-k$ is the charge on $L=L_1\cup ...\cup L_m$ defined by   $(2-k)_i=2-k_i$ for $i=1,...,m$.

\skipaline \noindent {\bf 1.8.  Proof of Lemma 1.3.}  In the case $\partial M=\emptyset$ Lemma 1.3 is essentially obvious. 
In this case $\partial \Sigma=\emptyset$ and the lemma amounts to saying that $c(e) (\mod 2)$ is a trivial homology class. 
Since $M$ is parallelizable, there is a non-singular vector field $u$ on $M$ homotopic to  $-u$.
Clearly, $ c([u])=1$ and the claim follows from the equality $c(e)\equiv c([u]) (\mod 2)$.

We
now prove Lemma 1.3 in the case  where
$M=\Z=S^1\times D^2$ is a   solid torus. Let $h $ be a generator of $H_1(\Z)$.  It is an
elementary topological exercise to check that $b_0(\partial \Sigma) \equiv h\cdot \Sigma (\mod 2)$. Therefore it remains
 to prove that for any  relative vector field $u$ on $\Z$ the class $c([u])$ is an odd power of $h$. 
We define two residues $n_+(u), n_-(u) \in \bold Z/2\bold Z$ as follows. 
Choose a trivialization  $v_1$ (resp.\ $v_2,v_3)$ of the
tangent bundle
   of $S^1$ (resp.\ of $D^2$) inducing the given orientation of $S^1$. This induces a trivialization  
$(v_1,v_2,v_3)$ of the tangent bundle  of    $  \Z$. Restricting $u$
to the disc $x \times D^2 \subset \Z$ with $x\in S^1$ and using the trivialization $(v_1,v_2,v_3)$  we obtain a map
$f_u$ from $D^2$ to the unit 2-sphere $S^2\subset  \bold R^3$  such that   $f_u(\partial D^2)$ does not meet the 
poles  
$(\pm 1,0,0)$ of
$S^2$. Let  $n_+(u)$ and $n_-(u)$ be the degrees mod 2 of $f_u$ with respect to these poles.   They depend only
on $[u]\in \vect (\Z)$.  It is clear that
$n_+(u)-n_-(u)=1$.  A simple  calculation shows that $n_{\pm} (h^k u) =n_{\pm }(u) +k (\mod 2)$ for any $k\in \bold Z$. The
restriction of
$v_1$   to
$\partial \Z$ is canonical and  can therefore be used   to describe  a vector field representing $[u]^{-1}$. On the regular
neighborhood of
$\partial
\Z$  used in the definition of  $[u]^{-1}$, this vector field never  belongs to  $-{\bold R}_+ v_1$. Therefore 
$n_-(u^{-1})= n_+ (u)=n_-(u)+1$.  Hence $[u]^{-1}= h ^{k} [u]$ with odd $k$ and  $c([u])=[u]/[u]^{-1}= h ^{-k}$ is
an odd power of $h$. 

Consider now the general case of Lemma 1.3. The intersection of $\Sigma$ with a component $T$ of $\partial M$
is a system of disjoint simple closed curves homotopic to each other.  We glue to $M$ a   solid torus  $ \Z_T$ along
$T$ so that the curves  forming $\Sigma \cap T$ bound disjoint embedded discs in $\Z_T$.
Denote the union of these discs  by $\Sigma_T$.
If $\Sigma\cap T=\emptyset$ then we glue a solid torus $\Z_T$ to $M$ along $T$ in an arbitrary way and set
$\Sigma_T=\emptyset$. Applying this procedure for   all components  $T$ of $\partial M$ we obtain a closed 3-manifold $\hat M$
and  a closed embedded surface $\hat \Sigma =\Sigma \cup  (\cup_T \Sigma_T)\subset \hat M$. 
Let $\hat e$ be the Euler structure on $\hat M$ obtained by gluing $e$ with  arbitrary Euler structures $\{e_T\in \vect
(\Z_T)\}_T$.  We have 
$c(\hat e)=c(e) \cdot \prod_T c (e_T)$. Therefore 
$$  c(\hat e) \cdot  \hat \Sigma=  c(e) \cdot \Sigma 
+\sum_T c (e_T) \cdot \Sigma_T \,(\mod 2). \eqno (1.8.a)$$
By the results of the previous paragraphs, 
$c(\hat e) \cdot  \hat \Sigma= 0 (\mod 2)$ and 
$$c (e_T) \cdot \Sigma_T =b_0 (\partial \Sigma_T)=b_0(\partial
\Sigma\cap T) \,(\mod 2).$$ Substituting this in (1.8.a), we obtain 
$     c(e) \cdot \Sigma 
= b_0(\partial \Sigma) (\mod 2)$.

\skipaline \noindent {\bf 1.9.  Remark.} The Euler structures on a closed oriented 3-manifold $M$ bijectively correspond to $Spin^c$-structures on the tangent vector
bundle $TM\to M$  of $M$.  The $Spin^c$-structure corresponding  to the class $[u]$ of  a non-singular  vector field $u$ on 
$M$ is defined as follows. 
  The  2-dimensional  quotient  vector bundle
$ TM/\bold  R u$  is oriented  and therefore  has the structure group $U(1)$. This   reduces the 
structure group of  $TM=(TM/\bold  R u) \oplus \bold  R u $ to $U(1)=U(1) \oplus (1) \subset U(2)=Spin^c(3)$.
Similarly, for  a compact oriented 3-manifold $M$ whose boundary
consists of  tori there is a   bijection of   $\vect (M)$   on the set of   relative $Spin^c$-structures on $M$,
see  [Tu6].

\skipaline
\skipaline \centerline {\bf  2.  Torsions of chain complexes}

  \skipaline \noindent {\bf   2.1.  Basic definitions. } Having two bases $c, c'$   of 
a finite dimensional vector space, 
we can expand  the vectors of $c$ as linear combinations of  the 
vectors of $c'$. The determinant of
the resulting square matrix is denoted by $[c/c']$. The bases $c,c'$ are   {\it equivalent} if  $[c/c']=1$.

 Let
$C=(C_m\to C_{m-1} \to  ... Ê\to C_0)$ be a finite dimensional chain complex  of length $m$  over a field
$F$.  Suppose that for all $i=0,1,...,m$ both $C_i$ and   $H_i(C)$ are based, i.e., have a distinguished basis. 
In this situation one defines the torsion $\tau(C)$ as follows (cf. [Mi1], [Mi2]).
Let $c_i$ be the given basis in $C_i$ and  let $h_i$ be a sequence of vectors in  
$\Ker (\partial_{i-1}: C_i \to C_{i-1})$ whose images under the projection 
$\Ker \, \partial_{i-1}\to H_i(C)$ form the given basis in $H_i(C)$. Let $b_i$ be a sequence of vectors in $C_i$ such that
$\partial_{i-1} (b_i)$ is a basis in $\Im\,  \partial_{i-1}$.  Then for every $i$,  the sequence
$\partial_i (b_{i+1}) h_i b_i$ is a basis in $C_i$.
Set
$$  \tau (C)=   \prod_{i=0}^m  \, [\partial_i (b_{i+1}) h_i b_i/c_i] ^{(-1)^{i+1}} \in F.$$
   The torsion  $  \tau (C)$ depends  on the equivalence classes of the given bases in $ C_i, H_i(C) $ and does not depend on
the choice of
$h_i, b_i$.  If $C$ is acyclic then  the definition of $\tau(C)$ simplifies since there is no need to fix a basis in $H_*(C)$.

We shall use  a \lq\lq sign-refined"   torsion   $\check \tau $ introduced in [Tu3].
Set
$$\alpha_i(C)= \dim\, C_i +\dim\, C_{i-1}+... + \dim\, C_0 \,(\mod 2) \in \bold Z/2\bold Z,$$
$$\beta_i(C)= \dim\, H_i (C)+\dim\, H_{i-1}(C) +... + \dim\, H_0(C) \,(\mod 2)\in \bold Z/2\bold Z,$$
$$ N(C)=\sum_{i=0}^m \alpha_i(C)\, \beta_i(C) \in \bold Z/2\bold Z.$$
Set   $\check \tau (C)=(-1)^{N(C)} \tau(C)  \in F$. 
If $C$ is acyclic then     $\check \tau (C)=\tau(C)$.

  \skipaline \noindent {\bf   2.2.  Relative torsions. }
 Let $C=(C_m\to ...Ê\to C_0)$ be a finite dimensional
chain complex over a field
$F$ and  let $C'=(C'_m\to ...Ê\to C'_0)$ be a chain subcomplex of $C$.   
Suppose that the homology of   $C,C'$ and $C''=C/C'$ are based. (The complexes themselves
are not assumed to be based).  We define   a {\it relative torsion} $ \tau ( C'\subset C)$   as follows.
Denote by
$\Cal H$   the
homology sequence of the pair $(C,C')$:
$$\Cal H= (H_m(C')\to H_m(C) \to H_m(C'')\to...\to H_0(C) \to H_0(C'')).  $$
Clearly, $\Cal H$ is a based acyclic chain complex over $F$. 
Set 
$$  \tau ( C'\subset C)=
(-1)^{\theta  (C,C')} \tau (\Cal H)\in F$$
where
$$\theta (C,C')=\sum_{i=0}^m [ (\beta_i(C) +1)( \beta_i(C')+\beta_i(C'')) + \beta_{i-1}(C') \beta_i(C'')] \in \bold Z/2\bold Z.$$
The relative torsion $\tau ( C'\subset C)$ can be expressed in terms of the torsions of $C,C',C''$.
Namely, provide 
for all
$i=0,1,...,m$ the vector spaces $C_i, C'_i, C''_i= C_i/C'_i$  with compatible bases  $c_i,c'_i, c''_i$ where the compatibility means
that   $c_i$ is equivalent to the basis $c'_i c''_i$ obtained as a juxtaposition      of   $c'_i$
with  a lift of  $c''_i$ to $C_i$. Then (see [Tu3], Lemma 3.4.2)
$$ \tau ( C'\subset C)= (-1)^{ \nu (C,C')} \frac { \check \tau (C)}{  \check \tau (C')\,\check \tau (C'') }  \eqno (2.2.a)$$
where $$\nu(C,C') =\sum_{i=0}^m \alpha_i(C'') \,\alpha_{i-1} (C')\in \bold Z/2\bold Z. \eqno (2.2.b)$$
If   $C, C'$ are acyclic,  then $\tau ( C'\subset C)=1$
and    (2.2.a) yields 
$$  \tau (C)=  (-1)^{ \nu (C,C')} \tau (C')\,  \tau (C''). \eqno (2.2.c)$$

\skipaline \noindent {\bf    2.3.  Lemma. }  
 {\sl  
 Let $C^1\subset C^2\subset C^3 $ be   finite dimensional
chain complexes of length $m$ over a field
$F$.  Suppose that the homology
of the chain complexes 
$C^1, C^2,C^3, C^2/C^1, C^3/C^1, C^3/C^2$ are    based. 
Then}
$$\tau( C^1\subset C^3 )
=  \frac
{\tau( C^1\subset C^2)\,
\tau(C^2\subset C^3)} { \tau(C^2/C^1\subset C^3/C^1)}.\eqno (2.3.a) $$

\skipaline {\sl Proof.} Choose in each   $C^3_i$ a basis whose first $\dim\, C^1_i$ vectors lie in $C^1_i$ and the next  $\dim\,
C^2_i -\dim\, C^1_i$ vectors lie in $C^2_i$. This basis determines bases in    $C^1_i, C^2_i , C^2_i/C^1_i, C^3_i/C^1_i,
C^3_i/C^2_i$ in the obvious way.   By (2.2.a), for any  $1\leq p <q \leq 3$,
$$\tau(C^p\subset C^q)= (-1)^{ \nu (C^q,C^p)} \frac
{\check \tau( C^q )} {  \check \tau( C^p ) \, \check \tau( C^q/C^p ) },$$
$$ \tau( C^2/C^1\subset C^3/C^1 )=
(-1)^{ \nu (C^3/C^1 ,  C^2/C^1 )} \frac
{\check \tau( C^3/C^1)} {  \check \tau( C^2/C^1)\,  \check \tau( C^3/C^2 ) }.$$
Substituting these expressions in (2.3.a) we obtain that (2.3.a) is equivalent to 
$${\nu (C^2,C^1 )+\nu (C^3,C^1) +\nu
(C^3 ,  C^2  )+
\nu (C^3/C^1,  C^2/C^1 )}=0.$$
This formula follows from (2.2.b) and the  
equalities
$\alpha_i(C^q/C^p)=\alpha_i(C^q) -\alpha_i(C^p)$ for  all $1\leq p <q \leq 3$. This implies (2.3.a).

\skipaline \skipaline \centerline {\bf  3.  Homology orientations of 3-manifolds}

\skipaline \noindent {\bf  3.1.   Homology orientations
}    A {\it homology orientation} $\omega$ of a  finite CW-pair $(X,Y)$ is an orientation of the vector space
$H_*(X,Y;\bold R)
=\oplus_{i\geq 0} H_i(X,Y;\bold R)$.  We   denote by $-\omega$ the opposite orientation.

The exterior, $E$,  of an oriented link $L=L_1\cup ...\cup L_m$ in an oriented 
 3-dimensional integral homology sphere  $N$ has  a 
canonical  homology orientation $\omega_L$  determined by the basis 
$([pt], t_1,...,t_m, g_1,...,g_{m-1})$  where $[pt]$ is the homology
 class of a point in $E$,  $t_1,...,t_m$ are  the meridional generators of $H_1(E;\bold R)$, and $g_1,...,g_{m-1}$
are the generators of
$H_2(E;\bold R)$ represented by   oriented boundaries of regular neighborhods of $L_1,...,L_{m-1}$, respectively. 
(We use the  \lq\lq outward vector first"  conveniton for the orientation of the boundary. The orientation in the regular
neighborhoods is induced by the one in $N$).  Note that $\omega_L$ does not depend on the numeration of the components of
$L$ and changes sign under inversion of the orientation of a component of $L$. 

\skipaline \noindent {\bf  3.2.   Induced homology orientations.}  Let $M$ be a   3-manifold
obtained by   gluing  $m$  directed solid
tori  $\Z_1=S^1\times D_1, ..., \Z_m=S^1\times D_m$ to a compact    3-manifold $E$. The aim of this section is to  
show that a homology orientation, $\omega$, of $E$ induces in a natural way a homology orientation, $\omega^M$, of $M$.

Fix an orientation of the 2-discs
$D_1, ... ,   D_m$ and   provide each   $\Z_i=S^1\times D_i$ with orientation obtained as the product of the given
orientation of its core and the fixed orientation in $D_i$.   Clearly, $$ H_2(\Z_i,\partial \Z_i;\bold R)=\bold R\, [D_i, \partial
D_i] ,
\,\,\,\,\, H_3(\Z_i,\partial \Z_i;\bold R)=\bold R\, [\Z_i, \partial \Z_i].
$$
 By excision, 
the vector
space
$H_*(M,E;\bold R)=\bold R^{2m}$ has a basis 
$$[D_1, \partial D_1],...,[D_m, \partial D_m], [\Z_1, \partial \Z_1],...,[\Z_m, \partial \Z_m].$$
Denote the homology orientation of $(M,E)$ determined by this basis by $\omega_{(M,E)}$.
 It is easy to check that $\omega_{(M,E)}$  depends neither on the choice of orientations in  $D_1, ...,   D_m$ nor on the
numeration of   $\Z_1,...,
\Z_m$. 

 There is a unique  homology orientation  $\tilde \omega$ 
of $M$  such that the torsion of
the exact homology sequence of the pair $(M,E)$
$$\Cal H = (H_3(E;\bold R)\to H_3(M;\bold R) \to H_3(M,E;\bold R)\to ... \to H_0(E;\bold R) \to H_0(M;\bold R))$$
has a positive sign. It is understood that the torsion  is taken with respect to arbitrary
bases in 
$H_*(E;\bold R), H_*(M;\bold R), H_*(M,E;\bold R)$ determining the orientations $\omega, \tilde \omega,  \omega_{(M,E)}$,
respectively. (The sign  of  $\tau (\Cal H)\in \bold R \backslash \{0\}$  depends   only on these orientations and does not depend
on the choice of the bases). 
 The homology orientation  $\omega^M$ of $M$ {\it induced by }
$\omega$ is defined by
$$\omega^M= (-1)^{m b_3(M)+ (b_0(E)+b_1(E)) (b_0(M)+b_1(M)+m +1) +b_3(E) (b_3(M)+1)}\, \tilde \omega. \eqno (3.2.a)
$$ The sign on the right-hand side is needed to ensure Lemma 3.3 below.

Clearly $(-\omega)^M=-(\omega^M)$. If $m=0$ then
$M=E$ and
$\omega^M=\tilde
\omega=\omega$. In the case of connected $E$,  formula (3.2.a) simplifies to 
$$\omega^M= (-1)^{m b_3(M) + (b_1(E)+1) (b_1(M)+m) }\, \tilde \omega \eqno (3.2.b) $$
We establish   two useful   properties of induced homology orientations.

\skipaline \noindent {\bf    3.3.  Lemma (transitivity). } 
 {\sl  
Let $E$ be a compact    3-manifold   whose boundary consists of tori.
Let $M\supset E$ be   obtained by gluing directed solid tori to $E$.  Let $M'\supset M$
be   obtained by gluing  directed solid tori to $M$.
Then for any homology orientation $\omega$ of $E$, } $$\omega^{M'} =(\omega^M)^{M'}.$$
 
\skipaline {\sl Proof.}        Assume that $M$ is 
obtained from $E$ by gluing $p$ directed solid tori  and $M'$ is 
obtained from $M$ by gluing $q$ directed solid tori. Fix a CW-decomposition of $M'$ such that
both $M$ and $E$ are its CW-subcomplexes. Consider the cellular chain complexes
$$c^1=C_*(E;\bold R) \subset c^2=C_*(M;\bold R) \subset  c^3=C_*(M';\bold R).$$
We  fix orientations of the meridional discs of the solid tori forming $M'\backslash E$ and provide  
 $H_*(c^2/c^1) ,  H_*(c^3/c^1)$,  $H_*(c^3/c^2) $
with    bases   as in Sect.\ 3.2.
Let us  provide  
$H_*(c^1), H_*(c^2) $ with   bases determining the homology orientations 
$\omega, \omega^M $, respectively.
The corresponding relative torsion 
$\tau (c^1\subset c^2)\in \bold R \backslash \{0\}$ depends on the choice of the bases but  its sign $\pm 1$, denoted 
$\tau_0 (c^1\subset c^2, \omega, \omega^M) $, depends only on $\omega, \omega^M$. 
The definition of   $\omega^M$ can
be reformulated by saying that
$$\tau_0 (c^1\subset c^2, \omega, \omega^M)
$$
$$= (-1)^{\theta ( c^2, c^1)+
p b_3(M)+ (b_0(E)+b_1(E)) (b_0(M)+b_1(M)+p +1) +b_3(E) (b_3(M)+1)} .$$
A direct computation shows that
$$\theta ( c^2, c^1)$$
$$=p b_3(M)+ (b_0(E)+b_1(E)) (b_0(M)+b_1(M)+p +1) +b_3(E) (b_3(M)+1) +p \, (\mod 2).$$
This follows from the definition of $\theta ( c^2, c^1)$ and the following equalities mod 2:
$$
\beta_i (c^2/c^1)=\cases p,~ {  {if}}\,\,\, i=2, \\
   0,~ {  {if}} \,\,\, i\neq 2 ,\endcases$$
 $$
 \beta_i (c^1)=\cases b_0(E),~ {  {if}}\,\,\, i=0, \\
  b_0(E)+b_1(E),~ {  {if}} \,\,\, i=1, \,\,\,   \\
  b_3(E),~ {   {if}} \,\,\, i= 2, \,\,\,  \\
0,~{  {if}} \,\,\,   i\geq 3,\endcases \,\,\,\,\, 
 \beta_i (c^2)=\cases b_0(M)=b_0(E),~ {  {if}}\,\,\, i=0, \\
  b_0(M)+b_1(M),~ {  {if}} \,\,\, i=1, \,\,\,   \\
  b_3(M),~ {   {if}} \,\,\, i= 2, \,\,\,  \\
0,~{  {if}} \,\,\,   i\geq 3.\endcases
$$
(Here we use that  $\chi(E)=\chi (\partial E)/2=0$ and $\chi(M)=\chi (\partial M)/2=0$.)
Thus, the definition of   $\omega^M$ can
be reformulated by saying that $$\tau_0 (c^1\subset c^2, \omega, \omega^M)= (-1)^p.\eqno (3.3.a)$$
Similarly,  
$$\tau_0 (c^2\subset c^3, \omega^M, (\omega^M)^{M'})=  (-1)^q,\,\,\, \tau_0 (c^1\subset c^3, \omega,
\omega^{M'})=  (-1)^{p+q}.$$ It is easy to compute that $\theta  (c^3/c^1, c^2/c^1)=0$ and 
 $\tau(c^2/c^1\subset c^3/c^1)=+1$. 
Applying  Lemma 2.3  and taking signs we obtain that 
$$\tau_0 (c^1\subset c^3, \omega, (\omega^M)^{M'})=
\tau_0 (c^1\subset c^2, \omega, \omega^M)\, \tau_0 (c^2\subset c^3, \omega^M, (\omega^M)^{M'}) $$
 $$= (-1)^{p+q}
=\tau_0 (c^1\subset c^3, \omega,
\omega^{M'}).$$
Therefore $(\omega^M)^{M'}=
\omega^{M'}$.

\skipaline \noindent {\bf    3.4.  Lemma. } 
 {\sl  Let $E$ be the exterior   of an oriented link $L=L_1\cup ...\cup L_m$ in an oriented 
 3-dimensional integral homology sphere  $N$. Let $E'\supset E$  be the exterior of a non-void sublink $L'$ of $L$.
Then $E'$ is obtained from $E$ by gluing directed solid tori and $(\omega_L)^{E'}= \omega_{L'}$.}

\skipaline {\sl Proof.}  Lemma 3.3 allows to deduce this lemma  by induction  from  the  case where $L'$ has one component
less than $L$.  We however  prefer to give a direct proof in the general case. Assume for concreteness that
$L'=L_1\cup ...\cup L_n$ with
$1\leq n\leq m$.  It is clear that $E'$ is obtained from $E$ by gluing $m-n$ directed solid tori $\Z_{n+1},..., \Z_m$ which are
regular neighborhoods of
$L_{n+1},..., L_m$, respectively.  We orient the meridional disc  of  each $\Z_i$ so that 
its intersection index  with $L_i$ equals $+1$.  This yields a basis in  
$H_* (E',E)$ as in Sect.\ 3.2. The homology $H_*(E), H_(E')$ have
  bases described  in Sect.\ 3.1. 
The exact homology sequence, $\Cal H$,  of the pair $(E',E)$   splits as a concatenation of two short exact sequences
$$ H_3(E',E)\to H_2 (E) \to H_2(E') , \,\,  H_2(E',E)\to H_1 (E) \to H_1(E') $$
and the  inclusion isomorphism $H_0 (E) \to H_0(E')$. 
The  torsions  of these pieces with respect to the  choosen bases are  equal to $(-1)^{mn+n}, (-1)^{mn+n}$ and $+1$,
respectively.  Hence   $\tau(\Cal H)=+1$ and 
therefore $ \tilde {\omega}_L = \omega_{L'}$. 
The numerical expression  appearing in 
(3.2.b)
equals 
$$  (m-n) b_3(E')+(b_1(E)+1) (b_1(E')+m-n)  \equiv 0\, (\mod 2).$$
Therefore the   sign   in (3.2.b) is $+$ and $(\omega_L)^{E'}= \tilde {\omega}_L  = \omega_{L'}$.

\skipaline \skipaline \centerline {\bf  4. Combinatorial Euler structures on 3-manifolds}

\skipaline \noindent {\bf   4.1. Combinatorial Euler structures.
}  In this and the next two subsections we recall the theory of combinatorial Euler structures and their torsions following [Tu4]. 
  Let $X$ be a finite   $CW$ space. 
An {\it Euler chain} in $X$  is a 1-dimensional singular chain $\xi$ in $X$ with 
$$\partial \xi = \sum_{a} (-1)^{\dim a} \alpha_a$$
where $a$ runs over all (open) cells of $X$ and 
  $\alpha_a$ is a point in   $a$. 
For two Euler chains $\xi, \eta$   we define  a \lq \lq quotient"   $\xi /\eta\in
H_1(X)$ as follows. Let  $$\partial \xi = \sum_{a } (-1)^{\dim a}
\alpha_a\,\,\,\,\, {\text {and}}\,\, \,\,\, \partial \eta = \sum_{a } (-1)^{\dim a}
\beta_a$$ where $\alpha_a, \beta_a\in a$. For each cell $a $, choose a path 
$\gamma_a:[0,1]\to a$ from $\alpha_a$ to $\beta_a$. Then $\xi /\eta\in
H_1(X)$ is the homology class of the   
singular 1-cycle $ \xi-\eta +\sum_{a} (-1)^{\dim a} \gamma_a$.  

We say that two Euler chains $\xi, \eta$ in $X$ define the same
{\it (combinatorial) Euler structure} on $X$ if  $\xi /\eta=1$. The Euler structure represented by   $\xi$ is denoted by
$[\xi]$. The set of Euler structures on
$X$ is denoted by
$\Eul (X)$. This set is non-void iff for every connected component $X_0$ of $X$ we have $\chi(X_0)=0$.
It is clear that $\Eul (X)=\prod_{X_0} \Eul (X_0)$ where $X_0$ runs over the components of $X$. 

The group $H_1(X)$ acts
on  
$\Eul (X)$: if $[h]\in H_1(X)$ is the homology class of a 1-cycle
$h$ and $\xi$ is an Euler chain  on $X$ then
$[h] [\xi]= [\xi +h]\in \Eul(X)$.  If $\Eul (X)\neq \emptyset$ then this action is free and transitive. 

For any cellular subdivision
$X'$ of $X$ there is a canonical $H_1(X)$-equivariant bijection
$\Eul(X)\to \Eul(X')$. This allows us to define
the set
of combinatorial Euler structures $\Eul (M)$ for a smooth  compact    manifold $M$. This set is obtained by identification of the
sets $\{\Eul (X)\}_X$ where $X$ runs over the $C^1$-triangulations of $M$. 
By [Tu4], there is a
canonical $H_1(M)$-equivariant bijection $\Eul (M)= \vect (M)$.

\skipaline \noindent {\bf   4.2. Torsions of Euler structures.
}   Let   $X$ be a   finite   connected $CW$ space with $\chi(X)=0$. Let $F$
be a field  and $\varphi: \bold Z[H_1(X )]\to F  $ be  a ring homomorphism.
For every Euler structure $e\in \Eul (X)$ and every homology orientation $\omega$ of $X$ we   define  a torsion  
$\tau^{\varphi} (X,e, \omega	)\in F$. 

  Consider  the maximal abelian covering $\tilde
X$ of $X$ with its induced $CW$ structure. The   group
$H=H_1(X)$   acts on $\tilde X$ via the covering transformations     permuting the
cells in $\tilde X$ lying over any given cell in $X$.  A family of cells in $\tilde X$ is
said to be {\sl fundamental}{} if   over each cell of $X$ lies exactly one cell of this
family.  Every  fundamental family of cells in  
$\tilde X$ gives rise to an Euler structure on $X$: consider a spider in $\tilde X$
formed by arcs  in $\tilde X$ connecting a   point $x\in \tilde X$ to   points
in these cells; the arc joining $x$ to a point of an odd-dimensional (resp.
even-dimensional) cell should  be 
 oriented towards $x$ (resp. out of $x$). Projecting this spider to $X$ we obtain an
Euler chain in $X$ representing an element of  $\Eul(X)$. It is clear that
any Euler structure on $X$ arises in this way from a fundamental family of cells in
$\tilde X$. We fix  a fundamental family of cells $\tilde e$ in
$\tilde X$ corresponding to $e\in \Eul (X)$ in this way.

We  orient and order the cells in the family $\tilde e$  in  an arbitrary way. This   yields  a basis for
the cellular  chain complex $C_*(\tilde X)=C_*(\tilde X;\bold Z)$ over   $\bold
Z[H ]$. Consider the induced basis in the chain complex  over $F$ $$C^{\varphi}_*(X)=F
\otimes_{\bold Z[H ]} C_*(\tilde X) $$ 
where $\bold Z[H ]$ acts on $F$ via $\varphi$.  If this based chain complex is acyclic, i.e.,   
 $H^{\varphi}_*(X)= H_*(C^{\varphi}_*(X))=0$, then
we  can consider  its torsion  $   \tau (C^{\varphi}_*(X)) \in F \backslash  0$.
Since the cells of $\tilde  e$
are in bijective correspondence with the cells of $X$, the choosen orientation and   order for
the cells of $\tilde  e$  induce   an  orientation and an order for the cells of  $X$. This
yields  a basis of the cellular  chain complex  $C_*(X;\bold R)$ over $\bold R$. Provide the
homology of $C_*(X;\bold R)$ with a 
basis  determining the  
homology orientation $\omega$. Consider the sign-refined torsion  $\check  \tau (C_*(X;\bold R))\in \bold
R\backslash \{0\}$ of the resulting based chain complex with based homology. We need
only its sign $ \pm 1$ denoted $\check  \tau_0 (C_*(X;\bold R))$.   Set
$$
\tau^{\varphi} (X,e, \omega	)=\cases \check  \tau_0 (C_*(X;\bold R))\, \, \tau (C^{\varphi}_*(X)) \in F\backslash \{0\},~ { 
{if}}\,\,\, H^{\varphi}_*(X)=0, \\
   0\in F ,~ {  {if}} \,\,\,H^{\varphi}_*(X)\neq 0.\endcases$$
 It turns out that $\tau^{\varphi} (X,e, \omega	)\in F$  does not depend on the auxiliary choices and  is invariant under cellular
subdivisions of $X$.   We have  $$\tau^{\varphi} (X,e, -\omega)=-\tau^{\varphi} (X,e, \omega)\,\,\, {\text {and}}\,\,\, 
 \tau^{\varphi} (X,he,   \omega)=\varphi(h)\,\tau^{\varphi} (X,e, \omega)  \eqno (4.2.a)$$  
 for any $e\in \Eul(X), h \in H_1(X)$,

We shall use the notation $\tau^{\varphi} (X,e)$ for $ \pm \tau^{\varphi} (X,e,\omega )$:
$$\tau^{\varphi} (X,e)= \pm \tau^{\varphi} (X,e,\omega )=\tau^{\varphi} (X,e,\pm \omega ).$$
We will  view  $\tau^{\varphi} (X,e)$ as  an element of $F$ defined up to multiplication by $-1$.  
 
The definition of $\tau^{\varphi} (X,e, \omega	)$ extends by multiplicativity to the case of    a  non-connected  finite CW space
$X$. All the connected components  of
$X$ should be  homology oriented. Then the torsion   corresponding to   $e\in \Eul (X)$ and 
a ring homomorphism
$\varphi: \bold Z[H_1(X )]\to F  $ is defined as the product of the torsions of the  components of $X$
corresponding to the restrictions of  $e$ and $\varphi$.  

In Section 6 we shall need a  version  $\check \tau^{\varphi}$  of $\tau^{\varphi}$ which is always non-zero. Assume that
the vector space $H^{\varphi}_*(X) $  is endowed with a basis $b$ and set 
$$\check \tau^{\varphi} (X,e, \omega; b	)=\check  \tau_0 (C_*(X;\bold R))\, \check \tau (C^{\varphi}_*(X), b)\in F
$$ where $\check  \tau_0 (C_*(X;\bold R))=\pm 1$ is the same sign as above and 
$\check \tau (C^{\varphi}_*(X), b)$ is the torsion $\check \tau$ of  $C^{\varphi}_*(X)$ corresponding to the basis $b$ in  
homology. The bases in the vector spaces of chains  are determined by $\tilde e$ as above. 
The torsion $\check \tau^{\varphi} (X,e, \omega; b	)$ does not depend on the auxiliary choices, never equals $0$ and  is
invariant under cellular subdivisions of $X$.   If $ H^{\varphi}_*(X)=0$ then  
$\check \tau^{\varphi} (X,e, \omega; \emptyset	)=\tau^{\varphi} (X,e, \omega)$.

  \skipaline \noindent {\bf    4.3. Example. } 
The 2-torus  $T=S^1\times S^1$ has a canonical  Euler structure,   $e^{\can}_T$, defined as follows. Set  $H=H_1(T)=\bold
Z^2$.
 Consider the inclusion $\mu$ of the group ring $\bold Z[H]$ into its 
field of quotients
$Q(H)$. 
It is  easy to check (using for instance a CW-decomposition
of $T$ consisting of one 0-cell, two 1-cells and one 2-cell)
that the chain complex $C^{\mu}_*(T)$ is acyclic and     
  for any $e\in \Eul (T)$,  $\tau^{\mu} (T,e )=\pm g_e$ with 
$g_e\in H$. By (4.2.a),
$g_{he}=hg_e$ for all $h\in H$. Therefore the formula $e\mapsto g_e$ establishes a bijection $\Eul (T)\to H$
 natural with respect to    diffeomorphisms $T \to T$.
The Euler structure $e\in \Eul (T)$ with $g_{e}=1\in H$ is  denoted $e^{\can}_T$. This  
is  the unique Euler structure  on $T$     invariant under all diffeomorphisms $T \to T$.  
Therefore under  the  identification $\Eul (T)=\vect (T)$, $e^{\can}_T$    corresponds to the class 
  $[w]$ described in Sect.\ 1.1.

One can check  using the CW-decomposition of $T$ as above, that for any field $F$ and any  
 ring homomorphism
$\varphi: \bold Z[H_1(T )]\to F  $ such that $\varphi (H_1(T ))\neq 1$, we have $$\tau^{\varphi} (T,e^{\can}_T )=\pm
1.
\eqno (4.3.a)$$

  \skipaline \noindent {\bf    4.4. Gluing of  combinatorial Euler structures. }
Let $M$ be a compact  3-manifold whose boundary consists of tori.
  Let  $T\subset M\backslash \partial M$ be a finite system of disjoint embedded 2-tori  splitting  $M$ into two 3-manifolds
$M_0,M_1$.
We define a gluing map $\cup:\Eul (M_0) \times \Eul (M_1)\to \Eul (M)$ as follows. Fix a CW-decomposition of $M$ so that
$M_0,M_1,T$ are its   CW-subcomplexes.  If $\xi_0, \xi_1$ are Euler chains on $M_0,M_1$, respectively, then
we set
$$ [\xi_0] \cup [\xi_1]= [\xi_0+\xi_1- \xi]\in \Eul (M)$$
where $\xi$ is an Euler chain on $T$ representing the canonical Euler structure  on each component of $T$, cf.  Sect.\ 4.3. 
Note that  $\xi_0+\xi_1- \xi$ is an Euler chain on $M$. It is clear that the gluing operation 
$\cup$ satisfies (1.5.b) where $e_0\in \Eul (M_0),e_1\in \Eul (M_1)$.  It can be deduced from definitions that the diagram
$$
\CD
 \Eul (M_0) \times \Eul (M_1)  @>\cup>>  \Eul (M) \\ 
  @V=VV    @VV=V        \\
\vect (M_0) \times \vect (M_1)  @>\cup>>  \vect (M)
\endCD $$
is commutative.

\skipaline \noindent {\bf    4.5.  Lemma. } 
 {\sl  Let $M, T, M_0, M_1$ be the same objects as in Sect.\ 4.4.  Let $F$
be a field  and $\varphi: \bold Z[H_1(M )]\to F  $ be  a ring homomorphism such that for  every  component $T'$ of $T$, we
have 
$\varphi (H_1(T'))\neq 1$.   Let} $\inc_i: \bold Z[H_1(M_i)]\to  \bold Z[ H_1(M)]$  {\sl  be the ring homomorphism induced by
the inclusion homomorphism $H_1(M_i)\to   H_1(M)$ where
$i=0,1$. Then for any}  $e_0\in \Eul (M_0) , e_1\in  \Eul (M_1)$, 
$$\tau^{\varphi} (M,e_0\cup e_1 )=\tau^{\varphi\,\inc_0} (M_0,e_0 ) \,\,\tau^{\varphi\,\inc_1} (M_1,e_1 ) \in F/\{\pm 1\}.
\eqno (4.5.a)$$

\skipaline {\sl Proof.}  Let us assume for simplicity that $M_0,M_1$ are connected and  $T$ is a 2-torus; the general case is  
similar.  Fix a CW-decomposition of
$M$ so that
$M_0,M_1,T$ are its   CW-subcomplexes. Denote the inclusion homomorphism $\bold Z[ H_1(T)]\to \bold Z[ H_1(M)]$ by
$\inc$.  We have the usual  short exact sequence of chain complexes
$$0\to C^{\varphi \,  \inc}_*(T)\to C^{\varphi\,\inc_0} (M_0)\oplus C^{\varphi\,\inc_1} (M_1) \to C^{\varphi}_*(M)\to
0.\eqno (4.5.b)$$ The  assumption  
$\varphi (H_1(T)) \neq 1$ implies   that 
$C^{\varphi \, \inc}_*(T)$ is acyclic. If $C^{\varphi}_*(M)$ is not acyclic then  at least  one of the complexes
$C^{\varphi\,\inc_0} (M_0),  C^{\varphi\,\inc_1} (M_1)$ is not acyclic and
both  sides of (4.5.a) are equal to 0. Assume that $C^{\varphi}_*(M)$  is  acyclic. Then all the complexes involved are acyclic
and the torsions in (4.5.a) are non-zero.

We now describe the gluing of combinatorial Euler structures 
in terms of fundamental families of cells. 
Let $\tilde T,  \tilde M_0, \tilde M_1, \tilde M$ be the maximal abelian coverings  of $T,M_0, M_1, M$, respectively. 
We have a commutative diagram
of cellular maps 
$$ 
\CD
\tilde T   @>f_0>>  \tilde M_0 \\ 
  @Vf_1VV                      @VVg_0V        \\
\tilde M_1   @>g_1>>  \tilde M.
\endCD $$
Each of these maps  is a lift of  the corresponding inclusion and is equivariant with  respect to the   inclusion
homomorphism in 1-homology.  For instance, the map $\tilde T  \to   \tilde M_0 $ is a lift of  the inclusion
$T  \hookrightarrow    M_0 $ and is equivariant with respect to the  inclusion homomorphism $H_1(T) \to H_1(M_0)$.
Choose a fundamental family of cells  $\tilde e$ in $\tilde T$ representing $e^{\can}_T$. 
It  is clear that $f_0(\tilde e)$ is a family of cells in $\tilde M_0$  such that over each cell of $T\subset M_0$ lies
exactly one cell
of this family.  We can  add to $f_0(\tilde e)$ a certain set  $\tilde  e_0$ of cells  in  $\tilde M_0 \backslash f_0(\tilde T)$
so that   $f_0(\tilde e)\cup \tilde  e_0$
is a fundamental family of cells in $\tilde M_0$ 
representing
$e_0$. Choose also a fundamental family of cells $\tilde  e_1$  in $\tilde M_1$ 
representing
$e_1$.   Then  $  g_0(\tilde  e_0) \cup g_1(\tilde  e_1)$
  is a fundamental family of cells in $\tilde M$. It follows from definitions that this family represents 
$e_0\cup e_1$. 

The  families of cells $\tilde e$,  $f_0(\tilde e)\cup \tilde  e_0$, $\tilde  e_1$, and 
$  g_0(\tilde  e_0) \cup g_1(\tilde  e_1)$ determine bases in the chain complexes
$C^{\varphi \,  \inc}_*(T),  C^{\varphi\,\inc_i} (M_i)$ with $i=0,1$, and $ C^{\varphi}_*(M)$, respectively.
It follows from (4.3.a) and definitions that the corresponding torsions are equal to $\pm 1$, $ 
\tau^{\varphi\,\inc_i} (M_i,e_i )$ with
$i=0,1$ and $\tau^{\varphi} (M,e_0\cup e_1 )$.
These bases in 
$C^{\varphi \,  \inc}_*(T),  C^{\varphi\,\inc_0} (M_0)\oplus C^{\varphi\,\inc_1} (M_1)$, and $ C^{\varphi}_*(M)$
are compatible in the sense of Sect.\  2.2 (at least up to sign).  Now, formula  (2.2.c) implies (4.5.a). 

\skipaline \noindent {\bf   4.6. Duality for torsions.} One of the fundamental properties 
  of the Reidemeister torsions is the duality due to Franz [Fr] and Milnor [Mi1].
We  state  the duality  theorem for torsions of Euler structures on 3-manifolds following [Tu4].  
Let $M$ be a    compact  orientable 3-manifold  whose boundary consists of tori. Let
$F$ be a field  with involution $f\mapsto \overline f:F\to F$ and      $\varphi: \bold Z[ H_1(M) ] \to
F  $ be  a ring homomorphism such that $\overline  {\varphi (h)}= \varphi (h^{-1})$ for
any $h\in H_1(M)$. Then for every $e\in \vect (M)=\Eul (M)$,
$$\overline {\tau^{\varphi} (M,e)}= \tau^{\varphi} \,(M,e^{-1})
= \varphi(c(e))^{-1}\, \tau^{\varphi} (M,e).\eqno (4.6.a)$$
All   torsions here are   elements of   $F$ defined up to multiplication by $-1$.
For closed 3-manifolds  formula (4.6.a) is established in [Tu4], Appendix B. Here we  use this result to prove   (4.6.a) in
the case
$\partial M\neq \emptyset$.  Denote by $K$ the double of $M$ obtained by gluing two copies of $M$
along the identity homeomorphism of the boundaries. Clearly, $K$ is a closed orientable 3-manifold.
Let $\psi: \bold Z[H_1(K)]\to  \bold Z[H_1(M)]$ be the ring homomorphism induced by the natural folding map $K\to M$ which is
the identity on both copies of
$M$.  The restrictions  of   $\varphi \psi $
to  both copies of $M$ are  equal to $\varphi$. 
Applying (4.6.a) to the Euler  structure $e\cup e$ on $K$ we obtain 
$$\overline {\tau^{\varphi  \psi} (K,e \cup e)}  = \tau^{\varphi  \psi} (K,(e \cup e)^{-1} )
= \tau^{\varphi  \psi} (K,e^{-1}\cup e^{-1} ).$$
By
Lemma
4.5,  
$$\tau^{\varphi  \psi} (K,e\cup e )=(\tau^{\varphi} (M,e ))^2, \,\,\,\,\,
 \tau^{\varphi  \psi} (K,e^{-1}\cup e^{-1} )=(\tau^{\varphi} (M,e^{-1}
))^2.$$
Thus,
$$\overline {\tau^{\varphi} (M,e )}^2
=    (\tau^{\varphi} (M,e^{-1}))^2 $$
and therefore  $\overline {\tau^{\varphi} (M,e )}=    \tau^{\varphi} (M,e^{-1} )$.

In Appendix 3 we   discuss a  more precise version of duality taking  into account homology orientations.

\skipaline \skipaline \centerline {\bf  5. The Torres formula for torsions}

\skipaline \noindent {\bf    5.1.  Lemma. } 
 {\sl  Let $E$ be a compact   connected  3-manifold   such that $\partial E$ consists of
tori. Let
$M$ be a   3-manifold obtained by   gluing a  directed solid
torus} $\Z$  {\sl  to $E$ and let  $h \in H_1(M)$ be the homology class of the core of} $\Z$. 
  {\sl  Let
$F$ be a field    and      $\varphi: \bold Z[ H_1(E)] \to
F $ be  a ring homomorphism  inducing a ring homomorphism  $\varphi^M: \bold Z[H_1(M) ]\to
F  $. Then for any} $e\in \Eul (E)  
$ {\it and any   homology orientation $\omega$ of $E$},
$$\tau^{\varphi} (E,e, \omega)=(\varphi^M(h)-1) \,\, {\tau^{\varphi^M} (M,e^M , \omega^M)}. \eqno (5.1.a)$$

\skipaline

The condition on $\varphi$   means that 
$\varphi$  maps the homology class  of the boundary of  the meridional disc
of $  \Z$ to $1$.  Lemma 5.1  is a version
of the   Torres identity for the Alexander polynomials of links in $S^3$.

We prove Lemma 5.1 in Sect.\ 5.3 after a few preliminary computations.

\skipaline \noindent {\bf   5.2.
Solid torus re-examined.}  Let
$\Z=S^1\times D^2$ be a directed solid torus  endowed with a CW-decomposition obtained from a  certain  
CW-decomposition of $\partial \Z$ by adjoining two cells:  a 2-cell $a^2=x\times D^2$ with $x\in S^1$ and a 3-cell $a^3$ with
interior
$(S^1\backslash
\{x\})
\times \Int (D^2)$. Let $ \tilde \Z=\bold R\times D^2$ be the universal   covering of $\Z$ with induced
CW-decomposition. We shall exhibit a  fundamental family  of cells in   $\tilde \Z$ representing the distinguished Euler structure 
$e_{\Z}$ on
$
 \Z$,  cf.  Sect.\ 1.4.

  Choose a fundamental family of cells   in  the maximal abelian covering of $\partial \Z$ representing  $e^{\can}_{\partial \Z} \in
\Eul (\partial \Z)$,  cf.  Sect.\  4.3.  Projecting this family to  $\partial \tilde \Z$ we obtain a family of cells, $\tilde e$, in $\partial
\tilde
\Z$ such that over each cell of
$\partial \Z
$ lies exactly one cell
of $\tilde e$. We lift   $a^2,a^3$ to cells $\tilde a^2, \tilde a^3$ of $\tilde \Z$ such that
$\partial \tilde a^3 =  \pm (h-1) \tilde a^2$ modulo 2-cells lying in $\partial \tilde \Z$. 
Here $h=h_{\Z}\in H_1(\Z)$ is the distinguished generator.
The family of cells   $ \tilde e, \tilde a^2, \tilde a^3 $ in $\tilde \Z$  is
fundamental and   represents  a certain Euler structure $
e\in \Eul (\Z) $.  We claim that $e=e_{\Z}$.

We  need to check that $c(e)=h^{-1}$.  One possible
proof consists in  a combinatorial computation of $c(e)$  from a fundamental
family of cells (or from an Euler chain) representing $e$. We  shall use another approach based on  (4.6.a).
Set $H=H_1(\Z)$ and   consider the inclusion $\mu$ of the group ring $\bold Z[H]$ into its 
field of quotients
$Q(H)$.  We claim that $$ \tau^{\mu} (\Z,e )=\pm   {1}/ {(h-1)}. \eqno (5.2.a)$$
This   would imply the equality $c(e)=h^{-1}$
since
by (4.6.a),
$$ \pm c(e)=\pm \mu (c(e))=  {\tau^{\mu} (M,e)}\, /\,{\overline {\tau^{\mu} (M,e)}  }= 
\pm   {(h^{-1}-1)}/ {( h-1)}=\pm h^{-1}  $$
 where the overline denotes the involution in $Q(H)$ sending each $g\in H$ to $g^{-1}$.

To prove (5.2.a) we consider a short exact sequence of acyclic chain complexes
$$0\to C^{\mu\,\inc}_* (\partial \Z)\to C^{\mu}_* (\Z) \to C^{\mu}_* (\Z,\partial \Z)\to 0$$
where $\inc: \bold Z[H_1(\partial \Z)] \to \bold Z[ H ]$ is the inclusion homomorphism.
We orient and order the   cells  of 
$\tilde e$ so that they determine  a basis of $C^{\mu\,\inc}_* (\partial \Z)$. 
By (4.3.a),  we have  $\tau (C^{\mu\,\inc}_* (\partial \Z))=\tau^{\mu\,\inc} (\partial \Z,e^{\can}_{\partial \Z} )=\pm 1$. 
The cells $\tilde a^2, \tilde a^3$ endowed with arbitrary orientations yield a basis  for   
 $C^{\mu}_* (\Z,\partial \Z)=  (Q(H)\, \tilde a^3 \to  Q(H) \,\tilde a^2)$. 
The torsion of this chain complex is equal to
$\pm (h-1)^{-1}$.  The family $ \tilde e, \tilde a^2, \tilde a^3 $  determines  a basis in $C^{\mu}_* (\Z)$. By   (2.2.c),  
$$\tau^{\mu} (\Z,e )=\pm \tau (C^{\mu}_* (\Z))=\pm  \tau (C^{\mu\,\inc}_* (\partial \Z))\, \tau(C^{\mu}_* (\Z,\partial \Z))= \pm
(h-1)^{-1}.$$
   
\skipaline \noindent {\bf   5.3.  Proof  of  Lemma 5.1.}  We fix a  CW-decomposition of
$E$ (such that 
$\partial E$ is a subcomplex) and extend it to a CW-decomposition  of $M$ by adjoining two cells $a^2,a^3\subset \Z$ as in Sect.\
5.2.   
 Choose a fundamental family of cells    in  the maximal abelian covering  of   $E$  representing $e\in \Eul (E)$. 
Projecting these cells to the maximal abelian covering, $\tilde M$, of $M$  we obtain a family of cells, $\tilde e$,  in $\tilde M$
such that over each cell of $E $ lies exactly one cell
of $\tilde e$.  We lift $a^2, a^3$ to cells $\tilde a^2, \tilde a^3\subset \tilde M$  
so that
$\partial \tilde a^3 =\pm (h-1)\, \tilde a^2$ modulo 2-cells lying over $ \partial \Z$. It is clear that
$\tilde e, \tilde a^2, \tilde a^3$ is a fundamental family of cells in $\tilde M$.  As in  Sect.\   4.5, one can deduce from
  the results of Sect.\ 5.2 that this family represents $e^M=e\cup  e_{\Z}$.
We orient and order the cells belonging to  $\tilde e$ in an arbitrary way. The cell $a^2$ is a meridional disc of $\Z$, we orient it
in an arbitrary way and   lift this orientation to $\tilde a^2$.   The product orientation of $\Z= S^1\times a^2$ induces 
orientations in $a^3$ and $\tilde a^3$. With   these orientations,  
$\partial \tilde a^3 =  (h-1)\, \tilde a^2$  modulo 2-cells lying over $ \partial \Z$. 

Consider the   chain complexes $C$,  $C'\subset C$ and $C''$ defined by 
$$C'= C^{\varphi }_* (E), \,\,C= C^{\varphi^M}_* (M),\,\, C''  =C/C'. \eqno (5.3.a)$$  The family of cells $\tilde e, \tilde a^2, \tilde
a^3$ determines    bases in $C,C',C''$ in a compatible way.   The 
non-zero part of   $C''$ amounts to the  boundary   homomorphism
$F\tilde  a^3 \to F \tilde  a^2$ sending $\tilde  a^3 $ to $(\varphi^M (h)-1) \,\tilde  a^2$. Now we consider two cases.

Case    $\varphi^M(h)\neq 1$. In this case 
$C''$ is acyclic and $\tau (C'')= (\varphi^M(h)-1)^{-1}$. The  complexes  $C$ and $C'$ are acyclic or not simultaneously. If
they are not acyclic then both  sides of (5.1.a) are equal to 0.  Assume that they are acyclic. 
By (2.2.c),
$$\tau(C )   =(-1)^{\nu (C,C')} \tau (C') \, \tau (C'') =(-1)^{\nu (C,C')}
 \tau (C')\, (\varphi^M(h)-1)^{-1}.\eqno (5.3.b)$$
Consider the   chain complexes $c$,  $c'\subset c$ and $c''$ defined by 
$$c'= C_* (E;\bold R) ,\,\, c= C_* (M;\bold R),\,\, c''=c/c'= C_* (M,E;\bold R). \eqno (5.3.c)$$ 
The same ordered families of  oriented cells as in the previous paragraph determine bases in $c', c, c''$.  We provide
their homology with bases determining the homology orientations $\omega, \omega^M, \omega_{(M,E)}$, respectively (cf.
Sect.\ 3.2).  By (2.2.a),
$$ \check \tau_0(c) = (-1)^{\nu (c,c')}\, \check \tau_0(c') \, \check \tau_0(c'')\, 
\tau_0 (c' \subset  c )$$
 where  the subindex $0$ indicates that we consider   the signs of the corresponding torsions. Observe  that $\nu (c,c')=\nu
(C,C') $ and by (3.3.a),  
$\tau_0 (c' \subset  c )=-1$. The non-trivial part of  $c'' $ amounts to 
 a zero boundary homomorphism
$ \bold R\, a^3  \to \bold R\, a^2$. The oriented cells   $a^2, a^3$ determine a 
  basis in   $H_*(c''_2)=H_*(M,E;\bold R)$ representing 
$\omega_{(M,E)}$.   Hence $\tau_0(c'')=1$,  $N(c'')= 1$, and 
$$\check \tau_0(c'')= (-1)^{N(c'')} 
\tau_0(c'')=-1.$$  We conclude that
$$\check \tau_0(c) =(-1)^{{\nu (C,C')}}\, \check \tau_0(c') .$$ Multiplying   this by    (5.3.b),
we obtain 
$$\tau^{\varphi } (E,e, \omega)= \check
\tau_0(c') \,\tau (C') =(\varphi^M(h)-1)\, \check \tau_0(c) \, \tau(C ) $$
$$=(\varphi^M(h)-1)\,  {\tau^{\varphi^M} (M,e^M ,
\omega^M)}.$$
 
Case    $\varphi^M(h)= 1$. We need to show that $\tau^{\varphi } (E,e, \omega)=0$ or equivalently that
$H_*^\varphi(E)\neq 0$.  Let us assume that 
 $H_*(C')=H_*^\varphi(E)= 0$ and look for a contradiction. It is clear that $H_i(C'')=F$ for $i=2,3$ and 
$H_i(C'')=0$ for $i\neq 2,3$. The exact homology sequence of the pair $(C,C')$ implies
that  $H_i(C) =F$ for $i=2,3$ and 
$H_i(C)=0$ for $i\neq 2,3$. If $\partial M\neq \emptyset$   then $M$ collapses onto a 2-dimensional
subcomplex which contradicts   $H_3(C)=F$. If $M$ is closed then it has a CW-decomposition with only one 3-cell.
For an appropriate choice of bases in $C_2,C_3$, the matrix of the boundary homomorphism $C_3\to C_2$ is then a row
$[\varphi^M(g_1)-1,...,
\varphi^M(g_m)-1]$ where $g_1,...,g_m$ are the generators of $H_1(M)$ dual to the 2-cells.
The equality $H_3(C)=F$ implies that $\varphi^M(g_1)=...=\varphi^M(g_m)=1$ so that $\varphi^M(H_1(M))=1$.
However, in this case $H_0(C)=F$ which is a contradiction.

\skipaline
 \skipaline \centerline {\bf  6.  Additivity of   torsions}

\skipaline \noindent {\bf   6.1. The setting.} 
Let $E$ be a compact  connected  oriented  3-manifold   such that $\partial E$ consists of
tori. Let $T$ be a component of $\partial E$ and   $\alpha$ be an 
 oriented  non-contractible
simple closed curve on $T$.  Fix   $\varepsilon =\pm 1$.
We can glue a   solid torus $\Z$
  to $E$ along   $  T$ so that 
$\alpha$ bounds a meridional disc in $\Z$. 
The resulting 3-manifold, $M$,  depends only on $\alpha$.  The orientation of
$E$ extends to  
$M$.  We make   $\Z=\overline {M\backslash E}$
directed by  first  orienting its meridional disc, $a^2$, so that  
$\partial a^2=\alpha$ in the oriented category and   then   orienting the core circle of $\Z$  
    so that  the   orientation in $M$ restricted to $\Z$   is $\varepsilon  $ times the product
orientation. Thus,  the manifold $M$ is obtained from $E$ by gluing a directed solid torus. 
We shall denote this manifold by  $M_{E,\alpha,\varepsilon}$.
Note that   the  manifolds $M_{E,\alpha,\pm \varepsilon}$ have the same underlying  manifold $M$ and differ   only by the
structure of a directed solid torus on $\overline {M\backslash E}$.

Fix   a ring homomorphism   $\varphi$ from $\bold Z[ H_1(E )] $  to  a field $F$ such that 
$\varphi(H_1(T))=1$. There is a  (unique) ring homomorphism $\varphi^M: \bold Z[ H_1(M)]\to 
F  $ whose composition with the inclusion homomorphism $\bold Z[ H_1(E)]\to  \bold Z[H_1(M)]$ is equal to $\varphi$. 

 \skipaline \noindent {\bf    6.2.  Lemma. } 
 {\sl  Let $E, T, \alpha, \varepsilon, M=M_{E,\alpha,\varepsilon}, F, \varphi , \varphi^M$  be   as in Sect.\  6.1.   Suppose
that
 $H_*^{\varphi^{M }}(M)=0$.  Let $\hat E\to E$ be the regular covering of $E$ corresponding to the
group} $H_1(E) \cap \varphi^{-1} (1)\subset H_1(E)$.  {\sl Then $T$ lifts to a torus $\hat T\subset \partial \hat E$ and
$\alpha$ lifts to an oriented simple closed  curve $\hat \alpha$ on $\hat  T$ such that their homology classes    $[ \hat\alpha] \in
H_1^\varphi (E) ,  [ \hat T] \in H_2^\varphi (E)$ form  a basis in 
$H_*^\varphi (E)$. (Here the orientation of $T, \hat T $ is induced by the one in $E$). 
For any 
    homology orientation
$\omega$ of
$E$ and any}
$e\in \Eul (E)$,
$$
 \tau^{\varphi^{M }}
(M ,e^{M },\omega^{M })
=-\varepsilon\, \tau^{\varphi} (E,e ,\omega; [ \hat \alpha]   ,  [ \hat T] ). \eqno (6.2.a)$$

\skipaline {\sl Proof.}  We  
fix a  CW-decomposition of
$E$ (such that  $\alpha$ and 
$T$ are subcomplexes) and extend it to a CW-decomposition  of $M$ by adjoining two cells $a^2,a^3\subset \Z$ as in Sect.\
5.2.  We can assume that $\partial a^2=\alpha$. We orient $a^2$ so that $\partial a^2=\alpha$ in the oriented category.
We provide $a^3$ with the orientation determined by the product orientation in $\Z$. (This orientation in $a^3$ depends on
$\varepsilon$). Let  
$\tilde e, \tilde a^2, \tilde a^3$ be a fundamental family of  ordered oriented cells in the maximal abelian covering $\tilde M$ of
$M$ constructed as in Sect.\  5.3 and  representing
$e^M=e\cup  e_{\Z}$.

Consider  the   chain complexes $C, C', C''$ defined by (5.3.a).  We shall  apply to them formula (2.2.a).   The family
of cells
$\tilde  e, \tilde  a^2, \tilde  a^3$ defines compatible bases     in  $C, C', C''$ in the usual way.  The homology 
of $C, C', C''$ are
provided with bases as follows. By assumption,  $H_*(C)=H_*^{\varphi^{M }}(M)=0$. Since $ \varphi^M(H_1(\Z))=1$, the 
non-trivial part of   $C''$ amounts to a zero   homomorphism  
$ F\tilde  a^3 \to F \tilde  a^2$.     Hence  $H_i(C'')=F [\tilde  a^i]$ for  $i=2,3$. We 
fix the basis  $ [\tilde  a^2],  [\tilde  a^3]$ in $H_*(C'')$.   Using the exact homology sequence of the pair $(C, C')$ 
we obtain that  $H_i (C' ) =0$ for
$i\neq 1,2$ and 
$H_i (C' ) =F [\partial  \tilde  a^{i+1}]$ for $i=1,2$. 
We fix the basis     $ [\partial \tilde  a^2],  [\partial \tilde  a^3]$ in $H_* (C' ) $.  
   A direct computation shows that $\check \tau
(C'')=\tau (C'\subset C)=-1$. Hence
$$\tau (C)= \check \tau (C)= (-1)^{\nu(C,C')}
\check \tau (C') \, \check \tau (C'') \, \tau (C'\subset C)$$
$$   = (-1)^{\nu(C,C')}
\check \tau (C';  [\partial \tilde  a^2],  [\partial \tilde  a^3])  \eqno (6.2.b)
$$
where we inserted $ [\partial \tilde  a^2],  [\partial \tilde  a^3]$    to keep track of the fixed basis in 
$H_*(C')=H_*^\varphi (E)$.

The same argument as in  Section 5.3 shows that 
$$\check \tau_0(C_* (M;\bold R))  =(-1)^{{\nu(C,C')}}\, \check \tau_0(C_* (E;\bold R)) $$
where the torsions are taken with respect
to
the   bases of chains determined by  $\tilde e,   a^2,  a^3$ and  bases  in homology determining the homology orientations
$\omega,
\omega^M$.  Multiplying by (6.2.b) we obtain,
$${\tau^{\varphi^M} (M,e^M , \omega^M)}=    \check \tau_0(C_* (M;\bold R)) \, \tau(C )$$
$$=  \check \tau_0(C_* (E;\bold R)) \,\check \tau (C';  [\partial \tilde  a^2],  [\partial \tilde  a^3])  =
 \tau^{\varphi } (E,e, \omega;  [\partial \tilde  a^2],  [\partial \tilde  a^3]). $$

It remains to give an interpretation of   $[\partial \tilde  a^2],  [\partial \tilde  a^3]$ in terms of  
 $\hat  E $.
The  complex $C'$     can be computed  from 
  $\hat E$.  Indeed, a natural projection from the maximal abelian covering $\tilde E$  of $E$ to $\hat E$ induces  the 
equalities 
$$C'=C^\varphi_* (E) = F\otimes_{\bold Z [H_1(E)]} C_*(\tilde E )=F\otimes_{\bold Z [G]} C_*(\hat E )$$
where $G= H_1(E)/(H_1(E) \cap \varphi^{-1} (1))= \varphi(H_1(E))$ is the group of covering transformations of  
$\hat E  $. The assumption
$\varphi(H_1(T))=1$ implies that   $T$ lifts to a  torus $\hat T \subset  \partial \hat E$.   
We provide $\hat T$ with the orientation       induced by the one on $T\subset \partial E$.
Since the projection
$\hat T
\to T$ is a homeomorphism, the curve $\alpha$ lifts to a curve $\hat \alpha$ on $\hat T$.
Then $\hat \alpha$ and $\hat  T$ represent certain  homology classes    $[ \hat\alpha] \in H_1^\varphi (E) , 
[ \hat T] \in H_2^\varphi (E)$.

The covering $\hat E\to E$ extends to a covering $\hat M\to M$ with the same group of covering transformations
$G$. The torus $\hat T$ bounds in $\hat M$ a solid torus, $\hat \Z$, projecting homeomorphically onto
$\Z$, and $\hat M=\hat E \cup_{g\in
G} g\hat \Z$. 
There is a covering  $\tilde M \to \hat M$   mapping   $\tilde  a^2$ onto a  meridional disc of $\hat  \Z$
bounded by
$\hat
\alpha$ and mapping $ \tilde a^3$ onto a   3-cell filling in $\hat \Z$.  
Then   $[\partial  \tilde   a^2]=[\hat \alpha]\in H_1^\varphi (E)$   and 
$[\partial  \tilde   a^3]=[\partial \hat \Z]=-\varepsilon\, [\hat T ]\in H_2^\varphi
(E)$.
Therefore
$${\tau^{\varphi^M} (M,e^M , \omega^M)}=    
 \tau^{\varphi } (E,e, \omega;  [\partial \tilde  a^2],  [\partial \tilde  a^3])=
-\varepsilon\, \tau^{\varphi } (E,e, \omega; [\hat \alpha],  [\hat T ]). $$

 \skipaline \noindent {\bf    6.3.  Remark. } 
The fact that the torsion $ \tau^{\varphi} (E,e ,\omega; [ \hat \alpha]   ,  [ \hat T] ) $   in Lemma 6.2 
does not depend on the choice of   $\hat T$ 
follows directly from the obvious equality
$\tau^{\varphi} (E,e ,\omega; g [ \hat \alpha]   ,  g [ \hat T] )=\tau^{\varphi} (E,e ,\omega; [ \hat \alpha]   ,  [ \hat T] )$
where 
$g\in G=\varphi(H_1(E))$.

\skipaline \noindent {\bf    6.4.  Lemma. } 
 {\sl  Let $E$ be a compact     3-manifold   such that $\partial E$ consists of
tori. Let $T$ be a component of $\partial E$.  Let $F$
be a field  and $\varphi: \bold Z [ H_1(E )]\to F  $ be a ring homomorphism such that 
$\varphi(H_1(T))=1$. Let  $\alpha_1,\alpha_2, \alpha_3$    be  
 oriented  non-contractible
simple closed curves on $T$ such that $\alpha_1$ is homological to   $\alpha_2 \alpha_3$. 
Fix  $\varepsilon_1, \varepsilon_2, \varepsilon_3=\pm 1$ and   set $M_i= M_{E,\alpha_r,\varepsilon_r}$ for  $r=1,2,3$. Then
for   any}
$e\in \Eul (E)$ {\sl and any   homology orientation
$\omega$ of $E$,}
$$
\varepsilon_1\, \tau^{\varphi^{M_1}}
(M_1,e^{M_1},\omega^{M_1})$$
$$
=\varepsilon_2\, \tau^{\varphi^{M_2}} (M_2,e^{M_2},\omega^{M_2})
+ \varepsilon_3 \,\tau^{\varphi^{M_3}} (M_3,e^{M_3},\omega^{M_3})
. \eqno (6.4.a)$$

\skipaline 

  {\sl Proof.} Let $\hat E\to E$ be the same covering of $E$ as in Lemma 6.2 with   group of covering transformations
   $G= \varphi (H_1(E))$.  Let $\hat T\subset \hat E$ be a lift of  $  T$  and   $\hat  \alpha_1,
\hat  \alpha_2,
\hat \alpha_3$  be the lifts
of $ \alpha_1, \alpha_2, \alpha_3$ to $\hat  T$.  Consider the elements 
$[\hat  \alpha_1], [ \hat  \alpha_2],
[\hat  \alpha_3]$ of  $H_1^\varphi(  E)=H_1(F\otimes_{\bold Z [G]}
C_*(\hat E ))$ represented by $\hat  \alpha_1, \hat  \alpha_2,
\hat  \alpha_3$, respectively. The assumptions of the lemma imply that  $[\hat  \alpha_1] = [ \hat  \alpha_2] +
[\hat  \alpha_3] $.

 If all three torsions entering  (6.4.a) are equal to $0$ then    (6.4.a) is
obvious.  Assume from now on that at least one of these  three torsions   is non-zero.  Then
the   proof of  Lemma  6.2 shows that $ H_i^\varphi (E)=0$ for $i\neq  1,2$ and  $H_1^\varphi (E),
H_2^\varphi (E)$ are 1-dimensional vector spaces over $F$. Moreover,  $H_2^\varphi (E)=F \,[\hat T]$. It is easy to see that
$[\hat
\alpha_r]\neq 0$   iff $H_*^{\varphi^{M_r}} (M_r )=0$ iff
$\tau^{\varphi^{M_r}} (M_r,e^{M_r},\omega^{M_r})\neq 0$.   By assumption, at least one of the classes $[\hat  \alpha_1], [ \hat  \alpha_2] ,
[\hat  \alpha_3]$ is non-zero.

 If all three   classes $[\hat  \alpha_1], [ \hat  \alpha_2] ,
[\hat  \alpha_3]$ are non-zero then by Lemma 6.2
$$\varepsilon_1 \,\tau^{\varphi^{M_1}}
(M_1,e^{M_1},\omega^{M_1})
=-  \tau^{\varphi} (E,e ,\omega; [ \hat \alpha_1]   ,  [ \hat T] )$$
$$
= -  \tau^{\varphi} (E,e ,\omega; [ \hat \alpha_2]   ,  [ \hat T] ) -  \tau^{\varphi} (E,e ,\omega; [ \hat \alpha_3]   ,  [ \hat T] )
$$
$$=\varepsilon_2 \,\tau^{\varphi^{M_2}} (M_2,e^{M_2},\omega^{M_2})
+ \varepsilon_3 \,\tau^{\varphi^{M_3}} (M_3,e^{M_3},\omega^{M_3})$$
where the second equality follows from    the formula  $[\hat  \alpha_1] = [ \hat  \alpha_2] +
[\hat  \alpha_3] $.
If $[\hat  \alpha_1]=0$ then  $[\hat  \alpha_2]=- [\hat  \alpha_3]\neq 0$.
Then  $ \tau^{\varphi^{M_1}}
(M_1,e^{M_1},\omega^{M_1})=0$ and
$$\varepsilon_2 \,\tau^{\varphi^{M_2}} (M_2,e^{M_2},\omega^{M_2})
=
 -  \tau^{\varphi} (E,e ,\omega; [ \hat \alpha_2]   ,  [ \hat T] ) $$
$$=  \tau^{\varphi} (E,e ,\omega; [ \hat \alpha_3]   ,  [ \hat T] )
=- \varepsilon_3\, \tau^{\varphi^{M_3}} (M_3,e^{M_3},\omega^{M_3})$$
which proves  (6.4.a). The cases where $[\hat  \alpha_2]=0$ or $[\hat  \alpha_3]=0$ are similar.

\skipaline \skipaline \centerline  {\bf 7.  The torsion $\tau$}
\skipaline
We discuss   the \lq \lq maximal abelian" torsion  $\tau$ introduced in [Tu2], see also  [Tu5], [Tu7].

\skipaline \noindent {\bf   7.1. Algebraic preliminaries.} For a unital commutative ring $K$ we denote by $Q(K)$ its {\sl classical
ring of fractions}, i.e., the localisation of $K$
by the multiplicative system of all non-zerodivisors. The natural map $K\to Q(K)$ is   injective. 

 For a finitely generated
abelian group $H$, set
$Q(H)=Q(\bold Z [H])=Q(\bold Q [H])$.  This commutative ring  splits (uniquely) as a direct sum of fields.
This is  obvious in the case of finite $H$ since then  $\bold Q[H]$ is a finite sum of cyclotomic fields so that
$Q(H)= \bold Q[H]$. In the general case such a splitting comes from 
a splitting of $H$ as a direct sum of $\Tors\,H$ and the free abelian group $G=H/\Tors\,H$. 
The ring $\bold Q[ \Tors\,H]$ splits as a sum of cyclotomic fields $\{K_r\}_r$. Hence
 $$\bold Q[H]= (\bold Q[ \Tors\,H]) [ G]= \bigoplus_r K_r [G].$$
  Each ring 
  $K_r [G]$ is a  domain and therefore  $$Q(H)= \bigoplus_r  F_r  \eqno (7.1.a)$$ 
where $F_r =Q(K_r [G])$ is the field of fractions of  $K_r [G]$. A useful fact: for any   $h\in H$
of infinite order, the element $h-1\in \bold Z[H]$ is a non-zerodivisor and therefore is invertible in $Q(H)$. 

One of the summands in (7.1.a) is the field $Q(H/\Tors\, H)$. 
The projection $Q(H)\to Q(H/\Tors\, H)$ is induced by the projection $H\to H/\Tors\, H$.
The inclusion $Q(H/\Tors\, H) \hookrightarrow Q(H)$ maps   $1\in Q(H/\Tors\, H)$ into 
  $\vert \Tors\,H\vert^{-1} \sum_{h\in \Tors\,H} h$.

Consider an increasing   filtration $$(\bold Z[H])_0=\bold Z[H] \subset (\bold Z[H])_1\subset (\bold Z[H])_2\subset ...\subset
Q(H)  \eqno (7.1.b)$$ where $(\bold Z[H])_k$ is   the set of $q\in Q(H)$ such that for any $h_1,...,h_k \in H$
we have $q\prod_{i=1}^k (h_i-1)  \in \bold Z[H]$.    We shall need the following easy fact:  any 
ring homomorphism $\varphi$ from $\bold Z[H]$ to a field $F$ such that $\varphi(H)\neq 1$ uniquely extends
to a ring homomorphism $\varphi_\#:\cup_k (\bold Z[H])_k \to F$.  If $\rank \, H\geq 2$ then 
the filtration (7.1.b) is
trivial, i.e.,  $(\bold Z[H])_k=\bold Z[H]$ for all $k$ (see [Tu5, Lemma 4.1.1]) and  $\varphi_\#=\varphi
:\bold Z[H]\to F$.

\skipaline \noindent {\bf   7.2. Definition of   $\tau$.} Let $X$ be a  finite connected $CW$-space. Set $H=H_1(X)$ and
denote by   $\varphi_r$ the composition of the inclusion 
$\bold Z[H] \hookrightarrow Q(H)$ and the projection $Q(H) \to F_r$ on the $r$-th term in   the splitting 
(7.1.a). 
By  Sect.\ 4.2, for any $e\in \Eul(X)$ and any homology
orientation $\omega$ of $X$, we have a torsion 
$\tau^{\varphi_r} (X, e, \omega) \in F_r$.
Set 
$$\tau (X,e, \omega)= \bigoplus_r  \tau^{\varphi_r} (X, e, \omega)\,\, \in \,\, \bigoplus_r F_r =Q(H).$$
This is a well defined element of   $Q(H)$. 
For  any $h\in H$,  
$\tau (X,h e, \omega)=h \,\tau (X,e, \omega)$.

\skipaline \noindent {\bf   7.3. The torsion $\tau$ for 3-manifolds.} 
 We shall need three basic lemmas  concerning the torsion $\tau
$ for 3-manifolds.

 \skipaline \noindent {\bf    7.3.1.  Lemma. } 
 {\sl   Let 
$M$ be a 
compact connected orientable 
3-manifold whose
boundary  is empty or  consists of  tori.  Let} $e\in \Eul (M)$ {\sl and $\omega$ be a homology
orientation of $M$. Set $H=H_1(M)$.
If $b_1(M) \geq 2 $, then $\tau (M,e,\omega) \in \bold Z [H]$}.  {\sl If $b_1(M) =1 $ and $\partial M\neq \emptyset$,
then $\tau (M,e,\omega) \in (\bold Z[H])_1$. If $b_1(M) =1 $ and $\partial M= \emptyset$, then $\tau (M,e,\omega) \in (\bold
Z[H])_2$.}
 
 \skipaline Lemma 7.3.1 follows from the results of [Tu5, Section 4]. Although we shall not need it, note that the inclusion $ 
\tau (M,e,\omega) \in Q(H)$ can be improved also in the case $b_1(M)=0$: In this case $\tau
(M,e,\omega) \in \vert \Tors\,H\vert ^{-1} \bold Z [H]$. 
 
 \skipaline \noindent {\bf    7.3.2.  Lemma. } 
 {\sl   Let  under the conditions of Lemma 7.3.1,  $\varphi$ be a ring homomorphism from $\bold Z[H]$ to a field $F$ such that
$\varphi(H)\neq 1$ and} $\char  F= 0$. {\sl 
Then 
$\tau^\varphi (M,e,\omega)=\varphi_\#(\tau (M,e,\omega))$.  }

\skipaline

  In the case 
$\tau^\varphi (M,e,\omega)\neq 0$ this lemma is   contained  in  [Tu7, Theorem 13.3] in a 
  general setting of   CW-spaces of any dimension. It remains only to show that 
if $\tau^\varphi (M,e,\omega)=0$, i.e., if $H_*^\varphi(M)\neq 0$, then   $\varphi_\#(\tau (M,e,\omega))=0$. This easily
follows from the computations in [Tu5, Sect.\ 4.1.2] and the assumption 
$\varphi(H)\neq 1$.

\skipaline \noindent {\bf    7.3.3.  Lemma. } 
 {\sl  Let $E$ be a compact   connected  orientable 3-manifold   such that $\partial E$ consists of   
tori and $b_1(E) \geq 2$. Let
$M$ be a   3-manifold with $b_1(M) \geq 1$ obtained by   gluing $m$  directed solid
tori to $E$ and let  $h_1,...,h_m\in H_1(M)$ be the homology classes of
the  core circles of these solid tori. Let}
      $\inc: \bold Z[ H_1(E) ] \to
\bold Z[ H_1(M) ] $  {\sl   be the inclusion homomorphism.  
Then for any} $e\in \Eul (E)  
$ {\it and any   homology orientation $\omega$ of $E$},
$$ {\inc} (\tau (E,e, \omega))= \prod_{i=1}^m (h_i-1)\,{\tau  (M,e^M , \omega^M)}. \eqno (7.3.a)$$

\skipaline 

  {\sl Proof.}   Consider the splitting 
$Q(H_1(M))= \bigoplus_r  F_r $ into a direct sum of fields. Denote  the projection
$Q(H_1(M))
\to F_r$ by
$p_r$. It suffices to prove that for all $r$, 
$$p_r ({\inc} (\tau (E,e, \omega)))=\prod_{i=1}^m (p_r (h_i)-1)\,  p_r (\tau  (M,e^M , \omega^M))
.$$ Let $\mu$ be   the  inclusion $\bold Z[H_1(M)] \hookrightarrow Q(H_1(M))$.  Applying    Lemmas  7.3.2  twice and Lemma
5.1 inductively $m$ times   we obtain that 
$$p_r ({\inc} (\tau (E,e, \omega)))= (p_r \,\mu\,\inc) (\tau  (E,e, \omega))
= (p_r \,\mu\,\inc)_\# (\tau  (E,e, \omega))=\tau^{p_r \mu\,\inc} (E,e, \omega)$$
$$  = \prod_{i=1}^m (p_r (h_i)-1)   \, \tau^{p_r \mu}  (M,e^M ,
\omega^M)=\prod_{i=1}^m (p_r (h_i)-1)  \, p_r (\tau  (M,e^M ,
\omega^M)).$$

\skipaline \skipaline \centerline  {\bf 8.  The Alexander-Conway function and derived invariants}

\skipaline \noindent {\bf    8.1.  The Alexander-Conway function. } Let 
$L=L_1\cup ...\cup L_m$  be 
  an ordered oriented link 
 in an oriented  3-dimensional integral homology sphere $N$.  The
Alexander-Conway function  $\nabla_L $
 of $L$ is a rational function on $m$ variables $t_1,...,t_m$ with integer coefficients. Considered up to sign, this 
function   is equivalent to the $m$-variable Alexander polynomial of $L$.  For links in $S^3$, the function  $\nabla_L
$   was introduced by Conway [Co] (see also  [Ha]);  it was extended to links in homology spheres in [Tu3]. 
We recall  here the definition of   $\nabla_L $ following 
[Tu3, Sect.\ 4].

  Let $E$ be  the exterior of  $L$. The group $H=H_1(E)$
is a free abelian group with
$m$ meridional generators  $t_1,...,t_m$. We have 
   $\bold Z [H]=\bold Z[t_1^{\pm 1},..., t_m^{\pm 1}]$  (the ring of  of Laurent polynomials) and  $Q(H)=
\bold Z (t_1,...,t_m)$ (the ring of  rational functions  on   $t_1,...,t_m$ with 
integer coefficients). Denote by  $\mu$   the inclusion 
$\bold Z [H] \hookrightarrow Q(H)$. 
Choose any  Euler structure $e$
  on
$E$ and set
$$A_e=A_e(t_1,...,t_m) =\tau^\mu (E,e,\omega_L)\in  Q(H)=\bold Z (t_1,...,t_m).$$ It is known (see for instance [Tu3, Sect.\
1.7]) that
$\overline A_e = (-1)^m t_1^{\nu_1}...t_m^{\nu_m} A_e $  
 where  the overbar  denotes the involution in $Q(H)$ sending each $t_i$ to $t_i^{-1}$ and    $(\nu_1,...,\nu_m)\in \bold
Z^m$  is a charge on
$L$ in the sense of Sect.\ 1.6. Then
$$\nabla_L=- t_1^{\nu_1}...t_m^{\nu_m} A_e (t^2_1,...,t^2_m)\in \bold Z (t_1,...,t_m). \eqno (8.1.a)$$
Clearly, $\overline {\nabla}_L=(-1)^m\,\nabla_L$. By (4.2.a), $\nabla_L$ does not depend on the choice of $e$.

We introduce a    version of    $\nabla_L $ depending on   a charge  
   $k=(k_1,...,k_m) \in \bold Z^m$ on $L$.  Set
$$\nabla (L,k) = -t_1^{k_1/2}... t_m^{k_m/2}\, \nabla_L (t_1^{1/2},...,t_m^{1/2}). $$
Note that although the  factors on the right-hand side can lie in 
$\bold Z (t^{1/2}_1,...,t^{1/2}_m)$, their product $\nabla (L,k)$  belongs to
$\bold Z (t_1,...,t_m)$. 
We claim that
$$\nabla (L,k)=\tau(E,e_k,\omega_L).\eqno (8.1.b)$$
  Indeed, 
 the splitting (7.1.a) consists here of only one summand and therefore 
 $\tau(E,e_k,\omega_L)=\tau^\mu (E,e_k,\omega_L)=A_{e_k}$.
By (4.6.a),
$$t_1^{\nu_1}...t_m^{\nu_m}=\pm \overline A_{e_k}/A_{e_k}=\pm c(e_k)^{-1}=\pm t_1^{-k_1}...t_m^{-k_m}$$
so that $\nu_i=-k_i$ for $i=1,...,m$. Substituting this in (8.1.a) with $e=e_k$ we obtain a formula  equivalent to (8.1.b). 

To end this subsection note that   for  $m\geq 2$ both $\nabla_L$ and $\nabla (L,k)$ are Laurent polynomials, i.e., belong to $
\bold Z[t_1^{\pm 1},..., t_m^{\pm 1}]$.  If $m=1$ then 
$\nabla_L$ and $ \nabla (L,k)$ are rational functions on $t=t_1$ which can be computed from the Alexander polynomial  of
$L=L_1$.  Let  $\Delta_L (t) \in 
\bold Z[t^{\pm 1}]$ be the   Alexander polynomial  of $L$  normalized (and symmetrized) in the canonical way  so that 
$\Delta_L (t^{-1})= \Delta_L (t)$ and $\Delta_L(1)=1$. Then $$\nabla_L=
 (t-t^{-1})^{-1} \Delta_L (t^2),\,\,\,\,\,\,\,\nabla (L,k)= -t^{(k+1)/2} (t-1)^{-1} \Delta_L (t) \eqno (8.1.c)$$
where $k=k_1\in \bold Z$ is the charge of $L$.

\skipaline \noindent {\bf   8.2. An invariant of sublinks.} Let $L=L_1\cup ...\cup L_m$ be an $m$-component
ordered  oriented link  
 in an oriented  3-dimensional integral homology sphere $N$. Let $L^I=\{L_i\}_{i\in I}$ be a    sublink of $L$
determined by a proper subset $I\subset \{1,...,m\}$. In this
subsection
we associate with  
$(L,L^I)$  an abelian group $H=H(L,L^I)$. For any charge   
$k $  on $L$ (see Sect.\ 1.6), we define an element 
$\nabla (L,L^I,k) \in Q(H)$. 

Set $\vert I \vert =\card (I)$.  The group $H=H(L,L^I)$
is defined as the quotient of the free (multiplicative) abelian group of rank $m$ with generators $t_1,...,t_m$
modulo the following $\vert I \vert$    relations
numerated by $i\in I$:
$$  \prod_{ j\neq i} (t_j t_i^{-1}) ^{lk(L_i,L_j)}=1   \eqno (8.2.a)$$ 
where  the product goes over all $j\in \{1,...,m\}\backslash \{i\}$. 
 The   ring $\bold Z[H]$ can be identified with the ring of polynomials on $m$ commuting invertible variables $t_1,...,t_m$
subject to
 (8.2.a).  There is a   group homomorphism from $H $ to the
infinite cyclic group sending  
$t_1,...,t_m$ to
the same generator. Therefore  each    $t_j$ has an infinite order in $H$  and  $t_j-1$ is a non-zerodivisor  
in
$\bold Z[H]$   invertible in $Q(H)$.  Set 
$$\nabla (L,L^I,k)=  \left ( {\prod_{i\in I} (t_i-1)^{-1} }\right ) \nabla (L,k)
\in Q(H).
\eqno (8.2.b)$$  
Here $\nabla (L,k)$
is viewed as a rational function  on   $t_1,...,t_m$   representing an element of $Q(H)$.   
In particular, if   $I=\emptyset$, then   $  \nabla
(L,L^I,k)=\nabla (L,k)$.

The group $H=H(L,L^I)$ and the element  $\nabla (L,L^I,k)\in Q(H)$ can be   computed from the linking 
numbers of the components of $L$ and $\nabla (L,k)$.  We give  a geometric interpretation 
of $H $ and   $\nabla (L,L^I,k) $.  
Let $E $ be the exterior of $L$ in $N$.  Let $P=P(L,L^I)$ be the  3-manifold   obtained from  $E$  by gluing  directed  
solid tori $\{\Z_i\}_{i\in I}$ such that the longitude  of $L_i$ with  $i\in I$ determined by the framing 
number $-lk (L_i, L\backslash
L_i)=- \sum_{j\neq i} lk (L_i,L_j)$ bounds a meridional disc 
in $\Z_i$.  The core of $\Z_i$ is oriented as in Sect.\ 1.7. The boundary of $P$ consists of $m-\vert I \vert\geq 1$ tori
corresponding to the components  of $L\backslash L^I$.  Identifying
$t_1,...,t_m$ with    the meridional generators of $H_1(E)$ we obtain an identification 
$H_1(P)=H $.
 We claim that
 $$\nabla (L,L^I,k) =\tau(P, e_k^P, \omega_L^P) \eqno (8.2.c)$$
where $e_k$ is the Euler structure  on $E$ determined by $k$,   $e_k^P= (e_k)^P$ is the induced Euler structure on $P$ (cf.
Sect.\ 1.5--1.7), $\omega_L$ is the    homology orientation  of $E$  determined by the
orientation of $L$, and $\omega_L^P=(\omega_L)^P$ is the 
  induced homology orientation  on $P$ (cf. Sect.\ 3.1, 3.2). 
Indeed, by      (7.3.a) and  (8.1.b)  
$$  \left ( \prod_{i\in I} (t_i-1)\right )  \tau(P, e_k^P, \omega_L^P) =\tau(E,e_k,\omega_L)= \nabla (L,k).$$
This and  (8.2.b) yield  (8.2.c).
This interpretation of  $\nabla
(L,L^I,k)$ implies that if $\rank \, H\geq 2$ then $\nabla (L,L^I,k)\in \bold Z[H]$. 
If $\rank \, H=1$ then $\nabla (L,L^I,k)\in (\bold Z[H])_1$. 
 
\skipaline \skipaline \centerline  {\bf  9.  A surgery formula for  $\varphi$-torsions}

\skipaline \noindent {\bf    9.1.  Setting and notation. }  
Let $M$ be a  closed    3-manifold
obtained by   surgery on a framed oriented link $L=L_1\cup ...\cup L_m$ in an oriented 
 3-dimensional integral homology sphere  $N$. Let $F$ be a field   of characteristic $0$. Fix a charge $k$ of
$L$   and consider  
 a ring homomorphism   $\varphi^M : \bold Z[H_1(M)]\to F$ such that $\varphi^M(H_1(M))\neq 1$. 
We shall give a surgery formula for  the torsion  $\tau^{\varphi^M} (M, e_k^M, \omega_L^M)\in F$
where $e_k^M=(e_k)^M \in \vect (M)$ is   induced by 
the   Euler structure $e_k  $ 
on the exterior $E$ of $L$  and  $\omega_L^M=(\omega_L)^M$ is the homology orientation   of $M$ 
  induced by the homology orientation $\omega_L$ of $E$ (see Sect.\  3.1).

Recall the notation $\overline I= \{1,...,m\}  \backslash I$ for a   set $I\subset  \{1,...,m\} $. 
Denote by   $\ell^{ I}(L)= (\ell^{ I}_{i,j}(L))$  the symmetric  square matrix whose rows and columns
are numerated by  
$i,j\in  I$ and whose entries are given by
$$
\ell^{ I}_{i,j}(L)=\cases lk(L_i,L_j),~ {  {if}}\,\,\, i\neq j, \\
 lk(L_i,L_i) + lk(L_i, L^{\overline I})= lk(L_i,L_i)+ \sum_{l\in  \overline I} \, lk(L_i,L_l),~ {  {if}} \,\,\, i=j .\endcases$$
 Here $lk(L_i,L_i)$ is the framing number of $L_i$.

Denote by $\varphi$ the composition of $\varphi^M$
with the  ring homomorphism $\bold Z[H_1(E)]\to \bold Z[H_1(M)]$ induced by the inclusion $E\subset M$. 
Let $t_1,...,t_m\in H_1(E)$ be the
homology classes of the meridians of $L_1,...,L_m$. 
Set
$$I(\varphi)= \{i  \,\vert \, \varphi (t_i)= 1\} \subset  \{1,...,m\} .$$
For any   $i\in  {I(\varphi)}$, the homomorphism $\varphi$ maps   the meridian $t_i$ and the longitude of
$L_i$ to $1\in F$.  This implies that for a    set $I\subset  \{1,...,m\} $, not contained in  $I(\varphi)$,   the formula
$\{t_i\mapsto
\varphi (t_i)\in F\}_{i\in I}$ 
defines  a ring homomorphism $
\bold Z[H(L^I,  L^{I\cap     {I(\varphi)}})]\to  F$ (cf. Sect.\ 8.2). It is  denoted $\varphi^I$. 
Recall the notation $\vert I \vert=\card (I)$.

\skipaline \noindent {\bf    9.2.  Theorem. } 
 {\sl  Under the assumptions of Sect.\ 9.1,}
$$ \tau^{\varphi^M} (M, e_k^M, \omega_L^M) \eqno (9.2.a)$$
$$=\left (\prod_{i\in \overline {I(\varphi)}} (\varphi(t_i)-1)^{-1}\right )
\sum_{  { I\subset  I(\varphi)}} (-1)^{\vert I \vert} \,\det (\ell^{ I}(L))\,
\varphi^{\overline I}_\#  (\nabla (L^{\overline I},L^{\overline I
\cap   I(\varphi)},k^{\overline I}))
 . $$

\skipaline  Formula (9.2.a) allows to compute the $\varphi$-torsion  of $M$ 
   in terms of the Alexander-Conway polynomials of $L$ and its  sublinks, and the linking and framing numbers of the components of
$L$. A curious feature of  this formula is a framing-charge separation: the factor  
$\det (\ell^{ I})$ does not depend on the charge  while  the factor  
$\varphi^{\overline I}_\#  (\nabla (L^{\overline I},L^{\overline I
\cap   I(\varphi)},k^{\overline I}))$ does not depend on the framing.  The product 
$\prod_{i\in \overline {I(\varphi)}} (\varphi(t_i)-1)^{-1}$ depends neither on the  charge  nor on the framing.

The terms   on the right-hand side of (9.2.a) corresponding to   proper subsets  $I\subset I(\varphi)$ can be  slightly simplified.  
Observe that  $$\rank \,  H (L^{\overline I},L^{\overline I
\cap   I(\varphi)}) \geq  \vert\overline {I} \backslash I (\varphi)\vert= \vert\overline {I(\varphi)}\vert + \vert{I(\varphi)}
\backslash I\vert.$$ The assumption $\varphi (H_1(E))\neq 1$ implies that $\vert \overline {I(\varphi)}\vert \geq 1$. Therefore
for  a proper subset  $I\subset I(\varphi)$ we have $\rank \,  H (L^{\overline I},L^{\overline I
\cap   I(\varphi)}) \geq  2$ and
$$\varphi^{\overline I}_\#  (\nabla (L^{\overline I},L^{\overline I
\cap   I(\varphi)},k^{\overline I}))=\varphi^{\overline I}   (\nabla (L^{\overline I},L^{\overline I
\cap   I(\varphi)},k^{\overline I})).$$ The term    corresponding to $I=I(\varphi)$ equals
$$\prod_{i\in \overline {I(\varphi)}} (\varphi(t_i)-1)^{-1}
  (-1)^{\vert I(\varphi)\vert} \,\det (\ell^{ I(\varphi) }(L))\,
\varphi^{\overline {I(\varphi)}}_\#  (\nabla (L^{\overline {I(\varphi)}},k^{\overline {I(\varphi)}})).$$
We can replace here $\varphi^{\overline {I(\varphi)}}_\# $ by $\varphi^{\overline {I(\varphi)}}$ provided 
$\vert \overline {I(\varphi)}\vert
\geq 2$.  

 If $I(\varphi)=\emptyset$, i.e.,  if  $\varphi (t_i)\neq 1$ for all $i=1,...,m$, then   (9.2.a) simplifies to  
$$ \tau^{\varphi^M} (M, e_k^M, \omega_L^M) 
 =\prod_{i=1}^m (\varphi(t_i)-1)^{-1}
\, \varphi_\#  (\nabla (L, k ))
.  \eqno (9.2.b)$$ This can be obtained   by applying (5.1.a) $m$ times followed by  Lemma 7.3.2 and (8.1.b):
$$ \tau^{\varphi^M} (M, e_k^M, \omega_L^M) 
 =\prod_{i=1}^m (\varphi(t_i)-1)^{-1}
\,
\tau^\varphi (E,e_k,\omega_L)$$
$$=
\prod_{i=1}^m (\varphi(t_i)-1)^{-1} \, \varphi_\#  (\tau(E,e_k,\omega_L))
=\prod_{i=1}^m (\varphi(t_i)-1)^{-1} \, \varphi_\#  (
\nabla
(L,k)).$$

\skipaline   {\sl Proof of  Theorem 9.2.} The proof goes by induction on $m$.  If $m=1$ then $I(\varphi)=\emptyset$  and  
(9.2.a) reduces to (9.2.b) proven above. Assume that for links with $<m$ components, Theorem 9.2 is true and  prove  it for
a link
$L$ with $m$ components.  

Consider first the case where for any 
$i\in  { I(\varphi)}$ the framing number  $lk(L_i,L_i)$   is equal to 
$- \sum_{l\neq i} \,  lk(L_i,L_l)$ (the framing numbers of $\{L_i\,\vert \,  {i\in \overline { I(\varphi)}}\}$ may be arbitrary). 
 Then for any   non-void  set
$I\subset  I(\varphi)$, the sum of  the columns of  the matrix $\ell^{  I}$ is     $0$ so that $\det (\ell^{ 
I}(L))=0$. Thus, the right-hand side of (9.2.a) contains only one possibly non-zero term corresponding to $I=\emptyset$. 
By (8.2.c), this term equals 
$$\prod_{i\in \overline {I(\varphi)}} (\varphi(t_i)-1)^{-1}
 \, \psi_\#  (\nabla (L,L^{    {I(\varphi)}},k))=
\prod_{i\in \overline {I(\varphi)}} (\varphi(t_i)-1)^{-1}
 \, \psi_\#  (\tau(P,e_k^P,\omega_L^P))
$$
where $P=P (L,L^{   {I(\varphi)}})$  is the 3-manifold  defined as in Sect.\ 8.2 and $\psi=\varphi^{\overline \emptyset}$ is the
ring homomorphism
$\bold Z[H_1(P)]=\bold Z[H(L,L^{    {I(\varphi)}})]
\to F$ induced by $\varphi$. By Lemma 7.3.2,
$$\prod_{i\in \overline {I(\varphi)}} (\varphi(t_i)-1)^{-1}
 \, \psi_\#  (\tau(P,e_k^P,\omega_L^P))= \prod_{i\in \overline {I(\varphi)}} (\varphi(t_i)-1)^{-1}
 \,
\tau^{\psi} (P,e_k^P,\omega_L^P).$$ 
The manifold  $M$ can be  obtained from $P$ by gluing directed solid tori corresponding to the components of
$L^{\overline {I(\varphi)}}$. Applying  (5.1.a) inductively and using Lemma 3.3, we obtain  
$$
\prod_{i\in \overline {I(\varphi)}} (\varphi(t_i)-1)^{-1}
 \,
\tau^{\psi} (P,e_k^P,\omega_L^P)= \tau^{\varphi^M} (M, (e_k^P)^M, (\omega_L^P)^M)=
\tau^{\varphi^M} (M, e_k^M,\omega_L^M).$$
This yields (9.2.a) in this case.

 Up to the rest of the proof, we denote the left  (resp. right)
hand side of (9.2.a)
by $\langle L, \varphi, k\rangle_1$ (resp.  $\langle L, \varphi, k\rangle_2$). We claim that
 for a framed link $L'$ obtained from $L$ 
by one   negative twist of the framing of   a component of $L^{  { I(\varphi)}} \subset L$,
$$\langle L, \varphi, k\rangle_1-\langle L, \varphi, k\rangle_2=\langle L', \varphi, k\rangle_1-\langle L', \varphi, k\rangle_2.
\eqno (9.2.c)$$ 
Hence,   formula (9.2.a) holds for $L$ if and only if  it holds  for $L'$. Together with the result of the previous paragraph this will
imply the claim of the theorem.

Now we check (9.2.c). Assume for concreteness that 
$m\in {  { I(\varphi)}}$ and that 
$L'=L_1\cup  ...\cup L_{m-1} \cup L'_m$ is   obtained from $L$ by one   negative twist of the framing of
$L_m$. Consider the closed 3-manifold, $W$,  obtained by   surgery on $N$ along
the framed link $L^{\overline m}= L_1\cup  ...\cup L_{m-1}$. The framed oriented 
knots
$L_m, L'_m$ lie in
$W$ and differ only by the   framing.  These knots have
  the same exterior in $W$, which we denote by $V$.  It is clear that  $\partial V$ is a 2-torus bounding
a regular neighborhood of $L_m$ in $W$ which can be identified with a regular neighborhood of $L_m$ in $N$. 
Let $\alpha_1\subset
 \partial V$ (resp. $\alpha_2\subset
 \partial V$) be a longitude of $L_m$ (resp.  of $L'_m$) determined by the framing.  Let 
$\alpha_3\subset
 \partial V$ be a meridian of $L_m$. By assumption, $\alpha_1$ is homological to $\alpha_2\alpha_3$ in $\partial V$. 
In the notation of Sect.\ 6.1 we have $M=M_{V, \alpha_1, -1}$
(the sign $-1$ appears because the orientation of $M$ induced from the one in $N$ is  {\sl opposite} to the product
orientation in the directed solid torus $M\backslash V$; see Sect.\ 1.7 and 6.1 for  the distinguished orientations
of the core and meridional disc).  Set
$M'=M_{V,
\alpha_2, -1}$. Observe that
$M_{V,
\alpha_3, +1}=W$. As in Sect.\ 1.7, the manifold $V$ is obtained from  the exterior  $E$ of 
$L $ by gluing $m-1$ directed solid tori.  Set   $e=(e_{k})^V\in \vect (V)$ and  provide $V$ with   homology orientation
$\omega =(\omega_{L})^V $.  Denote by ${\varphi^V}$ the ring homomorphism
$\bold Z[H_1(V)]\to F$ induced by $\varphi$; this   homomorphism   is the composition of the inclusion
homomorphism
$\bold Z[H_1(V)]\to \bold Z[H_1(M)]$ with $\varphi^M$.   By Lemma 6.4,   
$$
  \tau^{(\varphi^V)^{M}}
(M,e^{M},\omega^{M}) 
= \tau^{(\varphi^V)^{M'}} (M',e^{M'},\omega^{M'})
-  \tau^{(\varphi^V)^{W}} (W,e^{W},\omega^{W})
. \eqno (9.2.d)$$
By definition of $ \varphi^V$, we have $(\varphi^V)^M=\varphi^M$. By definition of $e_k^M$, we have $e^M=e_k^M$.
By Lemma 3.3, $\omega^{M}=((\omega_{L})^V)^M=\omega_L^M$. 
Thus the left-hand side of (9.2.d) equals $ \langle L, \varphi, k\rangle_1 $. 
The manifold $M'$ is  obtained by   surgery on  $N$ along  $L'$ and  similarly
$$\tau^{(\varphi^V)^{M'}} (M',e^{M'},\omega^{M'})= \langle L', \varphi, k\rangle_1. $$ Consider  now the exterior
$E^{\overline m}$ of  $L^{\overline m} =L_1\cup  ...\cup L_{m-1}$ in
$N$,  the ring homomorphism
$\psi: \bold Z[H_1(E^{\overline m})] \to F$ induced by $\varphi$, the charge $k^{\overline m}$ on
$L^{\overline m}$  induced by $k$ as in
Sect.\ 1.6 
and  the homology orientation  $\omega_{L^{\overline m}}$ of $E^{\overline m}$ nduced by the orientation of $L^{\overline m}
$. We claim that
$$\tau^{(\varphi^V)^{W}} (W,e^{W},\omega^{W})=
  \langle L^{\overline m}, \psi, k^{\overline
m}\rangle_1.\eqno (9.2.e)$$ Indeed, 
the manifold $W$ is  obtained by   surgery on $N$ along $L^{\overline m}$ and hence by gluing directed solid tori to 
$E^{\overline m}$. 
It follows from definitions that  $(\varphi^V)^{W}
=\psi^W$  (both homomorphisms map a meridional generator
$t_i\in H_1(W)$ with $i=1,...,m-1$ to $\varphi(t_i)$).  The remarks of Sect.\ 1.6  imply that
$(e_{k^{\overline m}})^W=((e_{k})^V)^W=e^{W}$.  Lemmas  3.3 and 3.4 imply that 
$(\omega_{L^{\overline m}})^W= ((\omega_{L})^V)^W=\omega^{W}$. Hence
$$\tau^{(\varphi^V)^{W}} (W,e^{W},\omega^{W})=
 \tau^{\psi^W} (W, (e_{k^{\overline m}})^W, (\omega_{L^{\overline m}})^W)= \langle L^{\overline m}, \psi, k^{\overline
m}\rangle_1. $$ 
Now formula (9.2.d) can be rewritten as
$$\langle L, \varphi, k\rangle_1-\langle L', \varphi, k\rangle_1= -  \langle L^{\overline m}, \psi, k^{\overline m}\rangle_1. \eqno
(9.2.e)$$ We claim that
similarly $$\langle L, \varphi, k\rangle_2-\langle L', \varphi, k\rangle_2=  - \langle L^{\overline m}, \psi, k^{\overline
m}\rangle_2.\eqno (9.2.f)$$
Indeed, it is clear that $I(\psi)=I(\varphi)\backslash \{m\}$ so that the complement of $I(\varphi)$ in $\{1,...,m\}$ coincides
with the complement of $I(\psi)$ in $\{1,...,m-1\}$. Hence the factors 
$
\prod_{i\in \overline {I(\varphi)}} (\varphi(t_i)-1)^{-1}$ which appear in the  three
terms  of (9.2.f) are identically the same.  For a set 
$  I \subset  I(\varphi) $ not containing $ m$, the corresponding summands in 
$\langle L, \varphi, k\rangle_2, \langle L', \varphi, k\rangle_2$ do not depend on the framing of the $m$-th component and
therefore are equal.   Consider  the summand  on the right-hand side of (9.2.a)  corresponding to a set 
$  I \subset  I(\varphi) $   containing $ m$.  The
matrices $ \ell^{ I} (L)$ and $ \ell^{ I} (L')$ coincide except at the $(m,m)$-entry
where  $\ell^{  I}_{m,m} (L)=\ell^{  I}_{m,m} (L')+1 $.
It is easy to see that  $$\det (\ell^{  I} (L))-\det (\ell^{  I} (L'))=\det (\ell^{ J} (L^{\overline m}) )$$
where    $J=J_I= I \backslash \{m\} \subset   I(\psi)$. The determinant $\det ( \ell^{J} (L^{\overline m}))$
is exactly the determinant appearing in $\langle L^{\overline m}, \psi, k^{\overline
m}\rangle_2$ in the summand corresponding to $J$. 
The factor 
$\varphi^{\overline I}_\#  (\nabla (L^{\overline I},L^{\overline I
\cap   I(\varphi)},k^{\overline I}))$ is the same for both links $L,L' $ and is equal to 
$\psi^{\tilde  J}_\#  (\nabla ((L^{\overline m})^{\tilde  {J}},(L^{\overline m})^{\tilde  {J} \cap I(\psi)},(k^{\overline
m})^{\tilde  {J}}))$ where   $\tilde  J=\{1,...,m-1\}\backslash J=\overline I$. Clearly,  $(-1)^{\vert I \vert}= - (-1)^{\vert J
\vert}$.  Since the formula $I \mapsto J_I= I \backslash \{m\}$ establishes a bijective correspondence between    sets 
$I\subset I(\varphi)$ containing $m$ and   subsets of $I(\psi)$, these equalities  yield (9.2.f).  Formulas  (9.2.e),   (9.2.f) and   the
induction assumption 
$\langle L^{\overline m}, \psi, k^{\overline m}\rangle_1 =   \langle L^{\overline m}, \psi, k^{\overline
m}\rangle_2$ imply (9.2.c). This accomplishes the proof of the theorem.

\skipaline \skipaline \centerline  {\bf  10.  A surgery formula for the Alexander polynomial}

\skipaline \noindent {\bf    10.1.  The refined Alexander polynomial. }  Let $M$ be a closed  connected 3-manifold  
with $b_1(M) \geq 1$.  Set $G=H_1(M)/\Tors H_1(M)$. 
The  Alexander polynomial $\Delta (M)$ of  $\pi_1(M)$ is defined   using  a 
finite presentation  of  $\pi_1(M)$ and the Fox differential calculus (see [CF]).  This polynomial is an element of the group ring
$\bold Z [G]$ defined up to multiplication by
$\pm  G$.  The  theory of torsions allows to  introduce  a refinement   $\Delta (M,e,\omega)  $ of  $\Delta (M)$ 
     depending on an Euler structure  $e\in \vect (M)$ and a homology orientation 
$\omega$ of $M$.  We define   $\Delta (M,e,\omega)  $
as follows. Denote by $Q(G)$ the field of fractions of the domain $\bold Z[G]$. Denote by $\mu$
the composition of the projection $\bold Z [H_1(M)] \to \bold Z [G]$ with the inclusion 
$\bold Z [G] \hookrightarrow Q(G)$. Set
$ \Delta (M,e,\omega) =\tau^\mu (M,e,\omega) \in Q(G)$. 
It follows from [Tu1] that  $ \Delta (M,e,\omega)$  is   a refinement of   $\Delta (M)$  in the following sense.  If
$b_1(M)
\geq 2$ then
$\Delta (M,e,\omega) \in \bold Z [G]$ and
$\Delta (M) =\pm G \, \Delta (M,e,\omega)$. If $b_1(M)=1$ and $M$ is orientable then 
$ \Delta (M,e,\omega) \in (t-1)^{-2} \bold Z [G]$ for any generator $t$ of   $G$ and  $\Delta (M) =\pm G
(t-1)^{2} \Delta (M,e,\omega)$. A similar formula holds for non-orientable $M$ (see [Tu1]) but we shall not need it.

\skipaline \noindent {\bf    10.2.  Theorem (a surgery formula for  $\Delta (M,e,\omega)$). } 
 {\sl  Let $M$ be a  closed    3-manifold with $b_1(M) \geq 1$ 
obtained by   surgery on a framed oriented link $L=L_1\cup ...\cup L_m$ in an oriented 
 3-dimensional integral homology sphere  $N$. Let   $k=(k_1,...,k_m)$ be a charge on $L$. 
Let} $[[t_1]],...,[[t_m]]\in G= H_1(M)/\Tors H_1(M)$  {\sl  be  represented by the meridians of $L_1,...,L_m$, respectively. 
Set
$$I_0= \{i  \,\vert \, [[t_i]] =1 \} \subset  \{1,...,m\} . \eqno (10.2.a)$$ Then }
$$ \Delta (M, e_k^M, \omega_L^M) \eqno (10.2.b)$$
$$=\prod_{i\in \overline {I_0}} ([[t_i]] -1)^{-1}
\sum_{  { I\subset  I_0}} (-1)^{\vert  I \vert} \,\det (\ell^{ I}(L))\,
\mu^{\overline I}_\#   (\nabla (L^{\overline I},L^{\overline I
\cap   I_0},k^{\overline I}))
 $$
{\sl where 
$\mu^{\overline I}_\# : (\bold Z[H(L^{\overline I},  L^{{\overline I}\cap     {I_0}})])_1 \to Q(G)$ is the natural extension of the
ring homomorphism 
$\mu^{\overline I} : \bold
Z[H(L^{\overline I},  L^{{\overline I}\cap     {I_0}})] \to \bold Z [G]$   sending the generator of
$H(L^{\overline I},  L^{{\overline I}\cap     {I_0}})$ represented by the meridian of $L_i$ to $[[t_i]]\in G$,  for  all $i \in
\overline I$.}

\skipaline

Formula (10.2.b) is obtained by a direct application of Theorem 9.2 to
the ring homomorphism  $\mu $ defined in Sect.\ 10.1.     Note that   all $[[t_i]]-1\in \bold Z[G]$ with $i\in \overline {I_0}$
are  non-zero and therefore
invertible in  
$Q(G)$. 
We have  $$\rank \, H(L^{\overline I},  L^{{\overline I}\cap     {I_0}})\geq \vert \overline I \backslash
({\overline I}\cap     {I_0})\vert=   m- \vert {I_0}\vert   \geq b_1(M)\geq 1. \eqno (10.2.c)$$ Hence
$\nabla (L^{\overline I},L^{\overline I
\cap   I_0},k^{\overline I})\in (\bold Z[H(L^{\overline I},  L^{{\overline I}\cap     {I_0}})])_1$ and   (10.2.b) makes sense.

If
  $m- \vert  {I_0}\vert  \geq  2$,  then  
$\nabla (L^{\overline I},L^{\overline I
\cap   I_0},k^{\overline I})\in  \bold Z[H(L^{\overline I},  L^{{\overline I}\cap     {I_0}})]$, and   (10.2.b) can be rewritten 
in a slightly simpler form:
$$ \Delta (M, e_k^M, \omega_L^M) \eqno (10.2.d)$$
$$=\prod_{i\in \overline {I_0}} ([[t_i]] -1)^{-1}
\sum_{  { I\subset  I_0}} (-1)^{\vert  I \vert} \,\det (\ell^{ I}(L))\,
\mu^{\overline I}   (\nabla (L^{\overline I},L^{\overline I
\cap   I_0},k^{\overline I})).
 $$
This  applies in particular  if $b_1(M) \geq 2$.

If   $m- \vert  {I_0}\vert   =1$ then it may   happen that $\rank \, H(L^{\overline I},  L^{{\overline I}\cap     {I_0}})=1$. In this
case we do need  
$\mu^{\overline I}_\# $. We shall see now that this  can happen  
  for only one  entry, namely the one corresponding to $I=I_0$. We begin with a  lemma which will be proven in Sect.\ 10.5.

\skipaline {\bf 10.3. Lemma.} {\sl   Assume that $I_0$ consists of all $i=1,...,m$ except a certain $n$. Then the  framing number 
of
$L_n$ is
$0$  and 
$lk(L_n, L_i)=0$ for all
$i\neq n$.   For sets  $I\subset J\subset I_0$, we have}   $\rank \, H(L^{\overline
I},L^{\overline I
\cap   J}) =1$ {\sl   if and only if  $I=J=I_0$.  }

\skipaline Using this lemma   we can  reformulate  (10.2.b) in the case  
   $ I_0=\{1,...,m\}\backslash \{n\}$.  Set $t= [[t_n]] \in G$. By assumption, $t\neq 1$ and therefore $t -1$ is invertible in
$Q(G)$. It follows from Theorem 10.2, Lemma 10.3,  and (8.1.c) that 
$$ \Delta (M, e_k^M, \omega_L^M)   
=  (t -1)^{-1}
\sum_{  { I\subset  I_0, I\neq I_0}} (-1)^{\vert  I \vert} \,\det (\ell^{ I}(L))\,
\mu^{\overline I}   (\nabla (L^{\overline I},L^{\overline I
\cap   I_0},k^{\overline I}))
 $$
$$+ \
  (-1)^{    m} \,  \det (\ell^{ I_0}(L)) \,\,
t^{(k_n+1)/2}\,(t -1)^{-2}\, 
\Delta_{L_n} (t) \eqno (10.2.e)$$
where $\Delta_{L_n}$ is the canonically normalized Alexander polynomial of  $L_n$ (see Sect.\ 8.1) 
and $\Delta_{L_n} (t)\in \bold Z [G]$ is obtained by computing $\Delta_{L_n}$ on $t$. 
   By definition of
a charge, the integer 
$k_n+1=k_n+1- lk (L_n,  L
\backslash L_n)$ is even.

  \skipaline {\bf  10.4. Special cases.}  
We describe   two  cases where (10.2.b)   simplifies.
If  $m\geq 2$ and all
$[[t_1]],..., [[t_m]]$  are non-trivial elements of $G$ then $I_0=\emptyset$ and the sum on the right-hand side of (10.2.b) has
only one term. We obtain  
$$\Delta (M,e_k^M,
\omega_L^M)= 
\prod_{i=1}^m ([[t_i ]]-1)^{-1} \,\, [[\nabla (L,k)]]
 $$ where  
for a Laurent polynomial $\nabla\in  \bold Z [t_1^{\pm 1},..., t_m^{\pm 1}]$ we denote by $[[\nabla]]$ 
 the element of $\bold Z [G]$  obtained from $\nabla  $ by the substitution  $t_i\mapsto [[t_i]] \in
G$ for $i=1,...,m$.

Consider now the case where $L=L_1\cup ... \cup L_m \subset N$  is an algebraically split link. 
This means that  $lk(L_i,L_j)=0$ for all $i\neq j$ where $lk$ is the linking number in $N$. 
A charge $k=(k_1,...,k_m)$ on $L$ is just an $m$-tuple of odd integers. It follows easily from the 
Torres formula (5.1.a) and definitions  
that    the Laurent polynomial $\nabla (L,k) $ is   divisible by $ \prod_{i=1}^m (t_i-1)$ in 
$\bold Z [t_1^{\pm 1},..., t_m^{\pm
1}]$.  Set
$$\check \nabla (L,k)= \nabla (L,k)/\prod_{i=1}^m (t_i-1)\in 
 \bold Z [t_1^{\pm 1},..., t_m^{\pm
1}].$$
Let  $M$ be a 3-manifold obtained by  surgery along $L$ endowed with certain framing numbers $ f_1,...,f_m \in \bold Z$. 
The set $I_0$ defined by (10.2.a) is computed by $I_0=\{i=1,...,m\,\vert \, f_i\neq 0\}$. 
Clearly, $b_1(M)=m- \vert  I_0 \vert$. If
$b_1(M) \geq 2$  then (10.2.d) implies
$$\Delta(M,e_k^M, \omega_L^M)= \sum_{I\subset I_0}\,\, (-1)^{\vert I \vert}\, \prod_{i\in I}   f_i\,\, [[ \check \nabla
(L^{\overline I} ,k^{\overline I})]] .  $$
If $b_1(M)=1$  then $I_0$ consists of all $i=1,...,m$ except a certain $n$ and (10.2.e) implies
$$ \Delta (M, e_k^M, \omega_L^M)   
=   
\sum_{  { I\subset  I_0, I\neq I_0}} (-1)^{\vert  I \vert} \, \prod_{i\in I}   f_i\,\, [[ \check \nabla
(L^{\overline I} ,k^{\overline I})]]
 $$
$$+ \
  (-1)^{    m} \, \big ( \prod_{i\neq n}   f_i \big ) \,\,
t^{(k_n+1)/2}\,(t -1)^{-2}\, 
\Delta_{L_n} (t) $$
where $t=[[t_n]] \in G$.
In particular,  if $L$ is a knot with framing number $0$ and  charge $k\in \bold Z$ then 
 $\Delta (M,e_k^M,
\omega_L^M) 
= -t^{(k+1)/2}\,(t-1)^{-2}\,  \Delta_L (t)$.

  \skipaline {\bf  10.5. Proof of Lemma 10.3}.  By assumption,   $[t_i] \in \Tors H_1(M)$ for all  $i\neq n$ and $[t_n]\in H_1(M)$ is an
element of infinite order.   Therefore the  group $H_1(M)/\{[t_i]=1\,\vert \,i\neq n\}$ is   infinite.  Since $(lk (L_i,L_j))_{i,j}$ is a
presentation matrix for $H_1(M)$ with respect to the generators $[t_1],..., [t_m]$,  this implies that
$lk(L_n, L_i)=0$ for all
$i$. If 
$I=J=I_0$  then  $H(L^{\overline I},L^{\overline I
\cap   J})= H(L_n,\emptyset)=\bold Z$. Conversely, assume that 
 $I\subset J\subset I_0$  and  $\rank \, H(L^{\overline
I},L^{\overline I
\cap   J})=1$.  The group $ H(L^{\overline I},L^{\overline I
\cap   J}) $ is generated by $\{t_i\}_{i\in \overline I}$ modulo the relations 
$$ \prod_{i\in \overline I} (t_i t_r^{-1})^{lk(L_i,L_r)}=1 $$ 
numerated by $r\in   \overline I
\cap   J $.  By the first part of the claim,  these relations do not involve the generator $t_n$.  Hence
$H(L^{\overline I},L^{\overline I
\cap   J})$ splits as a product of the infinite cyclic group generated by $t_n$ and a group $H'$ generated by 
$t_i$ with $i\in \overline I \backslash \{n\}=I_0\backslash I$. If $I\neq I_0$ then there is an epimorphism
$H'\to \bold Z$ sending all $t_i$ with $i\in I_0\backslash I$ to $1$. Then $\rank \, H(L^{\overline I},L^{\overline I
\cap   J})=\rank \,  H' +1\geq 2$ which contradicts our
assumptions. Thus,
$I=I_0$ and therefore also $J=I_0$.

\skipaline \skipaline \centerline  {\bf  11.  A surgery formula for  $\tau(M)$  in the case $b_1(M)\geq 1$}

\skipaline \noindent {\bf    11.1.  Setting and notation. } 
  Let $M$ be a  closed    3-manifold  
obtained by   surgery on a framed oriented link $L=L_1\cup ...\cup L_m$ in an oriented 
 3-dimensional integral homology sphere  $N$. 
 Let   $k=(k_1,...,k_m)$ be a charge on $L$.  Set $H=H_1(M)$. We shall need  the following notation. 

 For any   $a\in Q(H)$ we define the {\sl reduced inverse} $a^{-1}_{\red}$ as follows. Recall that $Q(H)$ splits
as a direct sum of fields, $Q(H)=\oplus_r F_r$. There is a unique expansion $a=\sum_r a_r$ with $a_r\in F_r$ for all $r$.
Set $a^{-1}_{\red}=\sum_{r, a_r\neq 0} a_r^{-1}$ where $a_r^{-1}\in F_r$ is the inverse of $a_r$ in $F_r$. 
Clearly, if $a$ is invertible in $Q(H)$ then $a^{-1}_{\red}=a^{-1}$. The reduced inverse is
defined for all elements of $Q(H)$. For instance, $0^{-1}_{\red}=0$. 
If $t$ is an element of $H$ of  infinite  order  then $(t-1)^{-1}_{\red}=(t-1)^{-1}$.  If $t$ is an element of $H$ of  finite  order
$n\geq 1$ then a simple computation shows that 
$$(t-1)^{-1}_\red=\frac {1+2t+3t^2+...+nt^{n-1}}{n} -  \frac {n+1}{2} \cdot \frac {1+t+...+t^{n-1}}{n} 
$$
$$=\frac {(1-n)+(3-n)t+(5-n)t^2+...+(n-1)t^{n-1}}{2n}.$$

Let  $[t_1],...,[t_m]\in H=H_1(M) $   be  represented by the meridians of $L_1,...,L_m$, respectively. 
Set
$$I_0= \{i  \,\vert \, [t_i] \in \Tors H \} \subset  \{1,...,m\} $$
(this is equivalent to (10.2.a)).  For any   sets $I\subset J\subset I_0$ 
we define an additive homomorphism  from $\bold Z[H(L^{\overline I}, L^{\overline I \cap J} )]$ to  $\bold
Q [H]$ called   {\sl transfer}.  The image of   $a\in \bold Z[H(L^{\overline I}, L^{\overline I \cap J} )]$ under the transfer  is
denoted  
$a^{\tr}$. By additivity, it suffices to define $a^{\tr}$ for $a \in H(L^{\overline I}, L^{\overline I \cap J} )$. Consider the group
$H_J=H/\{[t_j]=1\,\vert \, j\in J\}$.  Let $p$ be the  projection $H\to H_J$.   
  Sending each $t_i$ with $i\in \overline I$  to $p([t_i])\in H_J$, we obtain a group homomorphism,
$q:H(L^{\overline I}, L^{\overline I \cap J} )\to H_J$. 
For   $a \in H(L^{\overline I}, L^{\overline I \cap J} )$ we set
$a^{\tr}= \vert \Ker\, p \vert^{-1}\sum_{h\in p^{-1}  q (a)} h \in \bold Q[H]$
where $\vert \Ker\, p \vert =\card (\Ker\, p)$ and the addition on the right-hand side is  the one in the group ring $\bold
Q[H]$.  The sum on the right-hand side is finite: the inclusion $J\subset I_0$ ensures that
 $\Ker\, p$ is finite. The transfer $\bold Z[H(L^{\overline I}, L^{\overline I \cap J} )] \to \bold
Q [H]$ is multiplicative  but in general not a ring homomorphism since $1^\tr=
\vert \Ker\, p \vert^{-1}\sum_{h\in \Ker\, p } h$. 

   We give now a surgery formula for  $\tau (M,e_k^M,\omega_L^M)$ in the case   $b_1(M)\geq 1$,
see Appendix  1 for the case $b_1(M)=0$.

\skipaline \noindent {\bf    11.2.  Theorem. } 
 {\sl  If  }  $m- \vert I_0\vert \geq 2$, {\sl then }
$$ \tau (M, e_k^M, \omega_L^M) \eqno (11.2.a)$$
$$=
\sum_{  { I\subset J\subset  I_0}} (-1)^{\vert I \vert} \,\det (\ell^{ I}(L))\, \prod_{i\in \overline {J}} ([t_i] -1)_\red^{-1}
 \,  (\nabla (L^{\overline I},L^{\overline I
\cap   J},k^{\overline I}))^\tr. 
 $$
 {\sl  If  $I_0$ consists of all $i=1,...,m$ except a certain $n$, then }
$$ \tau (M, e_k^M, \omega_L^M) \eqno (11.2.b)$$
$$=
\sum_{  { I\subset J\subset  I_0}, I\neq I_0 } (-1)^{\vert I \vert} \,\det (\ell^{ I}(L))\, \prod_{i\in \overline {J}} ([t_i]
-1)_\red^{-1}
 \,  (\nabla (L^{\overline I},L^{\overline I
\cap   J},k^{\overline I}))^\tr
$$
$$+
  (-1)^{   m} \, \frac {\det (\ell^{ I_0}(L))}{ \vert \Tors H\vert }  \big ({\sum_{h\in \Tors H} h}\big )    
[t_{n}]^{(k_n+1)/2}\,([t_{n}] -1)^{-2}\, 
\Delta_{L_n} ([t_n]) 
 $$ 
{\sl where $\Delta_{L_n}$ is the canonically normalized Alexander polynomial of  $L_n$ (see Sect.\ 8.1) 
and $\Delta_{L_n} ([t_n])\in \bold Z [H]$ is obtained by computing $\Delta_{L_n}$ on $[t_n]$. }

\skipaline {\sl Proof.}  
Consider  first the case 
$m-\vert I_0 
\vert\geq 2$.  The sum  on the right-hand side of  (11.2.a)   goes over all subsets $I,J$  of $I_0$ such that $I\subset J$.  By
(10.2.c), we have 
$\rank\, H(L^{\overline I},L^{\overline I
\cap   J}) \geq 2$  so that $\nabla (L^{\overline I},L^{\overline I
\cap   J},k^{\overline I}) \in \bold Z[H(L^{\overline I},L^{\overline I
\cap   J}) ]$ and its transfer $ (\nabla (L^{\overline I},L^{\overline I
\cap   J},k^{\overline I}))^\tr$  is a well defined element of $\bold Q[H]\subset Q(H)$.  

Denote by   $\varphi_r$ the
composition of the inclusion 
$\bold Z[H] \hookrightarrow Q(H)$ and the projection $Q(H) \to F_r$ on the $r$-th term in   the splitting   
$Q(H)=\oplus_r F_r$. Set $$I_r=  I(\varphi_r)= \{i  \,\vert \, \varphi_r([t_i])= 1\} \subset  \{1,...,m\} $$ 
and  $\tau_r=
\tau^{\varphi_r} (M, e_k^M, \omega_L^M)\in F_r$.  Note that $I_r\subset I_0$ for all $r$ and therefore
$\vert \overline I_r \vert  \geq \vert \overline I_0\vert \geq 2$.  For
$I\subset  J \subset I_0$, set $\varphi_r^{\overline I}=(\varphi_r)^{\overline I}:\bold Z[H(L^{\overline I},L^{\overline I
\cap   J})] \to F_r$  (cf.  Sect.\ 9.1). We have
$$\tau (M, e_k^M, \omega_L^M) =\sum_r \tau_r
$$
$$=\sum_{  { r}} \sum_{  { I\subset   I_r }}\,\,  \prod_{i\in \overline {I_r}} ( \varphi_r([t_i]) -1)^{-1}
 (-1)^{\vert I \vert} \,\det (\ell^{ I}(L))\,
\varphi_r^{\overline I}(\nabla (L^{\overline I},L^{\overline I
\cap   I_r},k^{\overline I}))$$
$$= \sum_{  { I\subset  J \subset I_0 }}\,  (-1)^{\vert I \vert} \,\det (\ell^{ I}(L))\,\prod_{i\in \overline {J}} (
[t_i]-1)_\red^{-1}\, 
  \sum_{  { r,  J=I_r}} 
\varphi_r^{\overline I}(\nabla (L^{\overline I},L^{\overline I
\cap   J},k^{\overline I}))$$
$$= \sum_{  { I\subset  J \subset I_0 }}\,  (-1)^{\vert I \vert} \,\det (\ell^{ I}(L))\,
\prod_{i\in
\overline {J}} ( [t_i] -1)_\red^{-1}\, \sum_{  { r, J\subset I_r}}  \varphi_r^{\overline I}(\nabla (L^{\overline I},L^{\overline I
\cap   J},k^{\overline I})).$$ 
Here the first equality is   the definition of  $\tau (M, e_k^M, \omega_L^M)$. 
The second equality follows from Theorem 9.2. The third
equality is a tautology (the projections of both sides to each $F_r$ are equal). The fourth equality follows from the  fact that the
terms  on the right-hand side corresponding to  
 proper subsets $J\subset  I_r$ are equal to $0$. Indeed for $i\in I_r\backslash J$,  the projection of $[t_i] -1$ and hence of  $ (
[t_i] -1)_\red^{-1}$ to $F_r$ is equal to $0$. Therefore $( [t_i] -1)_\red^{-1} \, \Im (\varphi_r^{\overline I})=0$.
The theorem now follows from the next claim.
\skipaline

{\bf Claim.} For any sets $ I\subset  J \subset I_0 $  and any $ a\in \bold Z[H(L^{\overline I}, L^{\overline I \cap J} )]$, 
$$\sum_{  { r, J\subset I_r}}  \varphi_r^{\overline I}(a)=a^\tr.  $$

\skipaline {\sl Proof.} By additivity, it suffices to consider the case  where $a\in  H(L^{\overline I}, L^{\overline I \cap J} )$.  Let
$H=H_1(M)$ and $ H_J, p,q$ be the same objects as in Sect.\ 11.1.  Since $\Ker\, p$ is a   subgroup
of
$\Tors\, H$, the  epimorphism $p$ extends to a ring epimorphism $\tilde p:Q(H) \to  Q(H_J)$. Moreover, the  splitting of
$Q(H_J)
$ into a direct sum of fields is obtained from the splitting of $Q(H)$ by quotienting out all $F_r$ such that $\varphi_r(\Ker \,p)
\neq 1$.   Note  that  $\varphi_r(\Ker \,p)=1$ if and only if   $J\subset I_r$. 
Thus
$$Q(H_J) =\bigoplus_{r,  \varphi_r(\Ker \,p)=1} F_r=\bigoplus_{r,  J\subset I_r} F_r.$$
For each $r$ such that $J\subset I_r$,  the composition of $\tilde p:Q(H) \to  Q(H_J)$ with the projection $\psi_r:Q(H_J)\to
F_r$ is  equal to the projection $Q(H)\to F_r$.

  The identity
$\sum_r 
\varphi_r (h)=h$ for all
$h\in H$ implies that $$ a^\tr=  \vert
\Ker\, p
\vert^{-1}\sum_{h\in p^{-1}  q (a)} h =
\vert
\Ker\, p
\vert^{-1}
\sum_r  \varphi_r (\sum_{h\in p^{-1}  q (a)} h).$$ Observe  that if $\varphi_r (\Ker\, p)\neq 1$ then  $\varphi_r (\sum_{h\in
p^{-1}  q (a)} h)=0$. Therefore
$$ a^\tr=    \vert \Ker\, p \vert^{-1} \sum_{r,  \varphi_r(\Ker \,p)=1}  \varphi_r (\sum_{h\in p^{-1}  q
(a)} h)
$$
$$ = \vert \Ker\, p \vert^{-1} \,  \sum_{r,   J\subset I_r}  \psi_r  \tilde p (\sum_{h\in p^{-1}  q
(a)} h) = \sum_{r,   J\subset I_r} \psi_r(q(a)) =\sum_{   r,  J\subset I_r}  \varphi_r^{\overline I}(a)$$
where the last formula follows from the definition of $\varphi_r^{\overline I}=(\varphi_r)^{\overline I}$.

The proof of (11.2.b) is similar  with a separate  analysis of
the term corresponding to $I=J=I_0$  using  Lemma  10.3.  

\skipaline \noindent {\bf    11.3.  Remarks. } 1. It  follows from Lemma 10.3 that  $  {\det (\ell^{ I_0}(L))}/{ \vert
\Tors H\vert }=\pm 1$ in (11.2.b).  This sign will also appear in   Appendix 2.

2.  As an exercise, the reader may deduce  Theorem 10.2 from Theorem 11.2. Hint: use the fact that 
$\Delta  (M,e_k^M, \omega_L^M)$ is the image of $\tau(M,e_k^M, \omega_L^M)$ under the projection $Q(H)\to Q(H/\Tors\,
H)$ and that all summands corresponding to $J\neq I_0$ are annihilated by this projection.

\skipaline \noindent {\bf    11.4.  Special cases. } We discuss two    special cases of Theorem 11.2. 
We  use the notation of Sect.\ 11.1.  For a  Laurent polynomial $\nabla\in  \bold Z [t_1^{\pm 1},...,
t_m^{\pm 1}]$, we denote by $[\nabla]$ 
 the element of $\bold Z [H]$  obtained from $\nabla  $ by the substitution  $t_i\mapsto [t_i] \in H$ where 
$i=1,...,m$.
Here is the  simplest case of (11.2.a): if  all
$[t_1],..., [t_m]$   have infinite order in
$H$ then $I_0=\emptyset$ and (11.2.a) yields 
$$\tau(M,e_k^M, \omega_L^M)=  \prod_{i=1}^m
([t_i]-1)^{-1} \,  [\nabla (L,k)] .
\eqno (11.4.a)$$

Consider now the case where   $L$ is algebraically split. Let $f_i\in \bold Z$ be the framing number of $L_i$ for $i=1,..., m$.  
For $i\in I_0=\{i=1,...,m\,\vert \, f_i\neq 0\}$,  the element $[t_i]\in H$ has order $\vert f_i\vert=\sign (f_i) f_i$  where 
$\sign (f_i)=\pm 1$ is the sign of
$f_i$. Set
$$s_i=\sign (f_i) (1+ [t_i] +[t_i]^2+...+ [t_i]^{\vert f_i\vert-1}) \in \bold Z[H].$$
Assume first that $b_1(M)=m- \vert I_0 \vert \geq 2$. As we show below, formula (11.2.a) can be rewritten in this case as
follows:
$$\tau(M,e_k^M, \omega_L^M)= \sum_{I\subset I_0}\,\, (-1)^{\vert I \vert} \,  ( \prod_{i\in I}  s_i  )\,  [ \check \nabla
(L^{\overline I} ,k^{\overline I})] . \eqno (11.4.b)$$
Here
$k^{\overline I}$ is  just the restriction of $k$ to  $L^{\overline I}=\cup_{i\in \overline I} L_i$.  The
assumption $b_1(M)\geq 2$  ensures that $L^{\overline I}$ has $\geq 2$ components  for any $I\subset I_0$  so that 
$\check \nabla
(L^{\overline I} ,k^{\overline I})$ is a Laurent polynomial on $\{t_i\}_{i\in \overline I}$ and 
$[ \check \nabla
(L^{\overline I} ,k^{\overline I})]$ is a well defined element of $\bold Z[H]$. 
The sum on the right-hand side of (11.4.b)    contains a term corresponding to $I=\emptyset$
and equal to $ [ \check \nabla
(L  ,k  )]$. The other terms correspond to  proper sublinks of $L$. For instance, if  $f_i=0$ for all $i$ then $I_0=\emptyset$
and   we obtain a formula $\tau(M,e_k^M, \omega_L^M)=
[ \check \nabla
(L  ,k  )]$ equivalent to 
(11.4.a).  If $I_0$ has  only  one element $i$ then 
(11.4.b) gives 
$$\tau(M,e_k^M, \omega_L^M)=
 [ \check \nabla
(L  ,k  )]-   s_{i} \, [ \check \nabla
(L^{\overline {i}} ,k^{\overline  {i}})]$$
where $\overline {i} = \{1,...,m\} \backslash {i}$. 
If $I_0$ has  two elements $i, j$ then 
(11.4.b) gives 
$$\tau(M,e_k^M, \omega_L^M) =
 [ \check \nabla
(L  ,k  )]-   s_{i} \, [ \check \nabla
(L^{\overline {i}} ,k^{\overline  {i}})] -   s_{j} \, [ \check \nabla
(L^{\overline {j}} ,k^{\overline  {j}})]  +    s_{i} s_{j} \, [ \check \nabla
(L^{\overline {\{i, j\}}} ,k^{\overline  {\{i, j\}}})].$$

Let us   deduce (11.4.b) from (11.2.a). Observe that for  $I\subset J\subset  I_0$,
$$\det (\ell^{ I}(L))=\prod_{i\in I} f_i.
$$
  It follows from  
definitions that
 $\nabla (L^{\overline I},L^{\overline I
\cap   J},k^{\overline I})= \prod_{i\in \overline J} (t_i-1)\,\,  \check \nabla
(L^{\overline I}, k^{\overline I})$.
Hence 
$$(\nabla (L^{\overline I},L^{\overline I
\cap   J},k^{\overline I}))^\tr = \prod_{j\in J} \sigma_j \, \prod_{i\in \overline J} ([t_i]-1)\,\,  [\check \nabla
(L^{\overline I}, k^{\overline I})]$$
where
$$\sigma_j= \vert f_j \vert^{-1} (1+[t_j]+...+[t_j]^{\vert f_j\vert -1})=   f_j^{-1} s_j. \eqno (11.4.c) $$
Formula (11.2.a) yields then
$$ \tau (M, e_k^M, \omega_L^M)\eqno (11.4.d)  $$
$$=
\sum_{  { I\subset J\subset  I_0}} (-1)^{\vert I \vert}  \,\prod_{i\in I} f_i  \, \prod_{i\in \overline {J}}
([t_i] -1)_\red^{-1} \,\prod_{j\in J} \sigma_j
 \, \prod_{i\in \overline J} ([t_i]-1)\,\,  [\check \nabla
(L^{\overline I}, k^{\overline I})]. 
 $$
Observe   that for  $i\in \overline I_0$, we have  $([t_i] -1)_\red^{-1}
  ([t_i]-1)=([t_i] -1)^{-1}
  ([t_i]-1)=1$. If $i\in I_0$ then $([t_i] -1)_\red^{-1}
  ([t_i]-1)=1-\sigma_i$:  both sides  are mapped to $0$ by any character of $\Tors H$ mapping $t_i$ to $1$ and are 
 mapped to $1$ by all other characters of $\Tors H$. Therefore  the right-hand side of (11.4.d) can be rewritten as
$$\sum_{  { I\subset J\subset  I_0}} (-1)^{\vert I \vert}  \,\prod_{i\in I} f_i  \,\prod_{j\in J} \sigma_j \, \prod_{i\in I_0\backslash
{J}} (1-\sigma_i)\,  [\check \nabla
(L^{\overline I}, k^{\overline I})]  $$
$$=\sum_{  { I\subset    I_0}} (-1)^{\vert I \vert} \, \prod_{i\in I}
f_i
\,  \big (  \sum_{  {J,  I\subset J\subset    I_0}} \prod_{j\in J} \sigma_j \,
\prod_{i\in I_0\backslash {J}} (1-\sigma_i)\big ) \,  [\check \nabla
(L^{\overline I}, k^{\overline I})].$$
The sum over $J$ obviously equals to $\prod_{i\in I} \sigma_i$. Therefore 
$$ \tau (M, e_k^M, \omega_L^M) =\sum_{  { I\subset    I_0}} (-1)^{\vert I \vert} \, \prod_{i\in I}
f_i
\,   \prod_{i\in I} \sigma_i\,  [\check \nabla
(L^{\overline I}, k^{\overline I})] .$$
Substituting  
$\sigma_i=    f_i^{-1} s_i$ we obtain (11.4.b).

 Assume now that $L$ is algebraically split and  $b_1(M)=m- \vert I_0\vert =1$. Let $n$ be the only element of
$\overline {I_0}$. Then    (11.2.b) can be similarly rewritten   as follows:
$$ \tau (M, e_k^M, \omega_L^M) =
\sum_{I\subset I_0, I\neq I_0}\,\, (-1)^{\vert I \vert}\, \prod_{i\in I}  s_i \,\, [ \check \nabla
(L^{\overline I} ,k^{\overline I})] \eqno (11.4.e)$$
$$
+
  (-1)^{   m} \,\big (\prod_{i\in I_0}  s_i \big ) \,  
[t_{n}]^{(k_n+1)/2}\,([t_{n}] -1)^{-2}\, 
\Delta_{L_n} ([t_n]) .
 $$

\skipaline \skipaline \centerline  {\bf 12.  A surgery formula for the Seiberg-Witten invariants  of 3-manifolds}

\skipaline \noindent {\bf    12.1.  Torsion versus SW-invariants. }  Let $M$ be    a  closed, connected,
  oriented   3-manifold  with  homology orientation $\omega$. 
Set $H=H_1(M)$.  Assume first that   $b_1(M)\geq 2$. The
Seiberg-Witten invariant  of $M$ is a $\bold Z$-valued function $e\mapsto SW(e,\omega)$ on $\vect (M)$, see for
instance   [MT]. This function has a finite support. Note that  usually one considers the SW-invariants of $Spin^c$-structures on
$M$;   by  [Tu5], the set of 
 $Spin^c$-structures on $M$ can be   identified with $\vect (M)=\Eul (M)$. The SW-invariant is computed
from the torsion $\tau $ as follows.  Every element $a\in \bold Z[H]$ expands uniquely as a finite sum $\sum_{h\in H} a_h
h$ with $a_h\in \bold Z$. By [Tu6] (see also Appenidx 3), 
$\tau(M,e,\omega)\in \bold Z[H] $ for any 
$e\in
\vect (M)$
and
$$SW(  e, \omega)=\pm (\tau(M,e,\omega))_1\in  \bold Z  \eqno (12.1.a)$$ where $1$ is the neutral element of $H$
 and the sign $\pm$ does not depend
on
$e$. 
Together  with  
  (11.2.a)  this  implies    the following   formula   for $\pm SW(  e_k^M, \omega_L^M)$  in the setting of Sect.\
11.1:
$$ \pm SW(  e_k^M, \omega_L^M) =
\sum_{  { I\subset J\subset  I_0}} (-1)^{\vert I \vert}  \det (\ell^{ I}(L))  \left ( \prod_{i\in \overline {J}} ([t_i] -1)_\red^{-1}
 \,  (\nabla (L^{\overline I},L^{\overline I
\cap   J},k^{\overline I}))^\tr\right )_1  
 $$ where the sum runs over all  subsets $I,J$ of $I_0$ such that $I\subset J$.
This formula  computes   $\pm SW(  e_k^M, \omega_L^M)$  in terms of the  Conway polynomials of $L$ and its 
sublinks, and the linking and framing numbers of the components of
$L$.

Assume now that $b_1(M)=1$.  Choose an element $t\in H$  whose projection to  $H/\Tors\, H=\bold Z$ is  a  
generator.  
 An element $g\in H$ is said to be  $t$-positive  (resp. $t$-negative) if $g\in t^K \Tors\, H$ with $K>0$ (resp. with
$K<0$).  The Seiberg-Witten
invariant  of $M$ is a $\bold Z$-valued function $e\mapsto SW(e,\omega,t)$ on $\vect (M)$ depending on    $t (\mod  \Tors\,
H) \in H/\Tors\, H$. It  is computed from   $\tau(M,e,\omega)$ as follows.  Let $Q$ be the subring of $Q(H)$ generated by
$\bold Q[H]$ and elements $(g-1)^{-1}$ where $g$ runs over elements of $H$ of infinite order. 
Using the formula  $(g-1)^{-1}=-1-g-g^2-...$ for $t$-positive $g$ and  
the formula  $(g-1)^{-1}=-g^{-1} (g^{-1}-1)^{-1}= g^{-1}+g^{-2}+ ...$ for $t$-negative  $g$ 
we can uniquely expand any $a\in Q$ as a formal series $ 
a^t= \sum_{h\in H} a^t_h h$  with $a^t_h\in \bold Q$. The support of the map  
$h\mapsto a^t_h:H\to \bold Q$ is  
  essentially
$t$-positive   in the sense that it  meets the set of $t$-negative elements of $H$ in a finite  set. Theorem 11.2 or the results of 
 [Tu5, Sect.\ 4.2]  show that $\tau(M,e,\omega)\in Q$.   It follows from  [Tu6] (cf. Appendix 3 of the present paper) that
$$ \pm SW(e,\omega,t)=  (\tau(M,e,\omega))^{t^{-1}}_1= (\tau(M,e^{-1},\omega))^t_1.    \eqno (12.1.b)$$ 
Recall that $(e_k^M)^{-1}=e_{2-k}^M$, cf. Sect.\  1.7.  
Combining this with Theorem 11.2 we obtain 
 the following   formulas   for $\pm SW(  e_k^M, \omega_L^M, t)$  in the setting of Sect.\
11.1. If  
$m- \vert I_0\vert \geq 2$ (and $b_1(M)=1$),  then  $$  \pm SW(  e_k^M, \omega_L^M, t)  $$
$$=
\sum_{  { I\subset J\subset  I_0}} (-1)^{\vert I \vert} \,\det (\ell^{ I}(L))\left ( \prod_{i\in \overline {J}} ([t_i] -1)_\red^{-1}
 \,  (\nabla (L^{\overline I},L^{\overline I
\cap   J},k^{\overline I}))^\tr\right )^{t^{-1}}_1
 $$
$$=
\sum_{  { I\subset J\subset  I_0}} (-1)^{\vert I \vert} \,\det (\ell^{ I}(L))\left ( \prod_{i\in \overline {J}} ([t_i] -1)_\red^{-1}
 \,  (\nabla (L^{\overline I},L^{\overline I
\cap   J},(2-k)^{\overline I}))^\tr\right )^{t}_1. $$
Note that  $(2-k)^{\overline I} =2+(-k^{\overline I})$. 
   If  $I_0$ consists of all $i=1,...,m$ except a certain $n$, then   $t=[t_n]\in H$  projects to a generator of $H/\Tors\, H$.
Let  $\Delta_{L_n} (t)= \sum_l z_{l} t^{l} $ be the canonically normalized Alexander polynomial of
$L_n$.   Then 
$$  \pm SW(  e_k^M, \omega_L^M, t)    
 $$
$$=
\sum_{  { I\subset J\subset  I_0}, I\neq I_0  } (-1)^{\vert I \vert} \,\det (\ell^{ I}(L))\left ( \prod_{i\in \overline {J}} ([t_i]
-1)_\red^{-1}
 \,  (\nabla (L^{\overline I},L^{\overline I
\cap   J}, k^{\overline I}))^\tr\right )^{t^{-1}}_1
$$
$$+
  (-1)^{   m} \, \frac {\det (\ell^{ I_0}(L))}{ \vert \Tors H\vert }     (z_{(k_n-3)/2 }+2 z_{(k_n-5)/2} +3
z_{ (k_n-7)/2}+...)
 $$ 
$$=
\sum_{  { I\subset J\subset  I_0}, I\neq I_0  } (-1)^{\vert I \vert} \,\det (\ell^{ I}(L))\left ( \prod_{i\in \overline {J}} ([t_i]
-1)_\red^{-1}
 \,  (\nabla (L^{\overline I},L^{\overline I
\cap   J}, (2-k)^{\overline I}))^\tr\right )^{t}_1
$$
$$+
  (-1)^{   m} \, \frac {\det (\ell^{ I_0}(L))}{ \vert \Tors H\vert }     (z_{(k_n-3)/2 }+2 z_{(k_n-5)/2} +3
z_{ (k_n-7)/2}+...).
 $$ 
 
\skipaline \noindent {\bf    12.2.  The case of algebraically split links. }  The surgery formula  for $SW(  e_k^M,
\omega_L^M)$ can be made quite explicit when the link $L $ is algebraically split. Observe first that for an
algebraically split link with $m\geq 2$ components $L=L_1\cup ...\cup L_m$  in a homology sphere, the Laurent polynomial 
$\nabla_L $ is   divisible by $ \prod_{i=1}^m (t_i^2-1)$ in $\bold Z [t_1^{\pm 1},..., t_m^{\pm 1}]$. 
Thus we have a finite expansion
$$ \nabla_L (t_1,...,t_m) / \prod_{i=1}^m (t_i^2-1)= \sum_{l=(l_1,...,l_m)\in \bold Z^m} z_l (L) \,\,  t_1^{l_1}... t_m^{l_m}$$
where $z_l (L)\in \bold Z$.  

Assume   that  a 3-manifold 
$M$ is   obtained by   surgery on a    framed oriented algebraically split  link $L=L_1\cup ...\cup L_m$ in an oriented 
 3-dimensional integral homology sphere. 
 Let   $k=(k_1,...,k_m)$ be a charge on $L$. 
Let $f=(f_1,...,f_m)$ be the tuple of  the   framing numbers of $L_1,...,L_m$. 
Denote by
$J_0$  the subset of $\{1,...,m\}$ consisting of all 
$j $ such that  
$f_j=0$.  Note that $ \vert J_0 \vert=b_1(M)$. As usual we distinguish two cases $b_1(M)\geq 2$ and $b_1(M)=1$.

\skipaline \noindent {\bf    12.2.1. Case $b_1(M)\geq 2$.  } 
We claim that  $$   \pm  SW(e_k^M, \omega_L^M)  \eqno (12.2.a)
$$
$$  = \sum_{J_0\subset J \subset \{1,...,m\}}\,\, (-1)^{\vert J \vert }\,
\prod_{i\in \overline J}  \sign (f_i)  \,\,  
  \sum_{l\in \bold Z^{J}, l\equiv -k (\mod 2f)} z_{l} (L^{J}).  $$ 
 To see this, we analyze the terms  of (11.4.b) giving rise to entries of the  neutral element in $\tau(M, e_k^M, \omega_L^M)$. 
By definition (see Sect.\ 8.1 and 10.4),
 $$  \check \nabla
(L  ,k )
= \nabla (L ,k )/\prod_{i=1}^m (t_i-1)=-t_1^{k_1/2}... t_m^{k_m/2}\, \nabla_L (t_1^{1/2},...,t_m^{1/2})
/\prod_{i=1}^m (t_i-1)$$
$$=-  \sum_{l=(l_1,...,l_m)\in \bold Z^m} z_l (L) \,\,  t_1^{(k_1+l_1)/2}...
t_m^{(k_m+l_m)/2}\in \bold Z[t_1^{\pm 1}, ...,  t_m^{\pm 1}].$$ Similar formulas hold for $\check \nabla
(L^{\overline I} ,k^{\overline I})$ where $I$ runs over subsets of $I_0=\overline {J_0}$.  We substitute these formulas into
(11.4.b)  and use the obvious fact that a monomial $t_1^{n_1}... t_m^{n_m}$ represents   $1 \in H_1(M)$ 
iff $n_i \in \bold Z$ is divisible by
$f_i$ for all
$i=1,...,m$. This implies  that $1 \in H_1(M)$    appears in $\tau(M,  e_k^M,
\omega_L^M)$ with coefficient
$$- \sum_{I\subset I_0}\,\, (-1)^{\vert I \vert}\, \prod_{i\in I}  \sign (f_i) \,\, \sum_{l\in \bold
Z^{\overline I}, l\equiv -k (\mod 2f)} z_{l} (L^{\overline I}).$$
Setting  here $I=\overline J$ we obtain  
$$(-1)^{m-1} \sum_{J_0\subset J \subset \{1,...,m\}}\,\, (-1)^{\vert J  \vert}\,
\prod_{i\in \overline J}  \sign (f_i)  \,\,  
  \sum_{l\in \bold Z^{J}, l\equiv -k (\mod 2f)} z_{l} (L^{J}).$$
This and the results of  [Tu6] mentioned in Sect.\ 12.1 imply (12.2.a).

The  right-hand side of (12.2.a)     contains  two distinguished terms    corresponding to $J= \{1,...,m\}$ and to $J=J_0$. 
The term corresponding to $J= \{1,...,m\}$ 
 equals  $(-1)^m \sum_{l\in \bold Z^{m}, l\equiv -k (\mod 2f)} z_{l} (L)$. The term corresponding to $J= J_0$ 
  equals  $$(-1)^{ b_1(M)}  ( \prod_{j\in \overline J_0}  \sign (f_j)   ) \,  z_{-k^{J_0}} (L^{J_0})$$  where
$-k^{J}=\{-k_j\}_{j\in J}$ for any set $J\subset \{1,...,m\}$.   For  example, if all $f_j$ are equal to zero,
  then  $J_0= \{1,...,m\}$ and    $$SW(e_k^M, \omega_L^M)= \pm  \,
  z_{-k} (L). \eqno (12.2.b) $$  If $J_0$ consists of all indices $1,...,m$ except a certain $i$, i.e.,  
 $J_0=\overline {i} = \{1,...,m\} \backslash {i}$ 
then  (12.2.a)  gives 
$$   \pm  \, SW(e_k^M, \omega_L^M)=
 \sum_{l\in \bold Z^{m}, l\equiv -k (\mod 2f)} z_{l} (L)  -    \sign (f_{i})\,
z_{-k^{\overline {i}}} (L^{\overline {i}}).\eqno (12.2.c)$$  
The condition $  l\equiv -k (\mod 2f)$ means here  that $l_j=-k_j$ for $j\neq i$ and $l_i\equiv -k_i (\mod 2f_i)$. 
If $J_0=\overline  {\{i, j\}}$ is the complement of  two   indices $i, j$ then 
(12.2.a)  gives 
$$  \pm  \,  SW(e_k^M, \omega_L^M) =
 \sum_{l\in \bold Z^{m}, l\equiv -k (\mod 2 f)} z_{l} (L)  -   \sign (f_{i}) \,    \sum_{l\in \bold Z^{\overline {i}}, l\equiv -k (\mod 2 f)}
z_{l} (L^{\overline {i}}) $$
$$  -   \ \sign (f_{j}) \,    \sum_{l\in \bold Z^{\overline {j}}, l\equiv -k (\mod 2 f)}
z_{l} (L^{\overline {j}})
  +  \ \sign (f_{i}) \,\sign (f_{j})\,  z_{-k^{\overline {\{i, j\}}}} (L^{\overline {\{i, j\}}}) .$$

\skipaline \noindent {\bf    12.2.2. Case $b_1(M)=1$.  }   Let $n$ be the only element of
$  {J_0}$. Let $\Delta_{L_n} (t)=\sum_l z_l t^l$
be the normalized Alexander polynomial of $L_n$. Then it follows similarly from (11.4.e)
that 
$$  \pm   SW(e_k^M, \omega_L^M, [t_n])  $$
$$=
\sum_{n\in  J \subset \{1,...,m\}, \vert J\vert \geq 2}\,\, (-1)^{\vert J  \vert}\,
\prod_{i\in \overline J}  \sign (f_i)  \,\,  
  \sum_{l\in \bold Z^{J}, l\equiv 2-k (\mod 2f)} z_{l} (L^{J})
$$
$$ - \prod_{i\neq n}  \sign (f_i)   \, 
(z_{(k_n-3)/2 }+2 z_{(k_n-5)/2} +3
z_{ (k_n-7)/2}+...).$$

\skipaline \noindent {\bf    12.3.  Examples. } 1. Let a 3-manifold $M$ be obtained by surgery along an oriented knot
$L=L_1\subset S^3$ with framing $0$.
Then for any odd integer $k$, $$SW(e_k^M, \omega_L^M, [t_1])=\pm (z_{(k_n-3)/2 }+2 z_{(k_n-5)/2} +3
z_{ (k_n-7)/2}+...)$$ where $\Delta_{L_n} (t)=\sum_l z_l t^l$.

 2. Consider   the Borromean
  link  $L=L_1\cup L_2\cup L_3$. It is algebraically split, 
$z_l (L)=0$ for all $l\neq 0$ and $z_0(L)=\pm 1$ where the   sign depends on the orientation of $L$. 
All 2-component sublinks of $L$ are trivial and have a zero Alexander-Conway polynomial. Therefore  if $M $ is obtained by surgery
on $L$ with framing numbers $(f\in \bold Z,0,0)$ then  formulas  (12.2.b), (12.2.c)   show that there is  an Euler structure on
$M$ with  SW-invariant
$\pm 1$ and all other Euler structures on $M$ have a zero SW-invariant. This  was previously known for $f=0$ where 
$M=S^1\times S^1\times S^1$.
 
\skipaline \skipaline \centerline  {\bf  Appendix 1.  A surgery formula for  rational homology spheres}

\skipaline  
 We discuss here an analogue of Theorem 11.2  in the case where  $M$ is  a  rational
homology sphere, i.e., $b_1(M)=0$.  We use the notation of Sect.\ 11.1.  It is clear that $I_0=\{1,...,m\}$ and $Q(H)=\bold Q[H]$  
where $H=H_1(M)$.
For every $i=1,...,m$, we define  
  $\sigma_i\in \bold Q[H]$ by   (11.4.c) where $f_i\geq 1$ is the order of
$[t_i]\in H $.  For sets
$I\subset J\subset \{1,...,m\} $, set  $rk(I,J)= \rank \, H(L^{\overline I},L^{\overline I
\cap   J})$. 

The arguments of Sect.\ 11 apply here   with the following changes. One of the summands in the splitting of $\bold Q[H]$ into  a
direct sum of fields  corresponds to the augmentation homomorphism $\varphi_0:\bold Q[H] \to \bold Q$ which maps $H$ to
$1$. Clearly, $H_*^{\varphi_0}(M)=H_*(M;\bold Q)\neq 0$ so that   $\tau^{\varphi_0}(M)=0$. This allows us to proceed as in
the proof of Theorem 11.2 involving only $r$ such that $\varphi_r(H)\neq 1$. This leads to the condition 
$J\neq \{1,...,m\}$ or equivalently $
\vert J \vert \leq m-1$ in the
formulas below.

The second subtlety comes from the fact that the expression
$ (\nabla (L^{\overline I},L^{\overline I
\cap   J},k^{\overline I}))^\tr$   is defined only if  $rk(I,J) \geq 2$. Note that $rk(I,J)\geq  \vert 
  \overline  J \vert \geq 1$. We can have 
$rk(I,J)= 1$  only when
$J=\overline n=\{1,...,m\}\backslash \{n\}$ for a certain $n=1,...,m$. 
For  a  set  $I\subset J=\overline n$   three cases may occur:  
(i) $rk(I,J)\geq 2$; (ii)   $ rk(I,J)=1$ and $I\neq J$, and (iii) $I=J$. 
In the second case
$$([t_n]-1) \nabla (L^{\overline I},L^{\overline I
\cap   J},k^{\overline I}) \in \bold Z[ H (L^{\overline I},L^{\overline I
\cap   J} )]$$ and we can set
$$(\nabla (L^{\overline I},L^{\overline I
\cap   J},k^{\overline I}))^\tr= ([t_n] -1)_\red^{-1}\, (([t_n]-1) \nabla (L^{\overline I},L^{\overline I
\cap   J},k^{\overline I}))^{\tr} \in \bold Q[H].$$ Then the same arguments as in the proof of Theorem 11.2 show that 
$$ \tau (M, e_k^M, \omega_L^M)=  $$
$$ 
\sum_{  { I\subset J\subset \{1,...,m\}, \vert I \vert \leq m-2, \vert J \vert \leq m-1} } (-1)^{\vert I \vert} \,\det (\ell^{ I}(L))\,
\prod_{i\in
\overline {J}} ([t_i] -1)_\red^{-1}
 \,  (\nabla (L^{\overline I},L^{\overline I
\cap   J},k^{\overline I}))^\tr
$$
$$+
  (-1)^{   m} \,  \sum_{n=1}^m \,  \det (\ell^{\overline n}(L)) \,( \prod_{i\neq  n} \sigma_i)\, [t_{n}]^{(k_n- lk(L_n,
L^{\overline n})+1)/2}\, 
  (([t_n]
-1)_\red^{-1})^{2}     \,\Delta_{L_n} ([t_n]).  
 $$

In the case where $L$ is algebraically split and its components have non-zero framings $f_1,...,f_m$,   we obtain 
$$ \tau (M, e_k^M, \omega_L^M) =
\sum_{I\subset \{1,...,m\}, \vert I \vert \leq m-2}\,\, (-1)^{\vert I \vert}\, \prod_{i\in I}  s_i \,\, [ \check \nabla
(L^{\overline I} ,k^{\overline I})]  $$
$$
+
(-1)^{   m} \,  \sum_{n=1}^m \, \big (\prod_{i\neq n}  f_i\sigma_i \big ) \,  
[t_{n}]^{(k_n+1)/2}\, (([t_n]
-1)_\red^{-1})^{2} \, 
\Delta_{L_n} ([t_n]).   $$

\skipaline \skipaline \centerline  {\bf  Appendix 2.  Computation of $\omega_L^M$}

\skipaline \noindent {\bf     Canonical homology orientation. }  
Every closed oriented 3-manifold $M$ has a canonical  homology orientation $\omega_M$
determined by an arbitrary basis in the vector space $H_0(M;\bold R) \oplus H_1(M;\bold R)$ followed by the Poincar\'e dual basis
in 
$H_2(M;\bold R) \oplus H_3(M;\bold R)$.  Clearly,  $\omega_{-M}= (-1)^{b_0(M)+b_1(M)} \omega_M$.
The aim of this appendix  is to compute  the homology orientation $\omega_L^M=(\omega_L)^M$ of $M$
which appears above via  $\omega_M$.   It is obvious that $\omega_L^M=\pm  \omega_M$ and we shall   compute the sign
$\pm$ in this formula.  We begin with a   purely algebraic definition.

\skipaline \noindent {\bf      Sign of  the  determinant   as a torsion. }  Let $B$ be a symmetric $(m\times m)$-matrix
over $\bold R$ with $m=1,2,..$. Let
$b$ be the symmetric bilinear form on ${\bold R}^m$ determined by $B$.  Quotienting
${\bold R}^m$ by the annihilator 
$\Ann (b) $  we obtain a  non-degenerate symmetric bilinear form on  ${\bold R}^m/\Ann (b)$.  It  can be
represented  with respect to a   basis  of ${\bold R}^m/\Ann (b)$ by a  non-degenerate symmetric  matrix. The sign
of   its determinant  does not depend on the choice of the basis. This sign is   denoted by
$\det_0(B)$. If $B$ is non-degenerate then $\det_0(B)=\sign (\det  (B))$. 

We can interpret  $\det_0(B)$ in terms of torsions as follows.
Consider the   sequence of vector spaces and linear homomorphisms
$$
\CD
C= ( \Ann(b) \hookrightarrow {\bold R}^m   @>\ad (b)>> ({\bold R}^m)^* \to (\Ann(b))^*)
\endCD $$
where   $\ad (b): {\bold R}^m\to ({\bold R}^m)^*=\Hom ({\bold R}^m, \bold R)$ is  adjoint to $b$ and the homomorphism
$({\bold R}^m)^* \to (\Ann(b))^*$ is obtained by restricting linear functionals on ${\bold R}^m$ to $\Ann (b)$. 
Clearly,  $C$ is an acyclic chain complex.  We provide $\Ann (b), {\bold R}^m$ with arbitrary bases
and provide $  (\Ann(b))^*,({\bold R}^m)^*$ with dual bases.  It follows from definitions that the sign 
$\tau_0(C)=\pm 1$ of the  corresponding
torsion 
$\tau(C)\in \bold R$ depends only on $b$ and does not depend on the choice of   bases in $\Ann (b), {\bold R}^m$.
 Choosing these  bases
so that   the one in ${\bold R}^m$ extends the one in $\Ann (b)$ we  easily compute that
$$\tau_0(C)=  (-1)^{(m+1) \dim (\Ann(b))}\, \det_0 (B). \eqno (2.a)$$

\skipaline \noindent {\bf      Lemma. } 
 {\sl  Let $M$ be a  closed   3-manifold
obtained by  surgery on a framed oriented link $L=L_1\cup ...\cup L_m$ in an oriented 
 3-dimensional integral homology sphere  $N$. 
Let $\omega_M$ be the canonical homology orientation of $M$ determined by the orientation of $M$
  induced by the one in $N$.
Let $B$ be the  symmetric linking matrix $[ lk (L_i,L_j)]_{i,j=1}^m$
where    $
lk (L_i,L_i)$ is the framing number of $L_i$.
Then     }
$$\omega_L^M= (-1)^{b_1(M)+m+1 }\,\det_0(B) \,\omega_M. \eqno (2.b) $$

\skipaline {\sl Proof.}    We shall use the symbols $E,U_i, \Z_i$ introduced in  Sect.\ 1.7.  Let $D_i$ be  the meridional disc of
$\Z_i$. We  orient $D_i$ so that $\partial D_i\subset \partial U_i$ is homological to $L_i$ in $U_i$. 
As in Sect.\  3.2,  we provide
$H_*(M,E;\bold R)$ with the basis
$$d_1=[D_1,\partial D_1],...,d_m= [D_m,
\partial D_m], z_1=[\Z_1,
\partial \Z_1],..., z_m=[\Z_m, \partial \Z_m].$$ 
Note
that the product orientation in $\Z_i$ used to define $z_i$  is   opposite  to the orientation induced from the one in $M$.  

We provide $H_*(E;\bold R)$ with the basis $([pt], t_1,...,t_m,
g_1,...,g_{m-1})$ as in Sect.\ 3.1. Finally, we provide $H_*(M;\bold R)$ with a basis 
$[pt], h, h^*, [M]$ where $h$ is a basis in $H_1(M;\bold R)$ and $h^*$ is the  Poincar\'e dual basis in $H_2(M;\bold R)$.
The choosen bases   determine the   orientations  $\omega_L$, $\omega_M$, $\omega_{(M,E)}$ in 
$H_*(E;\bold R)$, $H_*(M;\bold R)$, 
$H_*(M,E;\bold R)$,  respectively. 
Let    $\Cal H$ be the  exact homological sequence of the pair $(M,E)$ as in Sect.\ 3.2.  Consider  
its torsion with respect to the choosen bases  
$ \tau(\Cal H)\in \bold R \backslash \{0\}$.  Let $\tau_0=\sign\, \tau(\Cal H)=\pm  1$. 
In the notation of Sect.\ 3.2 we have $\tilde \omega_L=\tau_0\, \omega_M$. Then
$$\omega_L^M=  (-1)^{m b_3(M)+ (b_1(E)+1) (b_1(M)+m) }\, \tilde \omega_L
=(-1)^{(m+1) b_1(M) + m} \,\tau_0 \,\omega_M. \eqno (2.c)$$

We  compute 
$\tau_0 $. 
The inclusion homomorphism $H_2(E;\bold R) \to H_2(M;\bold R)$ is zero and therefore   $\Cal H$ splits as a
concatenation of   three  acyclic chain complexes
$$H_0(E;\bold R) \to H_0(M;\bold R),\,\,\,\, H_3(M;\bold R) \to H_3(M,E;\bold R)\to H_2(E;\bold R) , $$
$$H_2(M;\bold R) \to H_2(M,E;\bold R)\to H_1(E;\bold R) \to H_1(M;\bold R).  $$
Therefore  $\tau_0=\varepsilon_1 \varepsilon_2 \varepsilon_3$ 
where $\varepsilon_1, \varepsilon_2, \varepsilon_3$ are the signs of the torsions of these three chain
complexes, respectively.  
The inclusion isomorphism $H_0(E;\bold R) \to H_0(M;\bold R)$ is given by the unit $(1\times 1)$-matrix  
and hence $\varepsilon_1=+1$.
The inclusion homomorphism $H_3(M;\bold R) \to H_3(M,E;\bold R)$ maps $[M]$ to 
$-(z_1+...+z_m)$ and the boundary homomorphism
$H_3(M,E;\bold R)\to H_2(E;\bold R)$ maps $z_1,...,z_m$ to $-q_1,...,-q_{m-1}, q_1+...+q_{m-1}$, respectively.
Now it is easy to compute that  $\varepsilon_2=-1$.     Using the basis $d_1,...,d_m$ we can identify $ H_2(M,E;\bold
R)$ with $\bold R^m$. We identify
$H_1(E;\bold R)$ with $(\bold R^m)^*=\Hom (\bold R^m, \bold R)$ so that the basis $t_1,...,t_m$ is dual to
$d_1,..., d_m$. With respect to these bases, the boundary homomorphism $ \partial:H_2(M,E;\bold R)\to H_1(E;\bold R) $
is presented by the matrix $B=[ lk (L_i,L_j)]_{i,j=1}^m$ and therefore equals $\ad(b)$
where $b$ is the symmetric bilinear form on ${\bold R}^m$ determined by $B$. 
This allows us to identify  $H_2(M;\bold R)=\Ker\,\partial, H_1(M;\bold R)=\Coker \,\partial$ with $\Ann (b), (\Ann(b))^*$,
respectively.  Under our choice of orientations, the integral
intersection index of $D_i$ with $t_j$ is equal to $-\delta_i^j$ where $\delta_i^j$ is  the Kronecker symbol. 
Therefore  the intersection  pairing between $H_2(M;\bold R)$ and $ H_1(M;\bold R)$ 
is equal to $-1$ times the Kronecker pairing between  $\Ann (b)$ and $ (\Ann(b))^*$.
By  formula (2.a), we have $$\varepsilon_3=(-1)^{b_1(M)} (-1)^{(m+1)b_1(M)} \,\det_0(B)=(-1)^{m b_1(M)}
\,\det_0(B).$$ Hence $\tau_0 = -(-1)^{m b_1(M)} 
\det_0(B)$. Substituting this in (2.c) we obtain (2.b).

\skipaline \noindent {\bf     Corollaries. }  We can rewrite  Theorems 9.2, 10.2, 11.2   
replacing  
$\omega_L^M$ with $\omega_M$ and simultaneously inserting  on the right-hand side the factor 
$(-1)^{b_1(M)+m+1 }\,\det_0(B)$.

\skipaline \skipaline \centerline  {\bf  Appendix 3. Corrections and additions to [Tu5], [Tu6]}

\skipaline \noindent {\bf    Corrections and additions to [Tu5]. } Let
$M$ be a compact  connected oriented 3-manifold  whose boundary  is either void or consists of tori.  
 Fix  a homology orientation $\omega$ of $M$. 
In the case 
$b_1(M)=1$ we fix also an element  $t\in H=H_1(M)$  whose projection to  $H/\Tors\, H=\bold Z$ is  a  
generator. 
Following [Tu5], we define a numerical      
{\sl torsion function} $T_\omega$ on the set of Euler structures on  $M$.   In the case
$b_1(M)=1$ this function  depends  on $t (\mod \Tors\, H)\in H/\Tors\, H$ and is denoted  $T_{\omega, t}$. Let
$e\in
\vect (M)$. If
$b_1(M)\geq 2$ then in the notation of Sect.\  12.1, 
$T_\omega(e)=(\tau(M,e,\omega))_1\in \bold Z$. 
If
$b_1(M)=0$ then
$T_\omega(e)=(\tau(M,e,\omega))_1\in \bold Q$. 
If $b_1(M)= 1$ then
$T_{\omega,t}(e)=(\tau(M,e,\omega))^t_1\in \bold Z$.  Clearly,  $T_{-\omega}=-T_\omega$.

The torsion function   $T$ in [Tu5] coincides with  $T_\omega$ under the following two restrictions on
$\omega, t$.  If
$\partial M=\emptyset$ then $\omega={\omega_M}$   should be the canonical homology orientation of
$M$ (see Appendix 2). If $\partial M\neq \emptyset$ and $b_1(M)=1$ 
then  $t=t(\omega)$ should be chosen so that  $\omega$ is determined by the  basis
$([pt], t)$ in homology. 

Several inaccuracies occured in [Tu5] in  the 
case of closed $M$ with $b_1(M)=1$. Their source   is a  mistake in the duality formulas   in  [Tu5,  Sect.\ 2.7, 3.4]. In
these formulas instead of $c(e)(=e/e^{-1})$ should be $(c(e))^{-1}$. As a result,  in a few formulas   the sign should be inverted. We
give here the correct formulas.  For $e\in \vect (M)$ denote by $K_t(e)$ the unique integer $K$ such that $c(e)\in t^K \Tors\,H$. 
Set $\Sigma=\Sigma_{h\in \Tors\, H} h\in \bold Z[ H]$. In
Theorem 4.2.3 of [Tu5] we should have  
 $$\tau_t (M,e) =\tau (M,e,\omega_M) +  \frac {K_t(e)+2}{2} ( 1-t)^{-1} \Sigma-( 1-t)^{-2} \Sigma\in \bold Z[H]. \eqno (3.a)$$ 
On p. 690,  line 11 from below   should be written  $\tau_{t^{-1}}(M,e)=\tau_{t}(M,e)-(K_t(e)/2)$.
On p. 691, line 18 from above should be  written $r=(-K-2)/2$.   
On   p. 694, lines 3 and 6   should be 
$T_t(e)=q^e_t(1)-K_t(e)/2$ and $T_{t^{-1}}(e)= T_t(e)+K_t(e)/2$ (where the torsion functions
$T_{t}, T_{t^{-1}}$ correspond to 
$\omega=\omega_M$).

We state here the duality property for $T_\omega$. Let $e\in \vect (M)$.
If $b_1(M)\neq 1$ then $$T_\omega(e)=(-1)^{b_0(\partial M)}\, T_\omega(e^{-1}).  \eqno (3.b)$$
If  $b_1(M)= 1$ then $$T_{\omega,t}(e)=(-1)^{b_0(\partial M)}\, T_{\omega,t^{-1}}(e^{-1}). \eqno (3.c)$$
These equalities result  from  the   formula, proven below, 
$$\overline {\tau(M,e,\omega)} =(-1)^{b_0(\partial M)}\, \tau(M,e^{-1}, \omega)  \eqno (3.d)$$
where the overbar  denotes the involution in $Q(H)$ sending any $h\in H$ to $h^{-1}$. 
  It is clear  that (3.d) implies (3.b) for $b_1(M)\neq 1$.
Let us check (3.c). If $\partial M=\emptyset$ then replacing if necessary $\omega$ by
$-\omega$ we can assume that  $\omega=\omega_M$.  By (3.a), (3.d),
$$\tau (M, e^{-1},\omega)= \overline {\tau(M,e,\omega)}= \overline {\tau_t (M,e)}-
\frac {K_t(e)+2}{2} ( 1-t^{-1})^{-1} \Sigma+( 1-t^{-1})^{-2} \Sigma.$$
Therefore 
$$T_{\omega,t^{-1}}(e^{-1})=( \overline {\tau_t (M,e)})_1- \frac {K_t(e)+2}{2}+1=(\tau_t (M,e))_1- \frac
{K_t(e)}{2}=T_{\omega,t}(e).$$ If $\partial M\neq \emptyset $ then the proof is similar using that by [Tu5],
Theorem 4.2.1,
  $$\tau (M,e,\omega) - ( 1-t(\omega))^{-1} \Sigma \in \bold Z[ H ].\eqno (3.e)$$ 
As an exercise, the reader can prove that if $b_1(M)= 1$ and $\partial M\neq \emptyset $ then 
$T_{\omega,(t(\omega))^{-1}} (e)=T_{\omega,t(\omega)} (e)-1$.

 \skipaline \noindent {\bf    Proof of (3.d). } In the case $\partial M=\emptyset$, formula (3.d) follows from the results of 
  [Tu4], Appendix B  combined with the computation of signs in [Tu3], Appendix.  We give
here a  proof in the case where
 $b_1(M)\geq 1$ and  $\partial M$ is possibly non-void.  This covers all
possible cases since in our setting $b_1(M)\geq 1$ 
  whenever $\partial M\neq \emptyset$. 

We first establish (3.d) when $M=E$ is the exterior of a link   $L=L_1\cup ...  \cup L_m \subset S^3$. 
Let $k=(k_1,...,k_m)$ be a charge on $L$. The key fact is the equality $\overline {\nabla}_L=(-1)^m\,\nabla_L$
(see Sect.\ 8.1). We have
$$\overline {\tau(E,e_k,\omega_L)}=\overline {\nabla (L,k)}=-t_1^{-k_1/2}
... t_m^{-k_m/2} \,  \nabla_L(t_1^{-1/2},... ,t_m^{-1/2})$$
$$=
-(-1)^{m}\, t_1^{-k_1/2}
... \, t_m^{-k_m/2} \,  \nabla_L(t_1^{1/2},..., t_m^{1/2})  =(-1)^{m}\, \nabla (L,-k)$$
$$=
(-1)^{m}\, \tau(E,e_{-k},\omega_L)=(-1)^{m} \,\tau(E,(e_k)^{-1},\omega_L).$$

We claim that if (3.d) holds for a 3-manifold $E$ with $b_1(E)\geq 2$ then it holds for  a 3-manifold $M$ with $b_1(M)\geq
1$ obtained from
$E$  by gluing $m$ directed solid tori whose cores represent   elements of infinite order, $h_1,...,h_m\in H_1(M)$.  
Indeed, let $e\in \vect (E)$ and $\omega$ be a homology orientation of $E$.  
Let
      $\inc: \bold Z[ H_1(E) ] \to
\bold Z[ H_1(M) ] $ be the inclusion homomorphism.  Assume first that $m=1$ and set $h=h_1 $.  By Lemma 7.3.3,
${\tau  (M,e^M , \omega^M)}= (h-1)^{-1} \, {\inc} (\tau (E,e, \omega))$.
Hence
$$\overline {\tau(M,e^M,\omega^M)}= (h^{-1}-1)^{-1} \,\overline {  \inc (\tau (E,e, \omega))}$$
$$
=(h^{-1}-1)^{-1}  \, {\inc} (\overline {  \tau (E,e, \omega)})=
-h (h-1)^{-1} \, (-1)^{b_0(\partial E)} \, {\inc} (   \tau (E,e^{-1}, \omega) )$$
$$
 =  (-1)^{b_0(\partial M)}\, h (h-1)^{-1} \,  \inc (c(e)^{-1})  \, {\inc} (   \tau (E,e, \omega) ) $$
$$
=      (-1)^{b_0(\partial M)} h \, \inc (c(e)^{-1}) \,   {\tau  (M,e^M , \omega^M)}$$
$$
=      (-1)^{b_0(\partial M)} h\, \inc (c(e)^{-1}) \,  c(e^M) \, {\tau  (M,(e^M)^{-1} , \omega^M)}.$$
It remains   to observe that $c(e^M)= \inc (c(e)) h^{-1}$  (cf. (1.5.c)).  The case $m>1$ is   similar.

Now we can accomplish the proof of (3.d) in the case
 $b_1(M)\geq 1$. The argument given in [Tu5], Sect.\ 3.9 shows that there is a framed link $L\subset M$ (with $\geq
2$ components)  whose components represent elements of infinite order in $H_1(M)$ and such that   the exterior, $E$,  of
$L$ in $M$ is   homeomorphic to the exterior of a link in $S^3$. The arguments above imply that (3.d) holds for $E$
and therefore   for
$M$. 

\skipaline \noindent {\bf    Corrections and additions to [Tu6]. } In [Tu6], Sect.\ 1 the definition of the distinguished 
Euler structure on a directed solid torus should be the same as in the present paper. Thus,  on the last line of [Tu6], Sect.\ 1.3 
instead of
$c(s_t)=t$ should be
$c(s_t)=t^{-1}$. The excision formula   in [Tu6], Sect.\ 2.3  should look like
$$\inc_*(v(E))=\pm \prod_{i=1}^m (1-[L_i]^{-1})\, v(M).   \eqno (3.f)$$
The corresponding changes should be implemented in  the corollaries of this formula   in   [Tu6], Sect.\   3.  However, all the statements
 in [Tu6], Sect.\ 2 and 3   remain valid.

In [Tu6], Sect.\  4.2 it is  claimed that the torsion function satisfies all the axioms stated in [Tu6], Sect.\ 2.  Here we give   more
details.  Let
$\Cal S$ be the class  of triples $(M,\omega,e)$ where  $M$ is a  compact connected  oriented 
3-manifold with $b_1(M)\geq 1$   whose boundary  is either void or consists of tori,
$\omega$ is a homology orientation on $M$, and $e$ is an  Euler structure
(= relative
$Spin^c$-structure) on $M$.
 In the case 
$\partial M=\emptyset $ we assume that $\omega=\omega_M$ is induced by the orientation of $M$. In the case 
$b_1(M)=1, \partial M=\emptyset $ the manifold $M$ is assumed to be endowed with  a homology class $t(M)\in 
H=H_1(M)
$  whose projection to 
$H/\Tors\, H$ is  a   generator.
 We define a $\bold Z$-valued function $v$ on $\Cal S$ by
$$
v(M,\omega,e)=\cases T_{\omega} (e^{-1}),~ {  {if}}\,\,\, b_1(M)\geq 2, \\
 T_{\omega,t } (e^{-1}),~ {  {if}} \,\,\, b_1(M)=1
\endcases \eqno (3.g)$$
where $t=t(M)$ in the case $b_1(M)=1, \partial M=\emptyset $ and $t=t(\omega)$ in the case $b_1(M)=1, \partial
M\neq\emptyset$.  We claim that   $v$ satisfies   the axioms stated in [Tu6], Sect.\ 2. By the uniqueness,  $v$ coincides
with the  Seiberg-Witten function $SW:\Cal S\to \bold Z$ at least up to sign depending only on the underlying 3-manifold. Combining
this   with (3.b), (3.c) we obtain   (12.1.a), (12.1.b). 

Before discussing the properties of $v$ we comment on the relations between $v$ and $\tau$. Assume first that  $b_1(M)\geq
2$. For any $e\in \vect   (M)$, $g\in H$,
$$(\tau(M,e,\omega))_g=(g^{-1}\tau(M,e,\omega))_1=(\tau(M, g^{-1}e,\omega))_1$$
$$=
v(M, (g^{-1}e)^{-1},\omega)=v(M, g e^{-1},\omega).$$
Hence $\tau(M,e,\omega)=\sum_{g\in H} v(M, g e^{-1},\omega) \, g$. The formal expression $$v(M)=\sum_{e\in
\vect  (M)} v(M,e,\omega)\, e\in \bold Z[\vect  (M)]$$ can therefore be computed as follows. Fix $e_0\in \vect  (M)$. Then 
$$v(M)=\sum_{e\in
\vect  (M)} v(M,e,\omega)\,  e=\sum_{g\in
H } v(M,g e^{-1}_0,\omega) \,ge^{-1}_0= \tau(M,e_0,\omega)\, e^{-1}_0.$$
Similarly, if
$b_1(M)=1$ then $(\tau(M,e,\omega))^t=\sum_{g\in H} v(M, g e^{-1},\omega) \,g$
and  
$$v(M) =\sum_{e\in
\vect  (M)} v(M,e,\omega)\,  e= (\tau(M,e_0,\omega))^t \,e^{-1}_0\in \bold Z[[\vect  (M)]]$$ where $t $ is as in (3.g).  

The function $v$ satisfies Axioms  1 - 4 stated in [Tu6], Sect.\ 2. Axiom   1 (topological invariance) is obvious. Axiom 2 (first part)
consists in the finiteness of the support of $v$ for any $M$ with $b_1(M)\geq 2$. This follows from the inclusion 
$\tau(M,e,\omega)\in \bold Z[H_1(M)]$.  Axiom 2 (second part)
claims that  the  support of $v$ is essentially $t$-positive for any   $M$ with $b_1(M)=1$ where 
  $t $ is as in (3.g). 
This   follows from   (3.a)  and 
(3.e). Axiom 3 amounts to the excision formula  (3.f) provided $b_1(E)\geq 2$.
Here $\inc_*$ is an additive homomorphism  of
abelian groups
$\bold Z[\vect  (E)]\to  \bold Z[\vect  (M)]$ sending any $e\in \vect  (E)$  to $e^M\in \vect  (M)$ in the notation of 
Sect.\ 1.5.  We check (3.f). Assume first that  $b_1(M)\geq 2$.  Fix   a homology orientation $\omega$ in $E$ and $e_0\in \vect  (E)$.
Note that
$\pm v(E)$ does not depend on the choice of   $\omega$.  Let
$\inc$ denote the inclusion homomorphism $\bold Z[H_1(E)]\to \bold Z[H_1(M)]$. By Lemma 7.3.3,
$$\pm \inc_* (v(E))=\pm \inc_* (\tau(M,e_0,\omega) \, e^{-1}_0 )=
\pm \inc (\tau(M,e_0,\omega)) (e^{-1}_0 )^M $$
$$ =  \pm  \prod_{i=1}^m ([L_i]-1)\,{\tau  (M,e_0^M , \omega^M)}  \prod_{i=1}^m  [L_i] ^{-1} (e_0^M)^{-1} 
=\pm  \prod_{i=1}^m (1-[L_i]^{-1})\, v(M).$$
In the case $b_1(M)=1$ the proof is similar using that 
$$\prod_{i=1}^m ([L_i]-1)\,{\tau  (M,e_0^M  , \omega^M)}=
\prod_{i=1}^m ([L_i]-1)\,{(\tau  (M,e_0^M , \omega^M))^t}.\eqno (3.h)$$
Formula (3.h)  follows 
from (3.a) (resp. from (3.e)) if $  \partial M=\emptyset$ and $m\geq 2$ (resp. if $  \partial M\neq \emptyset$).
If $b_1(M)=1, \partial M=\emptyset$ and $m= 1$ then the assumption $b_1(E)\geq 2$ implies that 
$[L_1]\in \Tors \, H_1(M)$. Then $[L_1]-1$  annihilates $\Sigma$ so that (3.h) also  follows 
from (3.a). 
Finally, Axiom 4 follows from the fact that for  the exteriors of links in $S^3$, the torsion $\tau$ coincides  with  the Milnor torsion
[Mi1]  which is  equivalent to the Alexander polynomial.

\skipaline

\centerline{\bf References}

 \skipaline

[Co] J. Conway,  {\it An enumeration of knots and links, and some of their algebraic properties},  1970 Computational Problems in
Abstract Algebra (Proc. Conf., Oxford, 1967)  Pergamon, Oxford, 329--358.

[CF] R. Crowell  and  R. Fox,  Introduction to knot theory.   Graduate Texts in Math.  57. Springer-Verlag, New York-Heidelberg,
(1977). 

[FS] R. Fintushel and R.J. Stern, {\it  Knots, links, and $4$-manifolds}, Invent. Math. 134 (1998),  363--400. 

[Fr] W. Franz,   {\it  Torsionsideale, Torsionsklassen und Torsion}, J.  Reine Angew. Math.
 176 (1937), 113--124.

[Ha] R. Hartley, {\it The Conway potential function for links}, Comment. Math. Helv. 58 (1983),  365--378. 

[HL] M. Hutchings and  Y.-J. Lee, {\it   Circle-valued Morse theory, 
Reidemeister torsion,
and Seiberg-Witten invariants of 3-manifolds},  Topology 38 (1999), no. 4, 861--888.

 [Li] Y. Lim, {\sl  Seiberg-Witten invariants for $3$-manifolds in the case $b\sb 1=0$ or $1$}, Pacific J. Math. 195 (2000),  
179--204.
 
[MT] G. Meng and C.H. Taubes, $\underline {SW} =$ {\it  Milnor torsion}, Math. Research
Letters 3 (1996),  661--674.

[Mi1] J. Milnor,   {\it  A duality theorem for Reidemeister torsion}, 
Ann. of Math. 76 (1962), 137--147. 

[Mi2] J. Milnor,  {\it Whitehead torsion}, Bull. Amer. Math. Soc. 72 (1966), 358--426.

 [Mo] J.W. Morgan, The Seiberg-Witten equations and applications to the topology of
smooth four-manifolds. Math. Notes 44, Princeton Univ. Press (1996),

 [MOY]  T. Mrowka, P. Ozsv\'ath, and B. Yu, {\it Seiberg-Witten monopoles on Seifert
fibered spaces},   Comm. Anal. Geom. 5 (1997), no. 4, 685--791. 

[OT] C. Okonek, A. Teleman,   Recent developments in Seiberg-Witten theory and complex geometry. Several complex
variables (Berkeley, CA, 1995--1996), 391--428, Math. Sci. Res. Inst. Publ., 37, Cambridge Univ. Press (1999).

 [Tu1] V.G. Turaev,  {\it The Alexander polynomial of a three-dimensional manifold},
Mat.  Sb. 97:3 (1975), 341-359;  English transl.: Math. USSR Sb. 26:3 (1975), 313--329.

 [Tu2] V.G. Turaev,  {\it Reidemeister torsion and the Alexander polynomial}, Mat.
Sb. 101:2 (1976), 252--270; English transl.: Math. USSR Sb.  30:2 (1976), 221--237.

[Tu3] V.G. Turaev,  {\it Reidemeister torsion in knot theory}, Uspekhi Mat. Nauk 41:1
(1986), 97--147;  English transl.: Russian Math. Surveys 41:1 (1986), 119--182.

[Tu4] V.G. Turaev,  {\it Euler structures, nonsingular vector fields, and
torsions of Reidemeister type}, Izvestia Ac. Sci. USSR 53:3  (1989); 
English transl.: Math. USSR Izvestia 34:3 (1990), 627--662. 

[Tu5] V.G. Turaev,  {\it Torsion invariants of $Spin^c$-structures
on 3-manifolds},  Math. Research Letters 4 (1997),  679--695.

 [Tu6]  V.G. Turaev, {\it  A combinatorial formulation for the Seiberg-Witten invariants of
3-manifolds},   Math. Research Letters,  5 (1998) 583--598.

 [Tu7]  V.G. Turaev, {\it  Introduction to Combinatorial Torsions},   Birkh\"auser (2001).

\skipaline
  Institut de Recherche Math\'ematique Avanc\'ee, Univ. Louis  
Pasteur  -  C.N.R.S., 7 rue Ren\'e Descartes, F-67084 Strasbourg, France

Max-Planck Institut f\"ur Mathematik, Vivatsgasse 7, D-53111 Bonn, Germany

\end